\documentclass{amsart}
\newtheorem{thm}{Theorem}[section]

\newtheorem{lem}[thm]{Lemma}

\theoremstyle{definition}

\theoremstyle{remark}

\numberwithin{equation}{section}

\begin{document}

\title{Vlasov-Maxwell-Boltzmann diffusive limit}%
\author{Juhi Jang}
\address{Department of Mathematics,
 Brown University, Providence,
RI 02912, USA} \email{juhijang@math.brown.edu}

\date{January 27, 2006}%
\begin{abstract}
We study the diffusive expansion for solutions around Maxwellian
equilibrium and in a periodic box to the Vlasov-Maxwell-Boltzmann
system, the most fundamental model for an ensemble of charged
particles. Such an expansion yields a set of dissipative new
macroscopic PDE's, the incompressible Vlasov-Navier-Stokes-Fourier
system and its higher order corrections for describing a charged
fluid, where the self-consistent electromagnetic field is
present. The uniform estimate on the remainders is established
via a unified nonlinear energy method and it guarantees the global
in time validity of such an expansion up to any order.
\end{abstract}
\maketitle

\section{Introduction and Formulation}\

The dynamics of charged dilute particles can be described by the
celebrated Vlasov-Maxwell-Boltzmann system:
\begin{equation}
\begin{split}
\partial _tF_{+}+v\cdot \nabla _xF_{+}+(E+v\times B)\cdot \nabla
_vF_{+}&=Q(F_{+},F_{+})+Q(F_{+},F_{-}),\\
\partial _tF_{-}+v\cdot \nabla _xF_{-}-(E+v\times B)\cdot \nabla
_vF_{-}&=Q(F_{-},F_{+})+Q(F_{-},F_{-}),
\end{split}
\end{equation}
with initial data $F_{\pm }(0,x,v)=F_{0,\pm }(x,v)$. For
notational simplicity we have set all physical constants to be
unity, see \cite{Guo1} for more background. Here $F_{\pm
}(t,x,v)\ge 0$ are the spatially periodic number density
functions for the ions (+) and electrons (-) respectively, at
time $t\ge 0$, position $x=(x_1,x_2,x_3)\in {\bf [}-\pi {\bf
,}\pi ]^3={\mathbb {T}}^3$ and velocity $v=(v_1,v_2,v_3)\in
{\mathbb{ R}}^3$. The collision between particles is given by the
standard Boltzmann collision operator $Q(G_1,G_2)$ with
hard-sphere interaction:
\begin{eqnarray}
&&Q(G_1,G_2)=\int_{{\mathbb{ R}}^3\times S^2}|(u-v)\cdot \omega
|\{G_1(v^{\prime })G_2(u^{\prime })-G_1(v)G_2(u)\}dud\omega,
\label{hard}
\end{eqnarray}
where $\text{ }v^{\prime }=v-[(v-u)\cdot \omega ]\omega\text{ }$
and $\text{ }u^{\prime }=u+[(v-u)\cdot \omega ]\omega.$

The self-consistent, spatially periodic electromagnetic field $%
[E(t,x),B(t,x)]$ in (1.1) is coupled with $F_{\pm }(t,x,v)$
through the Maxwell system:
\begin{equation}
\begin{split}
&\partial _tE-\nabla \times B=- \int_{{\mathbb{ R}}
^3}v(F_{+}-F_{-})dv,\;\;\nabla \cdot B=0,\\
&\partial _tB+\nabla \times E=0,\;\; \nabla \cdot E=
\int_{{\mathbb{ R}}^3}(F_{+}-F_{-})dv,
\end{split}
\end{equation}
with initial data $E(0,x)=E_0(x),\;B(0,x)=B_0(x)$.

It turns out that it is convenient to consider the sum and
difference of $F_{+}$ and $F_{-}$ as proposed in \cite{belm}.
Defining
\begin{equation}
F\equiv F_{+}+F_{-}\text{ and } G\equiv F_{+}-F_{-},\label{sd}
\end{equation}
(1.1) and (1.3) can be rewritten as following:
\begin{equation}
\begin{split}
&\partial _tF+v\cdot \nabla _xF+(E+v\times B)\cdot \nabla
_v G=Q(F, F),\\
&\partial _tG+v\cdot \nabla _xG+(E+v\times B)\cdot \nabla
_v F=Q(G,F),\\
&\partial _tE-\nabla \times B=-\int_{{\mathbb{ R}}
^3}v\; G\;dv,\;\;\nabla \cdot B=0,\\
&\partial _tB+\nabla \times E=0,\;\;\nabla \cdot E=\int_{{\mathbb{
R}}^3}G\;dv, \label{vmb}
\end{split}
\end{equation}
with initial data $F(0,x,v)=F_{0}(x,v)$, $G(0,x,v)=G_{0}(x,v)$,
$E(0,x)=E_0(x)$ and $B(0,x)=B_0(x)$.

Now we introduce the diffusive scaling to (\ref{vmb}): for any
$\varepsilon>0$,
\begin{equation}
\begin{split}
&\varepsilon\partial _tF^{\varepsilon}+v\cdot \nabla
_xF^{\varepsilon}+(E^{\varepsilon}+v\times B^{\varepsilon})\cdot
\nabla
_v G^{\varepsilon}=\frac{1}{\varepsilon}Q(F^{\varepsilon}, F^{\varepsilon}),\\
&\varepsilon\partial _tG^{\varepsilon}+v\cdot \nabla
_xG^{\varepsilon}+(E^{\varepsilon}+v\times B^{\varepsilon})\cdot
\nabla
_v F^{\varepsilon}=\frac{1}{\varepsilon}Q(G^{\varepsilon},F^{\varepsilon}),\\
&\varepsilon\partial _tE^{\varepsilon}-\nabla \times
B^{\varepsilon}=- \int_{{\mathbb{ R}}
^3}v \;G^{\varepsilon}\;dv,\;\;\nabla \cdot B^{\varepsilon}=0,\\
&\varepsilon\partial _tB^{\varepsilon}+\nabla \times
E^{\varepsilon}=0,\;\;\nabla \cdot E^{\varepsilon}=\int_{{\mathbb{
R}}^3}G^{\varepsilon}\;dv.\label{rvmb}
\end{split}
\end{equation}
For notational simplicity, we normalize the global Maxwellian as
\begin{equation}
\mu(v)=\frac{1}{(2\pi)^{3/2}}e^{-|v|^2/2}. \label{normalizedmu}%
\end{equation}
We consider the following formal expansion in $\varepsilon$ around
the equilibrium state $[F,G,E,B]=[\mu,0,0,0]$: for any $n\geq 1$,
\begin{equation}
\begin{split}
F^{\varepsilon}(t,x,v)&=\mu+\sqrt{\mu}\{\varepsilon
f_{1}(t,x,v)+\varepsilon
^{2}f_{2}(t,x,v)+...+\varepsilon^{n}f_{n}^{\varepsilon}(t,x,v)\},\\
G^{\varepsilon}(t,x,v)&=\sqrt{\mu}\{\varepsilon
g_{1}(t,x,v)+\varepsilon
^{2}g_{2}(t,x,v)+...+\varepsilon^{n}g_{n}^{\varepsilon}(t,x,v)\},\\
E^{\varepsilon}(t,x)&=\{\varepsilon E_{1}(t,x)+\varepsilon
^{2}E_{2}(t,x)+...+\varepsilon^{n}E_{n}^{\varepsilon}(t,x)\},\\
B^{\varepsilon}(t,x)&=\{\varepsilon B_{1}(t,x)+\varepsilon
^{2}B_{2}(t,x)+...+\varepsilon^{n}B_{n}^{\varepsilon}(t,x)\}.
\label{exp}
\end{split}
\end{equation}
To determine the coefficients $f_{1}(t,x,v),...,f_{n-1}(t,x,v);\;
g_{1}(t,x,v),...,g_{n-1}(t,x,v);\\ E_{1}(t,x),...,E_{n-1}(t,x);\;
B_{1}(t,x),...,B_{n-1}(t,x)$, we plug the formal diffusive
expansion
(\ref{exp}) into the rescaled equations (\ref{rvmb}):%
\begin{equation}
\begin{split}
(\varepsilon\partial_{t}+v&\cdot\nabla_{x})\{\varepsilon
f_{1}+...+\varepsilon
^{n}f_{n}^{\varepsilon}\}\\
+\frac{1}{\sqrt{\mu}}&\{\varepsilon(E_1+v\times
B_1)+...+\varepsilon^n(E_n^{\varepsilon}+v\times
B_n^{\varepsilon}) \}\cdot\nabla_{v}[\sqrt{\mu}\{\varepsilon
g_{1}+...+\varepsilon
^{n}g_{n}^{\varepsilon}\}] \\
=&\frac{1}{\varepsilon\sqrt{\mu}}Q(\mu+\sqrt{\mu }\{\varepsilon
f_{1}+...+\varepsilon^{n}f_{n}^{\varepsilon}\},\;\mu+\sqrt{\mu
}\{\varepsilon f_{1}+...+\varepsilon^{n}f_{n}^{\varepsilon}\}),\\
\ (\varepsilon\partial_{t}+v&\cdot\nabla_{x})\{\varepsilon
g_{1}+...+\varepsilon
^{n}g_{n}^{\varepsilon}\}\\
+\frac{1}{\sqrt{\mu}}&\{\varepsilon(E_1+v\times
B_1)+...+\varepsilon^n(E_n^{\varepsilon}+v\times
B_n^{\varepsilon}) \}\cdot\nabla_{v}[\mu+\sqrt{\mu}\{\varepsilon
f_{1}+...+\varepsilon
^{n}f_{n}^{\varepsilon}\}] \\
=&\frac{1}{\varepsilon\sqrt{\mu}}Q(\sqrt{\mu }\{\varepsilon
g_{1}+...+\varepsilon^{n}g_{n}^{\varepsilon}\},\;\mu+\sqrt{\mu
}\{\varepsilon
f_{1}+...+\varepsilon^{n}f_{n}^{\varepsilon}\}),\label{1/e}
\end{split}
\end{equation}
\begin{equation*}
\begin{split}
 \varepsilon\partial_{t}\{\varepsilon
E_{1}&+...+\varepsilon^{n}E_{n}^{\varepsilon}\}-\nabla\times
\{\varepsilon B_{1}+...+\varepsilon^{n}B_{n}^{\varepsilon}\}=-\int
_{\mathbb{R}^3} v\sqrt{\mu}\{\varepsilon
g_{1}+...\varepsilon^{n}g_{n}^{\varepsilon}\}dv,\\
\ \varepsilon\partial_{t}\{\varepsilon
B_{1}&+...+\varepsilon^{n}B_{n}^{\varepsilon}\}+\nabla\times
\{\varepsilon E_{1}+...+\varepsilon^{n}E_{n}^{\varepsilon}\}=0,\\
\ \nabla\cdot\{\varepsilon
E_{1}&+...+\varepsilon^{n}E_{n}^{\varepsilon}\}=\int_{\mathbb{R}^3}
\sqrt{\mu}\{\varepsilon
g_{1}+...\varepsilon^{n}g_{n}^{\varepsilon}\}dv,\\
\ \nabla\cdot\{\varepsilon
B_{1}&+...+\varepsilon^{n}B_{n}^{\varepsilon}\}=0.
\end{split}
\end{equation*}
To expand the right hand side $Q$ in the above, we define $L$ the
well-known linearized collision operator and $\mathcal{L}$
another linearized operater as
\begin{equation}
Lf\equiv-\frac{1}{\sqrt{\mu}}\{Q(\mu,\sqrt{\mu}f)+Q(\sqrt{\mu}f,\mu)\},
\label{lL1}
\end{equation}
\begin{equation}
\mathcal{L}g\equiv-\frac{1}{\sqrt{\mu}}Q(\sqrt{\mu}g,\mu),
\label{lL2}
\end{equation}
and the nonlinear collision operator $\Gamma$ as (non-symmetric)
\begin{equation}
\Gamma(f,g)\equiv\frac{1}{\sqrt{\mu}}Q(\sqrt{\mu}f,\sqrt{\mu}g).
\label{nlgamma}%
\end{equation}
Note that $Lf$ and $\mathcal{L}g$ can be written as following in
terms of $\Gamma$:
\begin{equation}
Lf =-\{\Gamma(\sqrt{\mu},f)+\Gamma(f,\sqrt{\mu})\},\;\;\;\;
\mathcal{L}g=-\Gamma(g,\sqrt{\mu}).\label{linear}
\end{equation}

Now we equate the coefficients on both sides of the equation
(\ref{1/e}) in front of different powers of the parameter
$\varepsilon$. Let
\[
f_{-1}= f_{0}\equiv0,\;\;g_{-1}=
g_{0}\equiv0,\;\;E_0\equiv0,\;\;B_0\equiv0
\]
to obtain
\begin{equation}
\begin{split}
&\partial_t f_{m}+v\cdot\nabla_x
f_{m+1}+\frac{1}{\sqrt{\mu}}\sum_{\substack{i+j=m+1\\i,j\geq1}}(E_i+v\times
B_i)\cdot\nabla_v(\sqrt{\mu}g_j)\\
&\;\;\;\;=-Lf_{m+2}+\sum
_{\substack{i+j=m+2\\i,j\geq1}}\Gamma(f_{i},f_{j}),\label{key}\\
&\partial_t g_{m}+v\cdot\nabla_x
g_{m+1}+\frac{1}{\sqrt{\mu}}\sum_{\substack{i+j=m+1\\i,j\geq1}}(E_i+v\times
B_i)\cdot\nabla_v(\sqrt{\mu}f_j)-E_{m+1}\cdot v\sqrt{\mu}\\
&\;\;\;\;=-\mathcal{L}g_{m+2}+\sum
_{\substack{i+j=m+2\\i,j\geq1}}\Gamma(g_{i},f_{j}),
\end{split}
\end{equation}
for $-1\leq m \leq n-3$ as well as
\begin{equation}
\begin{split}
 &\partial_t
E_{m}-\nabla\times B_{m+1}=-\int v
g_{m+1}\sqrt{\mu}dv,\;\;\nabla\cdot B_{m+1}=0,\\
&\partial_t B_{m}+\nabla\times E_{m+1}=0,\;\;\nabla\cdot
E_{m+1}=\int g_{m+1}\sqrt{\mu} dv,\label{keyeb}
\end{split}
\end{equation}
for $0\leq m \leq n-2$. Moreover, we can collect terms left in
(\ref{1/e}) with powers $\varepsilon^{n-1}$ or higher and divide
by $\varepsilon ^{n-1}$ to get the equations for the remainders
$f_{n}^{\varepsilon},g_{n}^{\varepsilon},
E_{n}^{\varepsilon},B_{n}^{\varepsilon}$. Note that all the
$\varepsilon^{m+1}$-th order terms vanish for $m\leq n-3$ because
of (\ref{key}). First, we write equations for $f_n^{\varepsilon}$
and $g_n^{\varepsilon}$:
\begin{equation}
\begin{split}
&\varepsilon^{2}\partial_{t}f_{n}^{\varepsilon}+\varepsilon
v\cdot\nabla
_{x}f_{n}^{\varepsilon}+Lf_{n}^{\varepsilon}=\\
&\{-\partial_{t}f_{n-2}-v\cdot\nabla_{x}f_{n-1}-
\frac{1}{\sqrt{\mu}}\sum_{\substack{i+j=n-1\\i,j\geq1}}(E_i+v\times
B_i)\cdot\nabla_v(\sqrt{\mu}g_j)+\sum_{\substack{i+j=n
\\i,j\geq1}}\Gamma(f_{i},f_{j})\}\\ \label{remainderf}
&+\varepsilon\{-\partial_{t}f_{n-1}-\frac{1}{\sqrt{\mu}}
\sum_{\substack{i+j=n\\i,j\geq1}}(E_i+v\times
B_i)\cdot\nabla_v(\sqrt{\mu}g_j)+\sum_{\substack{i+j=n+1
\\i,j\geq1
}}\Gamma(f_{i},f_{j})\}\\
&+\varepsilon^{n}\Gamma(f_{n}^{\varepsilon}
,f_{n}^{\varepsilon})+\sum_{i=1}^{n-1}\varepsilon^{i}\{\Gamma(f_{n}
^{\varepsilon},f_{i})+\Gamma(f_{i},f_{n}^{\varepsilon})\}+\sum_{i+j\geq
n+2}\varepsilon^{i+j-n}\Gamma(f_{i},f_{j})\\
&-\frac{\varepsilon^{n+1}}{\sqrt{\mu}}(E_n^{\varepsilon}+v\times
B_n^{\varepsilon})
\cdot\nabla_v(\sqrt{\mu}g_n^{\varepsilon})\\
&-\frac{1}{\sqrt{\mu}}\sum_{i=1}^{n-1}\varepsilon^{i+1}\{(E_i+v\times
B_i)
\cdot\nabla_v(\sqrt{\mu}g_n^{\varepsilon})+(E_n^{\varepsilon}+v\times
B_n^{\varepsilon}) \cdot\nabla_v(\sqrt{\mu}g_i)\}\\
&-\frac{1}{\sqrt{\mu}} \sum_{i+j\geq
n+1}\varepsilon^{i+j-n+1}(E_i+v\times
B_i)\cdot\nabla_v(\sqrt{\mu}g_j);\\
\\
&\varepsilon^{2}\partial_{t}g_{n}^{\varepsilon}+\varepsilon
v\cdot\nabla _{x}g_{n}^{\varepsilon}-\varepsilon
E_n^{\varepsilon}\cdot v\sqrt{\mu}
+\mathcal{L}g_{n}^{\varepsilon}=\\
&\{-\partial_{t}g_{n-2}-v\cdot\nabla_{x}g_{n-1}-
\frac{1}{\sqrt{\mu}}\sum_{\substack{i+j=n-1\\i,j\geq1}}(E_i+v\times
B_i)\cdot\nabla_v(\sqrt{\mu}f_j)+E_{n-1}\cdot v\sqrt{\mu}\\
&+\sum_{\substack{i+j=n
\\i,j\geq1}}\Gamma(g_{i},f_{j})\}
+\varepsilon\{-\partial_{t}g_{n-1}-\frac{1}{\sqrt{\mu}}
\sum_{\substack{i+j=n\\i,j\geq1}}(E_i+v\times
B_i)\cdot\nabla_v(\sqrt{\mu}f_j)\\
&+\sum_{\substack{i+j=n+1
\\i,j\geq1
}}\Gamma(g_{i},f_{j})\}+\varepsilon^{n}\Gamma(g_{n}^{\varepsilon}
,f_{n}^{\varepsilon})+\sum_{i=1}^{n-1}\varepsilon^{i}\{\Gamma(g_{n}
^{\varepsilon},f_{i})+\Gamma(g_{i},f_{n}^{\varepsilon})\}\\
&+\sum_{i+j\geq
n+2}\varepsilon^{i+j-n}\Gamma(g_{i},f_{j})-\frac{\varepsilon^{n+1}}
{\sqrt{\mu}}(E_n^{\varepsilon}+v\times
B_n^{\varepsilon})
\cdot\nabla_v(\sqrt{\mu}f_n^{\varepsilon})\\
&-\frac{1}{\sqrt{\mu}}\sum_{i=1}^{n-1}\varepsilon^{i+1}\{(E_i+v\times
B_i)
\cdot\nabla_v(\sqrt{\mu}f_n^{\varepsilon})+(E_n^{\varepsilon}+v\times
B_n^{\varepsilon}) \cdot\nabla_v(\sqrt{\mu}f_i)\}\\
&-\frac{1}{\sqrt{\mu}} \sum_{i+j\geq
n+1}\varepsilon^{i+j-n+1}(E_i+v\times
B_i)\cdot\nabla_v(\sqrt{\mu}f_j).\\
\end{split}
\end{equation}
Similarly, by (\ref{keyeb}), we get the remainders for
$E_n^{\varepsilon}, B_n^{\varepsilon}$:
\begin{equation}
\begin{split}
\varepsilon\partial_{t}E_n^{\varepsilon}-\nabla\times
B_n^{\varepsilon}=&-\partial_t E_{n-1}-\int_{\mathbb{R}^3} v
g_{n}^{\varepsilon}\sqrt{\mu}dv,\;\;\nabla\cdot
B_n^{\varepsilon}=0,\\\label{remaindereb}
\varepsilon\partial_{t}B_n^{\varepsilon}+\nabla\times
E_n^{\varepsilon}=&-\partial_t B_{n-1},\;\; \nabla\cdot
E_n^{\varepsilon}=\int_{\mathbb{R}^3}
g_n^{\varepsilon}\sqrt{\mu}dv.
\end{split}
\end{equation}

The fluid equations can be obtained through the conditions
(\ref{key}) and (\ref{keyeb}). We first recall that the operator
$L\geq0,$ and for
any fixed $(t,x),$ the null space of $L$ is generated by $[\sqrt{\mu}%
,v\sqrt{\mu},|v|^{2}\sqrt{\mu}]$. For any function $f(t,x,v)$ we
thus can decompose
\[
f=\mathbf{P_1}f+\{\mathbf{I-P_1}\}f
\]
where $\mathbf{P_1}f$ (hydrodynamic part) is the $L_{v}^{2}$
projection on the null space for $L$ for given $(t,x).$ We can
further denote
\begin{equation}
\mathbf{P_1}f=\{\rho_{f}(t,x)+v\cdot u_{f}(t,x)+(\frac{|v|^{2}}{2}-\frac{3}%
{2})\theta_{f}(t,x)\}\sqrt{\mu}. \label{hfield}%
\end{equation}
Here we define the \textit{hydrodynamic field} of $f$ to be
\[
\lbrack\rho_{f}(t,x),u_{f}(t,x),\theta_{f}(t,x)]
\]
which represents the density, velocity and temperature
fluctuations physically. For the velocity field $u_{f}(t,x)$, we
further define its \textit{divergent-free} part as
\[
P_{0}u_{f}=\text{ the divergent-free projection of }u_{f}%
\]
so that
\[
\nabla\cdot\{P_{0}u_{f}\}\equiv0.
\]
Similarly, one can show that $\mathcal{L}\geq 0$ and for any
fixed $(t,x),$ the null space of $\mathcal{L}$ is one dimensional
vector space generated by $[\sqrt{\mu}].$ Likewise, any $g(t,x,v)$
can be decomposed into
\[
g=\mathbf{P_2}g+\{\mathbf{I-P_2}\}g
\]
where $\mathbf{P_2}g$ (hydrodynamic part) is the $L_{v}^{2}$
projection on the null space for ${\mathcal{L}}$ for given
$(t,x).$ We can further denote
\begin{equation}
\mathbf{P_2}g=\sigma_{g}(t,x)\sqrt{\mu}.\label{hfield2}%
\end{equation}
Here $\sigma_{g}(t,x)$, the \textit{hydrodynamic field} of $g$,
can be interpreted as the concentration difference. For more
details about $\mathcal{L}$ and $\mathbf{P_2}$, we refer
\cite{belm}. Before going on, we state the coercivity of $L$ and
$\mathcal{L}$ which will be often used in the subsequent sections:
there exists a $\delta>0$ such that
\begin{equation}
\langle Lf,f\rangle\geq\delta|\mathbf{(I-P_1)}f|_{\nu}^{2},\;\;
\langle
\mathcal{L}g,g\rangle\geq\delta|\mathbf{(I-P_2)}g|_{\nu}^{2}.
\label{lower}
\end{equation}
See Lemma 1 in \cite{Guo1} for its proof. Note that the operator
$L$ defined in \cite{Guo1} is equivalent to $[L,\mathcal{L}]$ in
our case.

Now define $[ \rho_{m},u_{m},\theta_{m},\sigma_m] $ to be the
corresponding hydrodynamic field of the $m$-th coefficients
$f_{m}$ and $g_m$. As for the first coefficients $f_{1}(t,x,v)$
and $g_1(t,x,v)$, from (\ref{key})
\begin{equation}
\{\mathbf{I-P_1}\}f_{1}=0 \;\text{ and
}\;\{\mathbf{I-P_2}\}g_{1}=0\label{micro}
\end{equation}
which immediately yield that
\begin{equation}
B_1=0 \;\text{ and }\; E_1=\nabla\phi_1\label{e1b1}
\end{equation}
up to constant and for some function $\phi_1(t,x)$ satisfying
$\triangle\phi_1=\sigma_1$; in particular, $B_1$ may be assumed to
be zero physically in a sense that nonzero constants $B_1$ do not
cause the hydrodynamic equations
(\ref{navier-stokes})-(\ref{fourier}) to change. It will be shown
in Lemma \ref{fluide} that its velocity fluctuation $u_{1}(t,x)$
is incompressible:
\begin{equation}
\nabla\cdot u_{1}\equiv0\;\text{ or }\;u_{1}=P_{0}u_{1}, \label{incon}%
\end{equation}
and its density and temperature fluctuations $\rho_{1}(t,x)$ and
$\theta _{1}(t,x)$ satisfy the Boussinesq relation:
\begin{equation}
\rho_{1}+\theta_{1}\equiv0. \label{bous}%
\end{equation}
Moreover, $[u_{1},\theta_{1},\sigma_1]$ satisfies the nonlinear
incompressible Vlasov-Navier-Stokes-Fourier equations:
\begin{align}
\partial_{t}u_{1}+u_{1}\cdot\nabla u_{1}+\nabla p_{1}  &  =\eta\Delta
u_{1}+\sigma_1\nabla\phi_1,\label{navier-stokes}\\
\partial_t \sigma_1+u_1\cdot\nabla
\sigma_1&=\alpha\Delta\sigma_1-\alpha\sigma_1,\label{vlasov}\\
\Delta\phi_1&=\sigma_1,\label{poisson}\\
\partial_{t}\theta_{1}+u_{1}\cdot\nabla\theta_{1}  &
=\kappa\Delta\theta_{1},
\label{fourier}
\end{align}
where $p_{1}(t,x)$ is the pressure and $\eta$, $\kappa$,
$\alpha>0$ are physical constants.

As for the coefficients $f_{m}(t,x,v),\;g_m(t,x,v)$ for $m\geq2,$
by (\ref{key}), the microscopic part of $f_{m}$ and $g_m$ is
determined by:
\begin{align}
\{\mathbf{I-P_1}\}f_{m}=L^{-1}\{-\partial_{t}f&_{m-2}-v\cdot\nabla_{x}%
f_{m-1}+\sum_{\substack{i+j=m\\i,j\geq1}}\Gamma(f_{i},f_{j})\nonumber\\
-\frac
{1}{\sqrt{\mu}}&\sum_{\substack{i+j=m-1\\i,j\geq1}}(E_i+v\times
B_i) \cdot\nabla_v(\sqrt{\mu}g_j)\},\label{hmicro1} \\
\{\mathbf{I-P_2}\}g_{m}=\mathcal{L}^{-1}\{-\partial_{t}g&_{m-2}-v\cdot\nabla_{x}%
g_{m-1}+\sum_{\substack{i+j=m\\i,j\geq1}}\Gamma(g_{i},f_{j})\nonumber\\
-\frac
{1}{\sqrt{\mu}}&\sum_{\substack{i+j=m-1\\i,j\geq1}}(E_i+v\times
B_i) \cdot\nabla_v(\sqrt{\mu}f_j)+E_{m-1}\cdot
v\sqrt{\mu}\}.\label{hmicro2}
\end{align}
On the other hand, for the hydrodynamic field of $f_{m}\text{ and
}g_m$:
\begin{equation*}
\begin{split}
\mathbf{P_1}f_{m}&=\{\rho_{m}(t,x)+v\cdot
u_{m}(t,x)+\{\frac{|v|^{2}}{2}-\frac
{3}{2}\}\theta_{m}(t,x)\}\sqrt{\mu},\\
\mathbf{P_2}g_{m}&=\;\sigma_m(t,x)\sqrt{\mu},
\end{split}
\end{equation*}
we can deduce an $m$-th order incompressibility condition
\begin{equation}
\nabla\cdot\{I-P_{0}\}u_{m}=-\partial_{t}\rho_{m-1},\\ \label{hincon}%
\end{equation}
an $m$-th order Boussinesq relation
\begin{equation}
\begin{split}
\rho_{m}+    \theta_{m}=\Delta^{-1}\nabla\cdot \{-u_{1}\cdot
\nabla(P_{0}u_{m-1})-P_{0}u_{m-1}\cdot\nabla u_{1}
+E_1\sigma_{m-1}+E_{m-1}\sigma_1\\
+R_{m-1}^{u}\}
+\langle\frac{|v|^{2}\sqrt{\mu}}{3},L^{-1}(\{\mathbf{I-P_1}\}v\cdot\nabla
_{x}\mathbf{P_1}f_{m-1})\rangle-\frac{5}{2}\theta_{1}\theta_{m-1}-u_{m-1}\cdot
u_{1}, \label{hbous}%
\end{split}
\end{equation}
and an $m$-th order \textit{linear} Vlasov-Navier-Stokes-Fourier
system for $[P_{0}u_{m},\theta_{m},\sigma_m]$:
\begin{align}
&(\partial_{t}+u_{1}\cdot\nabla-\eta\Delta)P_{0}u_{m}+P_{0}u_{m}\cdot\nabla
u_{1}+\nabla p_{m}-(E_1\sigma_m+E_m\sigma_1)=R_{m}^{u},\label{hnavier-stokes}\\
& (\partial_{t}+u_{1}\cdot\nabla-\alpha\Delta+\alpha)\sigma_m
+P_{0}u_{m}\cdot\nabla\sigma_1=R_m^{\sigma},\label{hvlasov}\\
&  (\partial_{t}+u_{1}\cdot\nabla-\kappa\Delta)\theta_{m}+P_{0}u_{m}%
\cdot\nabla\theta_{1}=R_{m}^{\theta}, \label{hfourier1}
\end{align}
\begin{align}
& \nabla\times E_m=-\partial_t B_{m-1},\;\; \nabla\cdot
E_m=\sigma_m,
\label{helec}\\
& \nabla\times B_m=\int_{\mathbb{R}^3} \{\mathbf{I-P_2}\}\;g_{m}
v\sqrt{\mu}dv +\partial_t E_{m-1},\;\;\nabla\cdot
B_m=0,\;\;\;\;\;\;\;\;\;\;\;\;\;\;\;\;\;\;\;\;\label{hmag}
\end{align}
with compatibility conditions coming from conservation laws
\begin{equation}
\frac{d}{dt}\int_{\mathbb{T}^3} E_m dx=-\alpha\int_{\mathbb{T}^3}
E_m dx + \ell_{m-1},\;\;\;\int_{\mathbb{T}^3} B_m
dx=0.\label{compaEB}
\end{equation}
Here $R_{m}^{u}$, $R_{m}^{\sigma}$, $R_{m}^{\theta}$ and
$\ell_{m-1}$, defined precisely in (\ref{ru}), (\ref{rv}),
(\ref{rs}) and (\ref{ell}), essentially depend only on
$f_{j},\;g_j,\; E_j,\; B_j$ for $j\leq m-1,$ since
$\{\mathbf{I-P_1}\}f_{m}$, $\{\mathbf{I-P_2}\}g_{m}$,
$\{I-P_{0}\}u_{m}$, as well as $\rho_{m}+\theta_{m}$ have been
determined.

In order to state  our results precisely in the next section, we
introduce the following norms and notations. We use
$\langle\cdot\,,\cdot\rangle$ to denote the standard $L^{2}$
inner product in $\mathbb{R}_{v}^{3},$ while we use
$(\cdot\,,\cdot)$ to denote $L^{2}$ inner product either in
$\mathbb{T}^{3}\times\mathbb{R}^{3}$ or in $\mathbb{T}^{3}$ with
corresponding $L^{2}$ norm $||\cdot||.$ We use the standard
notation $H^{s}$ to denote the Sobolev space $W^{s,2}.$ For the
Boltzmann collision operator (\ref{hard}), define the collision
frequency to be
\begin{equation}
\nu (v)\equiv \int_{\mathbb{R}^3} |v-v'|\mu(v')dv', \label{nu}
\end{equation}
which behaves like $|v|$ as $|v|\rightarrow\infty.$ It is natural
to define the following weighted $L^{2}$ norm to characterize the
dissipation rate.
\[
|g|_{\nu}^{2} \equiv\int_{\mathbb{R}^{3}}g^{2}(v)\nu(v)dv,\;\;\;
\;||g||_{\nu}^{2}   \equiv\int_{\mathbb{T}^{3}\times\mathbb{R}^{3}}%
g^{2}(x,v)\nu(v)dv dx.
\]
Observe that for hard sphere interaction,
\begin{equation}
||(1+|v|)^{\frac{1}{2}}g||\leq C ||g||_{\nu}.\label{weight}
\end{equation}
In order to be consistent with the hydrodynamic equations, we
define
\begin{equation}
\partial_{\gamma}^{\beta}=\partial_{x_{1}}^{\gamma_{1}}\partial_{x_{2}%
}^{\gamma_{2}}\partial_{x_{3}}^{\gamma_{3}}\partial_{v_{1}}^{\beta_{1}%
}\partial_{v_{2}}^{\beta_{2}}\partial_{v_{3}}^{\beta_{3}} \label{derivative}%
\end{equation}
where $\gamma=[\gamma_{1},\gamma_{2},\gamma_{3}]$ is related to
the space derivatives, while
$\beta=[\beta_{1},\beta_{2},\beta_{3}]$ is related to the
velocity derivatives.

We now define instant energy functionals and the dissipation
rate.\\

\textbf{Definition 1 (Instant Energy) }\textit{For
}$N\geq8$,\textit{\ for some constant }$C>0,$\textit{\ an instant
energy functional
}$\mathcal{E}_{N}(f,g,E,B)(t)\equiv\mathcal{E}_{N}(t)$
\textit{satisfies: }
\begin{equation}
\frac{1}{C}\mathcal{E}_{N}(t)\leq\sum_{|\beta|+|\gamma|\leq N}%
||[\partial_{\gamma}^{\beta}f,\;\partial_{\gamma}^{\beta}g]||^{2}(t)+
\sum_{|\gamma|\leq
N}||[\partial_{\gamma}E,\;\partial_{\gamma}B]||^2(t)\leq
C\mathcal{E}_{N}(t).
\label{ehard}%
\end{equation}\

\textbf{Definition 2 (Dissipation Rate) }\textit{For
}$N\geq8$,\textit{\ the dissipation rate }$\mathcal{D}_{N}%
(f,g)(t)$\textit{\ is defined as}
\begin{equation}
\begin{split}
\mathcal{D}_{N}(f,g)(t)=&\sum_{|\gamma|\leq N}||
[\partial_{\gamma}\mathbf{P_1}f,\;\partial_{\gamma}\mathbf{P_2}g]||^{2}(t)\\
+\frac{1}{\varepsilon^{2}} &\sum_{|\beta|+|\gamma|\leq
N}||[\partial_{\gamma}^{\beta}\{\mathbf{I-P_1}\}f,\;\partial_{\gamma}^{\beta}\{\mathbf{I-P_2}\}g]
||_{\nu}^{2}(t). \label{dhard}
\end{split}
\end{equation}\

We remark that both the instant energy and the dissipation rate
are carefully designed to capture the structure of the rescaled
Vlasov-Maxwell-Boltzmann equation (\ref{rvmb}). First of all, the
electromagnetic field $[E,B]$ is included only in the instant
energy, which prevents the exponential decay on $\mathcal{E}_N$
unlike the pure Boltzmann case for hard potentials. See
\cite{Guo2} and \cite{St}.
Notice that there is no $\frac{1}%
{\varepsilon^{2}}$ factor in front of the hydrodynamic part
$[\mathbf{P_1}f,\mathbf{P_2}g]$ in the dissipation rate
$\mathcal{D}_{N}(f,g)$, since only the microscopic part
$[\mathbf{\{I-P_1\}}f,\mathbf{\{I-P_2\}}g]$ should vanish as
$\varepsilon\rightarrow0$. For notational simplicity, the
\textit{Einstein's summation convention }is used for Greek letter
up to order $N\geq8$ throughout the paper, unless otherwise
specified. We denote $\nabla=\nabla_{x}$ and use $C$ to denote a
constant independent of $\varepsilon.$ We also use $U(\cdot)$ to
denote a general positive polynomial with $U(0)=0.$\\

\section{Main Results}\

The first result is to determine the coefficients
$f_{1},f_{2},...f_{m};\;g_{1},g_{2},...g_{m};\;E_1,E_2,\\...E_m;\;
B_1,B_2,...B_m$ in a diffusive approximation (\ref{exp}).\

\begin{thm}
\label{cecoe}Let $m$ divergent-free vector-valued functions $[u_{1}%
^{0}(x),u_{2}^{0}(x),...,u_{m}^{0}(x)]$, $2m$ scalar functions
$[\theta _{1}^{0}(x),\theta_{2}^{0}(x)...
\theta_{m}^{0}(x);\sigma_1^0(x),\sigma_2^0(x),...\sigma_m^0(x)]$
be given such that
\[
||u_{1}^{0}||_{H^{2}}+||\theta_{1}^{0}||_{H^{2}}+
||\sigma_{1}^{0}||_{H^{2}} \leq M
\]
and
\[
\begin{split}
\int_{\mathbb{T}^{3}}\sigma_r^0(x)dx=0,\;\;\int_{\mathbb{T}^{3}}u_r^0(x)dx=-
\sum_{\substack{i+j=r\\i,j\geq1}}\int_{\mathbb{T}^{3}}
E_i^0(x)\times B_j^0(x)\;dx,\\
\int_{\mathbb{T}^{3}}\frac{3}{2}\theta_r^0(x)dx=-\frac{1}{2}
\sum_{\substack{i+j=r\\i,j\geq1}}\int_{\mathbb{T}^{3}}
E_i^0(x)\cdot E_j^0(x)+B_i^0(x)\cdot B_j^0(x)\;dx,
\end{split}
\]
for $1\leq r\leq m.$ Here $E_i^0,B_j^0$ are defined by
$E_i(0,x),B_j(0,x)$ which have been inductively determined at the
precedents since $i,j<r$, starting with the average conditions
$\int_{\mathbb{T}^3}B_1^0dx=0$ and
$\int_{\mathbb{T}^3}E_1^0dx=0$. Then for sufficiently small $M$
and given $m$ real  vectors $e_1(=0),e_2,...,e_m$, there exist
unique functions $f_{1}(t,x,v), f_{2}(t,x,v),...,\\f_{m}(t,x,v);$
$g_{1}(t,x,v),g_{2}(t,x,v),...,g_{m}(t,x,v)$;
$E_1(t,x),E_2(t,x),...,E_m(t,x)$ and \\
$B_1(t,x),B_2(t,x),...,B_m(t,x)$ with
\begin{equation}
\begin{split}
\int_{\mathbb{T}^{3}}\sigma_r(t,x)dx=0,\;
\int_{\mathbb{T}^{3}}P_0u_r(t,x)dx=-
\sum_{\substack{i+j=r\\i,j\geq1}}\int_{\mathbb{T}^{3}}
E_i(t,x)\times B_j(t,x)\;dx,\\
\int_{\mathbb{T}^{3}}\frac{3}{2}\theta_r(t,x)dx=-\frac{1}{2}
\sum_{\substack{i+j=r\\i,j\geq1}}\int_{\mathbb{T}^{3}}
E_i(t,x)\cdot E_j(t,x)+B_i(t,x)\cdot B_j(t,x)\;dx,
\label{conservation}
\end{split}
\end{equation}
such that initially $P_{0}u_{r}(0,x)  =u_{r}^{0}(x),\;
\theta_{r}(0,x)  =\theta_{r}^{0}(x),\; \sigma_{r}(0,x)
=\sigma_r^0(x)$ and  $\int_{\mathbb{T}^3}E_r(0,x)dx=e_r$, and
$f_{1}(t,x,v),g_1(t,x,v),E_1(t,x),B_1(t,x)$ satisfy
(\ref{micro})-(\ref{fourier}) and
$f_{r}(t,x,v),g_r(t,x,v),E_r(t,x),B_r(t,x)$ satisfy
(\ref{hmicro1})-(\ref{compaEB}) for $2\leq r\leq m.$ Moreover,
for $1\leq r\leq m,$ for any $\beta,$ and for all $s\geq 3$,
there exists a polynomial $U_{r,\beta,s} $ with
$U_{r,\beta,s}(0)=0$ such that
\begin{equation}
\begin{split}
\sum_{|\tau|\leq s }\{&||[\partial_{\tau}^{\beta}f_{r},\;
\partial_{\tau}^{\beta}g_r||_{\nu}(t)+
||[\partial_{\tau}E_r,\;\partial_{\tau}B_r]||(t)\}\\
\leq&\;
e^{-\lambda t}U_{r,\beta,s}(\sum_{1\leq j\leq r}\{||u_{j}^{0}%
||_{H^{2s+4(r-j)}}+||\theta_{j}^{0}||_{H^{2s+4(r-j)}}
+||\sigma_{j}^{0}%
||_{H^{2s+4(r-j)}}\}), \label{rdecay}%
\end{split}
\end{equation}
where space-time derivatives
\[
\partial_{\tau}=\partial_{t}^{\tau_{0}}\partial_{x_{1}}^{\tau_{1}}%
\partial_{x_{2}}^{\tau_{2}}\partial_{x_{3}}^{\tau_{3}}%
\]
and $\lambda$ can be chosen as
$\frac{1}{4}\min\{\eta,\kappa,\alpha\}$ for sufficiently small
$M$.
\end{thm}\

We now turn to the most important question about the remainder
estimates for $f_{n}^{\varepsilon}$, $g_{n}^{\varepsilon}$,
$E_{n}^{\varepsilon}$ and $B_{n}^{\varepsilon}$. We first study
the classical case for the first order remainders
\[
f^{\varepsilon}\equiv
f_{1}^{\varepsilon},\;\;g^{\varepsilon}\equiv
g_{1}^{\varepsilon},\;\;E^{\varepsilon}\equiv
E_{1}^{\varepsilon},\;\;B^{\varepsilon}\equiv B_{1}^{\varepsilon}
\]
which satisfy the nonlinear Boltzmann type equations:
\begin{align}
\partial_{t}f^{\varepsilon}+\frac{1}{\varepsilon}v\cdot\nabla_{x}%
f^{\varepsilon}+\frac{1}{\varepsilon^{2}}Lf^{\varepsilon}=\frac{1}%
{\varepsilon}\Gamma(f^{\varepsilon},f^{\varepsilon})
-\frac{1}{\sqrt{\mu}}(E^{\varepsilon}+v\times
B^{\varepsilon})\cdot\nabla_v (\sqrt{\mu}g^{\varepsilon}),
\label{kinetic}\\
\partial_{t}g^{\varepsilon}+\frac{1}{\varepsilon}v\cdot(\nabla_{x}%
g^{\varepsilon}-\sqrt{\mu}E^{\varepsilon})+\frac{1}{\varepsilon^{2}}
\mathcal{L}g^{\varepsilon}=\frac{1}
{\varepsilon}\Gamma(g^{\varepsilon},f^{\varepsilon})-
\frac{1}{\sqrt{\mu}}
(E^{\varepsilon}+v\times B^{\varepsilon})\cdot\nabla_v
(\sqrt{\mu}f^{\varepsilon}),\nonumber
\end{align}
\begin{equation}
\begin{split}
&\varepsilon\partial_t E^{\varepsilon}-\nabla\times
B^{\varepsilon}=-\int g^{\varepsilon}v\sqrt{\mu}\;dv,
\;\;\nabla\cdot B^{\varepsilon}=0,\\
&\varepsilon\partial_t B^{\varepsilon}+\nabla\times
E^{\varepsilon}=0,\;\;\nabla\cdot E^{\varepsilon}=\int
g^{\varepsilon}\sqrt{\mu}\;dv.\label{kinetic1}
\end{split}
\end{equation}\

\begin{thm}
\label{nslimit}Let $N\geq8.$ Let
$f^{\varepsilon}(0,x,v)=f_{0}^{\varepsilon
}(x,v),\;g^{\varepsilon}(0,x,v)=g_{0}^{\varepsilon }(x,v)$ and
$E^{\varepsilon}(0,x)=E_0^{\varepsilon}(x),\;
B^{\varepsilon}(0,x)=B_0^{\varepsilon}(x)$ satisfy the mass,
momentum and energy conservation laws:
\begin{equation}
\begin{split}
(f_{0}^{\varepsilon},\sqrt{\mu})=0,\;\;
(g_{0}^{\varepsilon},\sqrt{\mu})=0,\label{con1}\\
(f_{0}^{\varepsilon},v\sqrt{\mu})+\varepsilon \int_{\mathbb{T}^3}
E_0^{\varepsilon}\times B_0^{\varepsilon} dx =0,\\
(f_{0}^{\varepsilon},|v|^2\sqrt{\mu})+\varepsilon
\int_{\mathbb{T}^3}
|E_0^{\varepsilon}|^2+ |B_0^{\varepsilon}|^2 dx =0.\\
\end{split}
\end{equation}
Then there exists an instant energy functional $\mathcal{E}%
_{N}(f^{\varepsilon},g^{\varepsilon},
E^{\varepsilon},B^{\varepsilon})(t)$ such that\\ if
$\;\mathcal{E}_{N}(f^{\varepsilon},g^{\varepsilon},
E^{\varepsilon},B^{\varepsilon})(0)$ is sufficiently small, then
\begin{equation}
\frac{d}{dt}\mathcal{E}_{N}(f^{\varepsilon},g^{\varepsilon},
E^{\varepsilon},B^{\varepsilon})(t)+\mathcal{D}_N%
(f^{\varepsilon},g^{\varepsilon})(t)\leq0. \label{en}%
\end{equation}
In particular,
\begin{equation}
\sup_{0\leq
t\leq\infty}\mathcal{E}_{N}(f^{\varepsilon},g^{\varepsilon},
E^{\varepsilon},B^{\varepsilon})(t)\leq
\mathcal{E}_{N}(f^{\varepsilon},g^{\varepsilon},
E^{\varepsilon},B^{\varepsilon})(0). \label{enbound}%
\end{equation}
Moreover, for $k\geq 1$ there exists $C_{N,k}>0$ such that
\begin{equation}
\mathcal{E}_{N}(f^{\varepsilon},g^{\varepsilon},
E^{\varepsilon},B^{\varepsilon})(t)\leq
C_{N,k}\mathcal{E}_{N+k}(f^{\varepsilon},g^{\varepsilon},
E^{\varepsilon},B^{\varepsilon})(0)\left\{1+\frac{t}{k}\right\}^{-k}.
\label{decayrate}
\end{equation}
\end{thm}

We remark that our initial data
$f_{0}^{\varepsilon},g_{0}^{\varepsilon},E_0^{\varepsilon},
B_0^{\varepsilon}$ are general and can contain initial layer for
the Vlasov-Navier-Stokes-Fourier limit. We can easily include
(one) temporal derivative $\partial_{t}$ in our definition of
instant energy and dissipation rate, and obtain the same uniform
bound for such a new
norm. With such a modification, the boundedness of $\partial_{t}%
f_{0}^{\varepsilon},\;\partial_t g_0^{\varepsilon},\;\partial_t
E_0^{\varepsilon},\;\partial_t B_0^{\varepsilon}$ automatically
removes the formation of any initial layer.

For higher order remainders
$f_{n}^{\varepsilon},\;g_{n}^{\varepsilon},\;E_{n}^{\varepsilon},\;
B_{n}^{\varepsilon}$ with $n\geq2,$ we have\

\begin{thm}
\label{high}Let $N\geq8.$ Given
$f_{1,}f_{2,}...f_{n};\;g_{1,}g_{2,}...g_{n};\;E_{1,}E_{2,}...E_{n};\;
B_{1,}B_{2,}...B_{n}$ constructed in Theorem \ref{cecoe} and let
\begin{equation}
\begin{split}
|||[&\mathbf{u}_{n}^{0},\mathbf{\theta}_{n}^{0},\mathbf{\sigma}_{n}^{0}]
|||_{N}(t)\\
&\equiv\sum_{1\leq j\leq
n}\{||u_{i}^{0}||_{H^{2N+10+4(n-j)}}+||\theta_{i}^{0}||_{H^{2N+10+4(n-j)}}
+||\sigma_{i}^{0}||_{H^{2N+10+4(n-j)}}\}. \label{0norm}
\end{split}
\end{equation}
And let
\begin{equation}
\begin{split}
F^{\varepsilon}(0,x,v)&\equiv\mu+\sqrt{\mu}\{\varepsilon f_{1}%
(0,x,v)+...+\varepsilon^{n-1}f_{n-1}(0,x,v)+\varepsilon^{n}f_{n}^{\varepsilon
}(0,x,v)\},\\ \label{initialdata}%
G^{\varepsilon}(0,x,v)&\equiv\sqrt{\mu}\{\varepsilon g_{1}%
(0,x,v)+...+\varepsilon^{n-1}g_{n-1}(0,x,v)+\varepsilon^{n}g_{n}^{\varepsilon
}(0,x,v)\},\\
E^{\varepsilon}(0,x)&\equiv\varepsilon E_{1}%
(0,x)+...+\varepsilon^{n-1}E_{n-1}(0,x)+\varepsilon^{n}E_{n}^{\varepsilon
}(0,x),\\
B^{\varepsilon}(0,x)&\equiv\varepsilon B_{1}%
(0,x)+...+\varepsilon^{n-1}B_{n-1}(0,x)+\varepsilon^{n}B_{n}^{\varepsilon
}(0,x)
\end{split}
\end{equation}
be given initial data satisfying the following conservation laws:
\begin{equation}
\begin{split}
\int_{\mathbb{T}^3}\int_{\mathbb{R}^3}\{F^{\varepsilon}(0,x,v)-\mu(v)\}
dvdx=0,\;\;
\int_{\mathbb{T}^3}\int_{\mathbb{R}^3}G^{\varepsilon}(0,x,v)dvdx=0,\label{con2}\\
\int_{\mathbb{T}^3}\int_{\mathbb{R}^3}vF^{\varepsilon}(0,x,v)dvdx+
\int_{\mathbb{T}^3}E^{\varepsilon}(0,x)\times B^{\varepsilon}(0,x)dx=0,\\
\int_{\mathbb{T}^3}\int_{\mathbb{R}^3}|v|^2\{F^{\varepsilon}(0,x,v)
-\mu(v)\} dvdx+
\int_{\mathbb{T}^3}|E^{\varepsilon}(0,x)|^2+|B^{\varepsilon}(0,x)|^2dx=0.
\end{split}
\end{equation}
Then there exist an instant energy functional $\mathcal{E}_{N}$
and a positive polynomial $U$ with $U(0)=0$ such that if both
$\varepsilon$ and
\[
\mathcal{E}_{N}(f_{n}^{\varepsilon}-f_{n},g_{n}^{\varepsilon}-g_{n}
,E_{n}^{\varepsilon}-E_{n},B_{n}^{\varepsilon}-B_{n})(0)
\]
are sufficiently small, then
\begin{equation}
\begin{split}
\sup_{0\leq t\leq\infty}&\mathcal{E}_{N}(f_{n}^{\varepsilon}-f_{n},%
g_{n}^{\varepsilon}-g_{n},E_{n}^{\varepsilon}-E_{n},
B_{n}^{\varepsilon}-B_{n})(t)\\
\leq&\{e^{U(|||[\mathbf{u}_{n}^{0},\mathbf{\theta}_{n}^{0},
\mathbf{\sigma}_{n}^{0}]|||_{N}%
)}\mathcal{E}_{N}(f_{n}^{\varepsilon}-f_{n},g_{n}^{\varepsilon}-g_{n}
,E_{n}^{\varepsilon}-E_{n},B_{n}^{\varepsilon}-B_{n})(0)\\
&+\varepsilon^{2}%
U(|||[\mathbf{u}_{n}^{0},\mathbf{\theta}_{n}^{0},
\mathbf{\sigma}_{n}^{0}]|||_{N})\}. \label{hen}%
\end{split}
\end{equation}
Moreover, for $k\geq 1$ there exists $C_{N,k}>0$ such that
\begin{equation}
\begin{split}
\mathcal{E}_{N}(f_n^{\varepsilon}-f_n&,g_n^{\varepsilon}-g_n,
E_n^{\varepsilon}-E_n,B_n^{\varepsilon}-B_n)(t)\\
\leq
C_{N,k}\{&e^{U(|||[\mathbf{u}_{n}^{0},\mathbf{\theta}_{n}^{0},
\mathbf{\sigma}_{n}^{0}]|||_{N}%
)}\mathcal{E}_{N+k}(f_{n}^{\varepsilon}-f_{n},g_{n}^{\varepsilon}-g_{n}
,E_{n}^{\varepsilon}-E_{n},B_{n}^{\varepsilon}-B_{n})(0)\\
&+\varepsilon^{2}%
U(|||[\mathbf{u}_{n}^{0},\mathbf{\theta}_{n}^{0},
\mathbf{\sigma}_{n}^{0}]|||_{N})\}\left\{1+\frac{t}{k}\right\}^{-k}.
\label{hdecayrate}
\end{split}
\end{equation}
\end{thm}\

Note that the conservation laws (\ref{con2}) in Theorem \ref{high}
imply the conservation laws for
$f_{n}^{\varepsilon},\;g_{n}^{\varepsilon},\;E_{n}^{\varepsilon},\;
B_{n}^{\varepsilon}$ due to (\ref{conservation}). See
(\ref{hmme}) for the precise description. In addition, we remark
that the positivity for the initial data $(F^{\varepsilon}\pm
G^{\varepsilon})(0,x,v)\geq 0$ which is equivalent to
$F_{\pm}^{\varepsilon} (0,x,v)\geq 0$ through the relation
(\ref{sd}) can be verified in the same way
as discussed in the Appendix of \cite{Guo2}.\\

There are several macroscopic fluid models for classifying the
dynamics of a charged fluid, but none has been derived from the
Boltzmann theory mathematically. This is because the construction
of the global solution to the important Vlasov-Maxwell-Boltzmann
system had been open for a long time until only a few years ago,
in \cite{Guo1}, a unique global in time classical solution near a
global Maxwellian for such a master system was constructed.
However, it still remains a major open problem: to construct
global renormalized solutions to the same system.

Singular limit problems emanating from the Boltzmann equations
have been studied by many people for decades
\cite{BGL1,BGL2,BGL3,DL,GL,GS,Guo2,LM,LM1,LM2,MS,Sa,UA}-- a great
overview of the issue is given in \cite{Guo2,M,V}. In particular,
in the recent work \cite{Guo2},  higher order approximations with
the unified energy method have now been shown to give rise to a
rigorous passage from the Boltzmann equation to the
Navier-Stokes-Fourier systems beyond the Navier-Stokes
approximation.

In this article, we rigorously establish the global in time
validity of the diffusive expansion (\ref{exp}) to the rescaled
Vlasov-Maxwell-Boltzmann equations (\ref{rvmb}) for any order.
This not only gives the uniform estimates for the higher order
corrections, but also leads to the mathematical derivation of new
dissipative hydrodynamic equations, which we call
Vlasov-Navier-Stokes-Fourier System. An interesting feature of
our result is that for the chosen incompressible regime, the
nontrivial magnetic effect in the limit first occurs in the
second order expansion $(n=2)$, while the electric effect is
important at any order. We believe our result opens a new line of
research for those macroscopic approximations both physically and
mathematically.

The method of this paper is based on the improvement of the
recently developed nonlinear energy method in \cite{Guo1,Guo2}.
We use the reformulation (\ref{vmb}), by introducing new unknowns
(\ref{sd}), to simplify both the character of the
Vlasov-Maxwell-Boltzmann system and the analysis presented in this
paper endowed with the clear, concrete linearized collision
operators $L$ and $\mathcal{L}$. Macroscopic equations and local
conservation laws are proven to be key tools to build the crucial
positivity of $L$ and $\mathcal{L}$. The most important
analytical difficulty lies in deriving the uniform estimates on
$f^{\varepsilon}$ and $g^{\varepsilon}$. Due to the singular
behavior of the time derivative in our problem, the positivity
estimate for purely spatial derivatives is invoked: see
(\ref{diff}) in Lemma \ref{abc}. Notice that an additional term
$\frac{d}{dt}G(t)$ is needed. Such a differential form can still
yield decay estimates with the notion of equivalent instant
energies. The local conservation laws are used to estimate the
more singular temporal derivative for the hydrodynamic field
directly in terms of purely spatial derivatives of the
microscopic part. We also use the trick of the integration by
parts in the time variable and, in turn, by using such a uniform
estimate, avoid encountering
$[\partial_t\mathbf{\{I-P_1\}}f^{\varepsilon},\;\partial_t
\mathbf{\{I-P_2\}}g^{\varepsilon}]$. This complication is another
reason for introducing the notion of equivalent instant energies.

As pointed out in \cite{Guo1}, the Vlasov-Maxwell-Boltzmann system
is tricky to handle due to the hyperbolic nature of Maxwell
equation and indeed, it is the most intriguing and critical part
of this article to control electromagnetic fields. For that, we
first utilize the macroscopic equation (\ref{e}) to estimate
electric fields and then the Maxwell equation itself for magnetic
fields. Noting that the macroscopic equation (\ref{e}) is simpler
than the one considered in \cite{Guo1} owing to the reformulation
(\ref{vmb}), we remark that  our method together with the
positivity estimate (\ref{diff}) provides another lucid and
concise way of proving the global in time classical solution for
the Vlasov-Maxwell-Boltzmann equations, especially without using
the temporal derivatives.

On the other hand, it is delicate to establish the well-posedness
of the new hydrodynamic equations because of their complexity,
mainly stemming from electromagnetic fields. It turns out that the
compatibility conditions (\ref{compaEB}) for averages of
electromagnetic fields are necessary and sufficient conditions in
order for the Vlasov-Navier-Stokes-Fourier systems to have the
unique solution. We should point out that those conditions are
natural restrictions in that they can be derived from
conservation laws. In turn, the exponential decay of hydrodynamic
variables as well as the solvability of hydrodynamic equations
lead to the solvability of kinetic equations and almost the same
decay for approximate solutions. While the same exponential decay
rate is obtained for the pure Boltzmann case, we are able to
obtain only the polynomial decay rate (\ref{decayrate}) and
(\ref{hdecayrate}) for solutions of the Vlasov-Maxwell-Boltzmann
equations  (compare with (2.6) and (2.13) in \cite{Guo2}) by
applying the method proposed in \cite{St}. In particular, for
higher order remainders the more sophisticated continuity
argument is employed in order to get the desired polynomial decay
rate (\ref{hdecayrate})
 from (\ref{hegd}).

The paper will proceed as follows: we will derive the high order
linear Vlasov-Navier-Stokes-Fourier system in Section 3; we will
prove Theorem \ref{cecoe}  for coefficients $f_{1}(t,x,v),
...,f_{m}(t,x,v),\;
g_{1}(t,x,v),...,g_{m}(t,x,v),\;E_1(t,x),...,E_m(t,x),\;
B_1(t,x),..,\\B_m(t,x)$ in Section 4; Section 5 will be devoted to
the positivity of $L$ and $\mathcal{L}$; in the last two
sections, the first and higher order remainder estimates--Theorem
\ref{nslimit} and Theorem \ref{high}-- will be proven respectively.\\

\section{High Order Vlasov-Navier-Stokes-Fourier System:\\
 Formal Derivation}\

In this section, we derive the microscopic equations
(\ref{hmicro1}), (\ref{hmicro2}) and hydrodynamic equations
(\ref{hincon})-
(\ref{compaEB}). We shall use many results from \cite{Guo2}.

\begin{lem}
\label{fluide}Assume that the expansion (\ref{exp})
 satisfies (\ref{rvmb}) and such that\\
for  $|\tau|+|\beta|\leq N$,
\begin{equation}
\begin{split}
\sum_{m=1}^{n-1}\{||[\partial_{\tau}^{\beta}f_{m},
\partial_{\tau}^{\beta}g_m]||_{\nu}+&||[\partial_{\tau}E_m,
\partial_{\tau}B_m]||\}\\&+ ||[\partial_{\tau}^{\beta}f_{n}^{\varepsilon},
\partial_{\tau}^{\beta}g_n^{\varepsilon}]||_{\nu}+||[\partial_{\tau}E_n
^{\varepsilon},
\partial_{\tau}B_n^{\varepsilon}]||<\infty. \label{hassume}%
\end{split}
\end{equation}
Then there exist
$\;f_{n}(t,x,v),\;g_n(t,x,v),\;f_{n+1}(t,x,v),\;g_{n+1}(t,x,v)\;$
and $\;E_n(t,x),\;\\B_n(t,x)$ such that (\ref{key}) and
(\ref{keyeb}) are valid for all $m\leq n-1.$ Moreover, $f_{1}$ and
$g_1$ satisfy (\ref{micro}), the incompressibility condition
(\ref{incon}), the Boussinesq relation (\ref{bous}), and the
first order Vlasov-Navier-Stokes-Fourier system
(\ref{navier-stokes})-
(\ref{fourier}); $E_{1}$ and $B_1$ satisfy (\ref{e1b1}). For
$m\geq2,$ $f_{m},g_m,E_m,B_m$ satisfy the microscopic equation
(\ref{hmicro1}) and (\ref{hmicro2}), the m-th order
incompressibility condition (\ref{hincon}), the m-th order
Boussinesq relation (\ref{hbous}), and the m-th order
Vlasov-Navier-Stokes-Fourier system
(\ref{hnavier-stokes})-(\ref{compaEB}) with

\begin{align}
  R_{m}^{u}\equiv&\;\langle v\cdot\nabla_{x}L^{-1}\{\partial_{t}%
\{\mathbf{I-P_1\}}f_{m-1}-\sum_{\substack{i+j=m+1, \\i,j>1}}\Gamma(f_{i}%
,f_{j})\},v\sqrt{\mu}\rangle\nonumber\\
&  +\langle
v\cdot\nabla_{x}L^{-1}\{\{\mathbf{I-P_1\}(}v\cdot\nabla
_{x}\{\mathbf{I-P_1\}}f_{m}),v\sqrt{\mu}\rangle\nonumber\\
&  -\langle
v\cdot\nabla_{x}L^{-1}\{\Gamma(f_{1},\{\mathbf{I-P_1\}}f_{m})
+\Gamma(\{\mathbf{I-P_1\}}f_{m}%
,f_{1})\},v\sqrt{\mu}\rangle\nonumber\\
& +\langle v\cdot\nabla_{x}L^{-1}\{\frac{1}{\sqrt{\mu}}
\sum_{\substack{i+j=m,
\\i,j\geq1}}(E_i+v\times B_i)\cdot\nabla_v(g_j\sqrt{\mu})\},
v\sqrt{\mu}\rangle\label{ru}\\
& -(\partial_{t}+u_{1}\cdot\nabla-\eta\Delta)\{I-P_{0}\}u_{m}-\{I-P_{0}%
\}u_{m}\cdot\nabla u_{1}\nonumber\\
&  -(\nabla\cdot\{I-P_{0}\}u_{m})u_{1}+\frac{\eta}{3}\nabla(\nabla
\cdot\{I-P_{0}\}u_{m})\nonumber\\
&+\sum_{\substack{i+j=m+1,
\\i,j>1}}\{E_i\nabla\cdot E_j-(\partial_t E_{j-1}-\nabla\times B_j)
\times B_i\};\nonumber
\end{align}
\begin{align}
  R_{m}^{\sigma}\equiv&\;\langle v\cdot\nabla_{x}\mathcal{L}^{-1}
  \{\partial
_{t}\{\mathbf{I-P_2\}}g_{m-1}-\sum_{\substack{i+j=m+1,
\\i,j>1}}\Gamma
(g_{i},f_{j})\},\sqrt{\mu}\rangle\nonumber\\
&  +\langle
v\cdot\nabla_{x}\mathcal{L}^{-1}\{\{\mathbf{I-P_2\}(}v\cdot\nabla
_{x}\{\mathbf{I-P_2\}}g_{m}),\sqrt{\mu}\rangle
\nonumber\\
&  -\langle
v\cdot\nabla_{x}\mathcal{L}^{-1}\{\Gamma(g_{1},\{\mathbf{I-P_2\}}g_{m})
+\Gamma(\{\mathbf{I-P_2\}}g_{m}%
,g_{1})\},\sqrt{\mu}\rangle\label{rv}\\
& +\langle v\cdot\nabla_{x}\mathcal{L}^{-1}\{\frac{1}{\sqrt{\mu}}
\sum_{\substack{i+j=m,
\\i,j\geq1}}(E_i+v\times B_i)\cdot\nabla_v(f_j\sqrt{\mu})\},
\sqrt{\mu}\rangle\nonumber\\ &
-\{\nabla\cdot(I-P_{0})u_{m}\}\sigma_{1}-(I-P_{0})u_{m}\cdot\nabla
\sigma_{1};\nonumber
\end{align}
\begin{align}
  R_{m}^{\theta}\equiv&\;\langle v\cdot\nabla_{x}L^{-1}\{\partial
_{t}\{\mathbf{I-P_1\}}f_{m-1}-\sum_{\substack{i+j=m+1,
\\i,j>1}}\Gamma
(f_{i},f_{j})\},\frac{|v|^{2}\sqrt{\mu}}{5}\rangle\nonumber\\
&  +\langle
v\cdot\nabla_{x}L^{-1}\{\{\mathbf{I-P_1\}(}v\cdot\nabla
_{x}\{\mathbf{I-P_1\}}f_{m}),\frac{|v|^{2}\sqrt{\mu}}{5}\rangle
\nonumber\\
&  -\langle
v\cdot\nabla_{x}L^{-1}\{\Gamma(f_{1},\{\mathbf{I-P_1\}}f_{m})
+\Gamma(\{\mathbf{I-P_1\}}f_{m}%
,f_{1})\},\frac{|v|^{2}\sqrt{\mu}}{5}\rangle\nonumber\\
& +\langle v\cdot\nabla_{x}L^{-1}\{\frac{1}{\sqrt{\mu}}
\sum_{\substack{i+j=m,
\\i,j\geq1}}(E_i+v\times B_i)\cdot\nabla_v(g_j\sqrt{\mu})\},
\frac{|v|^{2}\sqrt{\mu}}{5}\rangle\label{rs}\\ &
+(\frac{2}{5}\partial_{t}-u_{1}\cdot\nabla+\kappa_{1}\Delta)\{\rho
_{m}+\theta_{m}\}\nonumber\\
&
-\{\nabla\cdot(I-P_{0})u_{m}\}\theta_{1}-(I-P_{0})u_{m}\cdot\nabla
\theta_{1}\nonumber\\
&-\frac{2}{5}\sum_{\substack{i+j=m+1,\\i,j\geq1}} E_i\cdot
(\partial_t E_{j-1}-\nabla\times B_j);\nonumber
\end{align}
\begin{align}
 p_{m}\equiv&\;(\rho_{m+1}+\theta_{m+1})-\langle\frac{|v|^{2}\sqrt{\mu}}%
{3},L^{-1}\{\{\mathbf{I-P_1\}(}v\cdot\nabla_{x}\mathbf{P_1}f_{m})\}\rangle
\;\;\;\;\;\;\;\;\;\;\;\;\;\;\;
\label{press}\\
&  +\frac{5}{2}\theta_{1}\theta_{m}+(u_{m}\cdot u_{1}).\nonumber
\\ \nonumber
\end{align}
Also the conservation laws (\ref{conservation})  are valid for
each $m$.
\end{lem}\

\begin{proof}
Notice that (\ref{key}) is clearly valid for $m<n-2$ under the
assumption (\ref{hassume}) since it can be derived by letting
$\varepsilon\rightarrow 0$ after dividing (\ref{1/e}) by
$\varepsilon^{m+1}$. By the same token, (\ref{keyeb}) holds for
$m<n-1$. We now show the existence of the coefficients
$f_{n},g_n,E_n,B_n$ and $f_{n+1},g_{n+1}$ so that the equations
(\ref{key}) for $m=n-2,\;n-1$ and (\ref{keyeb}) for $m=n-1$ are
satisfied. Firstly, by (\ref{hassume}), up to a subsequence,
there exist $f_{n},g_n,E_n,B_n$ such that
\begin{align*}
[f_{n}^{\varepsilon},g_{n}^{\varepsilon}]\rightharpoonup
[f_n,g_n]\;\;\text{weakly in } ||\cdot||_{\nu}\;\text{ and }\;\;
[E_{n}^{\varepsilon},B_{n}^{\varepsilon}]\rightharpoonup
[E_n,B_n]\;\;\text{weakly in } ||\cdot||.
\end{align*}
The assumption (\ref{hassume}) guarantees the existence of their
derivatives. By subtracting off $Lf_{n},\;\mathcal{L}g_n$ on both
sides of (\ref{remainderf})  we can isolate the zeroth order
terms in the remainder equations and deduce (\ref{key}) for
$m=n-2$ at least in the sense of distributions. Moreover, letting
$\varepsilon\rightarrow 0$ of (\ref{remaindereb}), we also deduce
(\ref{keyeb}) for $m=n-1$ at least in the sense of distributions.

Next, in order to find $f_{n+1}\text{ and }g_{n+1}$, we use
(\ref{key}) for $m=n-2$ to simplify (\ref{remainderf}) as
\begin{align}
 L&\left\{  \frac{f_{n}^{\varepsilon}-f_{n}}{\varepsilon}\right\}
\nonumber\\
&=-\varepsilon\partial_t f_n^{\varepsilon}-v\cdot\nabla_x
f_n^{\varepsilon}\nonumber\\
&+\{-\partial_{t}f_{n-1}-\frac{1}{\sqrt{\mu}}
\sum_{\substack{i+j=n\\i,j\geq1}}(E_i+v\times
B_i)\cdot\nabla_v(\sqrt{\mu}g_j)+\sum_{\substack{i+j=n+1
\\i,j\geq1
}}\Gamma(f_{i},f_{j})\}\nonumber\\
&+\varepsilon^{n-1}\Gamma(f_{n}^{\varepsilon}
,f_{n}^{\varepsilon})+\sum_{i=1}^{n-1}\varepsilon^{i-1}\{\Gamma(f_{n}
^{\varepsilon},f_{i})+\Gamma(f_{i},f_{n}^{\varepsilon})\}
+\sum_{i+j\geq
n+2}\varepsilon^{i+j-n-1}\Gamma(f_{i},f_{j})\label{iso2f}\\
&-\frac{\varepsilon^{n}}{\sqrt{\mu}}(E_n^{\varepsilon}+v\times
B_n^{\varepsilon})
\cdot\nabla_v(\sqrt{\mu}g_n^{\varepsilon})\nonumber\\
&-\frac{1}{\sqrt{\mu}}\sum_{i=1}^{n-1}\varepsilon^{i}\{(E_i+v\times
B_i)
\cdot\nabla_v(\sqrt{\mu}g_n^{\varepsilon})+(E_n^{\varepsilon}+v\times
B_n^{\varepsilon}) \cdot\nabla_v(\sqrt{\mu}g_i)\}\nonumber\\
&-\frac{1}{\sqrt{\mu}} \sum_{i+j\geq
n+1}\varepsilon^{i+j-n}(E_i+v\times
B_i)\cdot\nabla_v(\sqrt{\mu}g_j)\; ;\nonumber\\
\nonumber
\end{align}
\begin{align}
 \mathcal{L}&\left\{  \frac{g_{n}^{\varepsilon}-g_{n}}{\varepsilon}\right\}
\nonumber\\
&=-\varepsilon\partial_{t}g_{n}^{\varepsilon}- v\cdot\nabla
_{x}g_{n}^{\varepsilon}+ E_n^{\varepsilon}\cdot v\sqrt{\mu}\nonumber\\
&+\{-\partial_{t}g_{n-1}-\frac{1}{\sqrt{\mu}}
\sum_{\substack{i+j=n\\i,j\geq1}}(E_i+v\times
B_i)\cdot\nabla_v(\sqrt{\mu}f_j)+\sum_{\substack{i+j=n+1\\i,j\geq1
}}\Gamma(g_{i},f_{j})\}\nonumber\\
&+\varepsilon^{n-1}\Gamma(g_{n}^{\varepsilon}
,f_{n}^{\varepsilon})+\sum_{i=1}^{n-1}\varepsilon^{i-1}\{\Gamma(g_{n}
^{\varepsilon},f_{i})+\Gamma(g_{i},f_{n}^{\varepsilon})\}\label{iso2g}\\
&+\sum_{i+j\geq
n+2}\varepsilon^{i+j-n-1}\Gamma(g_{i},f_{j})-\frac{\varepsilon^{n}}
{\sqrt{\mu}}(E_n^{\varepsilon}+v\times B_n^{\varepsilon})
\cdot\nabla_v(\sqrt{\mu}f_n^{\varepsilon})\nonumber\\
&-\frac{1}{\sqrt{\mu}}\sum_{i=1}^{n-1}\varepsilon^{i}\{(E_i+v\times
B_i)
\cdot\nabla_v(\sqrt{\mu}f_n^{\varepsilon})+(E_n^{\varepsilon}+v\times
B_n^{\varepsilon}) \cdot\nabla_v(\sqrt{\mu}f_i)\}\nonumber\\
&-\frac{1}{\sqrt{\mu}} \sum_{i+j\geq
n+1}\varepsilon^{i+j-n}(E_i+v\times
B_i)\cdot\nabla_v(\sqrt{\mu}f_j). \nonumber
\end{align}

Take the inner product of (\ref{iso2f}) with
$\{\mathbf{I-P_1}\}\left\{ \frac
{f_{n}^{\varepsilon}-f_{n}}{\varepsilon}\right\}.$ By
(\ref{lower}), the LHS of (\ref{iso2f}) is bounded from below by
\[
\delta||(\mathbf{I-P_1})\left\{
\frac{f_{n}^{\varepsilon}-f_{n}}{\varepsilon }\right\}
||_{\nu}^{2}.
\]
On the other hand, for the inner products in the RHS of
(\ref{iso2f}), by (\ref{weight}),
 we have%
\begin{align*}
&  (-\varepsilon\partial_{t}f_{n}^{\varepsilon}-v\cdot\nabla_{x}%
f_{n}^{\varepsilon}-\partial_{t}f_{n-1},\frac{1}{\varepsilon}\{\mathbf{I-P_1}%
\}\{f_{n}^{\varepsilon}-f_{n}\})\\
&\;\;\;\;\leq\{\varepsilon||\partial_{t}f_{n}^{\varepsilon}||+||
\nabla_{x}f_{n}^{\varepsilon}||_{\nu}+||\partial_{t}f_{n-1}||\}\times\frac
{1}{\varepsilon}||\{\mathbf{I-P_1}\}\{f_{n}^{\varepsilon}-f_{n}\}
||_{\nu};\\
&(-\frac{1}{\sqrt{\mu}}
\sum_{\substack{i+j=n\\i,j\geq1}}(E_i+v\times
B_i)\cdot\nabla_v(\sqrt{\mu}g_j),\frac{1}{\varepsilon}\{\mathbf{I-P_1}%
\}\{f_{n}^{\varepsilon}-f_{n}\})\\
&\;\;\;\;\leq\{\sum_{\substack{i+j=n\\i,j\geq1}}(||E_i||+||B_i||)
(||g_j||_{\nu}+||\nabla_v g_j||)\}\times\frac
{1}{\varepsilon}||\{\mathbf{I-P_1}\}\{f_{n}^{\varepsilon}-f_{n}\}
||_{\nu}.
\end{align*}
For the remaining terms, we use Lemma 3.3 in \cite{Guo2} to get
the upper bound as
\[
\begin{split}
\{\{(\sum_{i=1}^{n-1}||g_{i}||_{\nu}+||g_{n}^{\varepsilon}||_{\nu})
\cdot (\sum_{i=1}^{n-1}(||E_i||+||B_i||)+||E_n^{\varepsilon}||+
||B_n^{\varepsilon}||)\}\\
+(\sum_{i=1}^{n-1}||f_{i}||_{\nu}+||f_{n}^{\varepsilon}||_{\nu}
)^2 \}\times \frac{1}{\varepsilon}||\{\mathbf{I-P_1}\}
\{f_{n}^{\varepsilon}-f_{n}\}||_{\nu}
\end{split}
\]
We therefore conclude that $\frac{1}{\varepsilon}||\{\mathbf{I-P_1}%
\}\{f_{n}^{\varepsilon}-f_{n}\}||_{\nu}$ is uniformly bounded.
Similarly, one can show that
$\frac{1}{\varepsilon}||\{\mathbf{I-P_2}%
\}\{g_{n}^{\varepsilon}-g_{n}\}||_{\nu}$ is also uniformly
bounded. Hence there exist $f_{n+1},\;g_{n+1}$ such that
\[
\{\mathbf{I-P_1}\}\left\{
\frac{f_{n}^{\varepsilon}-f_{n}}{\varepsilon }\right\}
\rightharpoonup f_{n+1},\;\{\mathbf{I-P_2}\}\left\{
\frac{g_{n}^{\varepsilon}-g_{n}}{\varepsilon }\right\}
\rightharpoonup g_{n+1}\;\;\text{weakly in }||\cdot||_{\nu}.
\]
Subtracting off $Lf_{n+1}$, $\mathcal{L}g_{n+1}$ on both sides of
(\ref{iso2f}) and (\ref{iso2g}), we again isolate the zeroth
order term. Letting $\varepsilon\rightarrow0$ again, one can
deduce (\ref{key}) for $m=n-1.$

We now turn to the derivation of various hydrodynamic equations
based on (\ref{key}) and (\ref{keyeb}). For the pure Boltzmann
case, it is well known that the case $m=1$ yields the celebrated
nonlinear incompressible Navier-Stokes equations and $m>1$, the
Navier-Stokes-Fourier system of which derivation can be found in
 \cite{Guo2} (p.21-25).  We take the same path to
derive new dissipative PDE's from the rescaled
Vlasov-Maxwell-Boltzmann system. We will use numerous results
from \cite{Guo2}. First consider (\ref{key}) when $m=-1$. It is
equivalent to (\ref{micro}) and hence we get
\begin{equation}
\begin{split}
f_1=\{\rho_1+v\cdot
u_1+(\frac{|v|^2-3}{2})\theta_1\}\sqrt{\mu}\;\;\text{ and }\;\;
g_1=\;\sigma_1\sqrt{\mu}.\label{first}
\end{split}
\end{equation}
When $m=0$, electromagnetic field equations become $\nabla\times
B_1=0,\;\nabla\cdot B_1=0$, $\nabla\times E_1=0,$ and
$\nabla\cdot E_1=\sigma_1.$ Thus we may assume $B_1=0$ and
$E_1=\nabla\phi_1$ for some $\phi_1$ which satisfies
$\Delta\phi_1=\sigma_1$ under the restriction $\int B_1
dx=0\text{ and } \int E_1dx=0$; indeed, a nonzero constant $B_1$
does not produce any new macroscopic terms in (\ref{key}) when
$m=1$, since
\[
\begin{split}
v\times B_1\cdot\nabla_v(\sqrt{\mu}g_1)=v\times
B_1\cdot\nabla_v(\sigma_1\mu)=v\times B_1\cdot(-\sigma\mu)v=0,\\
\langle \frac{1}{\sqrt{\mu}}v\times
B_1\cdot\nabla_v(\sqrt{\mu}f_1),\sqrt{\mu}\rangle=
\langle\frac{1}{\sqrt{\mu}} v\times B_1\cdot\nabla_v(v\cdot
u_1\mu),\sqrt{\mu}\rangle=0,
\end{split}
\]
where the last equality is due to the integration by parts in $v$.
Recalling the kernel of $L$ and $\mathcal{L}$, by collision
invariant property,
\begin{equation}
\begin{split}
\langle v\cdot\nabla_x f_1 ,[1,v,|v|^{2}/2]\sqrt{\mu
}\rangle=0\;\;\text{ and }\;\; \langle v\cdot\nabla_x
g_1-E_1\cdot v\sqrt{\mu}, \sqrt{\mu}\rangle=0.\label{first2}
\end{split}
\end{equation}
Plugging (\ref{first}) into (\ref{first2}), the very first
equation gives rise to the incompressibility (\ref{incon}) and
Boussinesq relation (\ref{bous}). Microscopic part of $f_2$ and
$g_2$ can be written as following by solving (\ref{key}) when
$m=1$ for $Lf_2$ and $\mathcal{L}g_2$:
\begin{equation}
\begin{split}
\{\mathbf{I-P_1}\}f_2&=L^{-1}\{-v\cdot\nabla_x
f_1+\Gamma(f_1,f_1)\},\\
\{\mathbf{I-P_2}\}g_2&=\mathcal{L}^{-1}\{-v\cdot\nabla_x
g_1+\Gamma(g_1,f_1)+E_1\cdot v\sqrt{\mu}\}\label{mi1}\\
&=\mathcal{L}^{-1}(v\sqrt{\mu})\cdot(-\nabla_x\sigma_1+E_1)+\sigma_1
f_1,
\end{split}
\end{equation}
because of the following identity obtained from (\ref{linear}) and
(\ref{micro}):
\begin{equation}
\mathcal{L}^{-1}\Gamma(g_1,f_1)=\sigma_1
\mathcal{L}^{-1}\Gamma(\sqrt{\mu},f_1)=\sigma_1
\mathcal{L}^{-1}[-Lf_1+\mathcal{L}f_1]=\sigma_1 f_1.
\end{equation}
Notice that $\{\mathbf{I-P_1}\}f_2$ and $\{\mathbf{I-P_2}\}g_2$
are completely determined by already known functions. By
collision invariant property, we get for $m=1$
\begin{equation}
\begin{split}
\langle \partial_t f_1+v\cdot\nabla_x f_2+\frac{1}{\sqrt{\mu}}
E_1\cdot\nabla_v(\sqrt{\mu}g_1),[1,v,|v|^{2}/2]\sqrt{\mu }\rangle&=0,\\
\langle \partial_t g_1+v\cdot\nabla_x g_2+\frac{1}{\sqrt{\mu}}
E_1\cdot\nabla_v(\sqrt{\mu}f_1)-E_2\cdot v\sqrt{\mu},
\sqrt{\mu}\rangle&=0.\label{second}\\
\end{split}
\end{equation}
The process is similar to the derivation from pure Boltzmann to
Navier-Stokes except that now we have new terms to deal with; for
instance, see \cite{BGL1,M} for the computation of $\langle
\partial_t f_1+v\cdot\nabla_x f_2,[1,v,|v|^{2}/2]\sqrt{\mu
}\rangle$. Here we pay attention to those extras such as electric
field related terms. Note that magnetic fields are not involved
at this stage. The first equation of (\ref{second}) is equivalent
to
\begin{equation*}
\begin{split}
&\partial_t\rho_1+\nabla\cdot u_2+\langle\frac{1}{\sqrt{\mu}}
E_1\cdot\nabla_v(\sqrt{\mu}g_1),
\sqrt{\mu }\rangle=0,\\
&\partial_t u_1+u_1 \cdot\nabla u_1+\nabla p_1-\eta\Delta u_1 +
\langle\frac{1}{\sqrt{\mu}} E_1\cdot\nabla_v(\sqrt{\mu}g_1),
v\sqrt{\mu }\rangle=0,\\
&\frac{3}{2}\partial_{t}\{\rho_{1}+\theta_{1}\}+\frac{5}{2}u_1\cdot
\nabla\theta_1-\frac{5}{2}\kappa\Delta\theta_1+\frac{5}{2}
\nabla\cdot u_{2}+\langle\frac{1}{\sqrt{\mu}}
E_1\cdot\nabla_v(\sqrt{\mu}g_1), \frac{|v|^2}{2}\sqrt{\mu
}\rangle =0.\\
\end{split}
\end{equation*}
Electric field related terms can be taken care of by integrating
by parts in $v$:
\begin{align*}
&\langle\frac{1}{\sqrt{\mu}} E_1\cdot\nabla_v(\sqrt{\mu}g_1),
\sqrt{\mu }\rangle = E_1 \cdot\int \nabla_v(\sigma_1\mu)dv =0,\\
&\langle\frac{1}{\sqrt{\mu}} E_1\cdot\nabla_v(\sqrt{\mu}g_1),
v\sqrt{\mu }\rangle = E_1 \cdot\int \nabla_v(\sigma_1\mu) v
dv =-E_1\int\sigma_1 \mu dv=-E_1\sigma_1,\\
&\langle\frac{1}{\sqrt{\mu}} E_1\cdot\nabla_v(\sqrt{\mu}g_1),
\frac{|v|^2}{2}\sqrt{\mu }\rangle = E_1 \cdot\int
\nabla_v(\sigma_1\mu)\frac{|v|^2}{2}dv=-E_1\cdot(\sigma_1\int
v\sqrt{\mu}dv)=0.
\end{align*}
Thus the first equation in (\ref{second}) leads to (\ref{hincon})
for $m=2$ as well as (\ref{navier-stokes}) and (\ref{fourier}).
Before going any further, let us integrate (\ref{second}) over
${\bf T}^3$ and then, based on the above computation, we get
\[
\frac{d}{dt}\int_{\mathbb{T}^3}u_1(t,x)dx
=\frac{d}{dt}\int_{\mathbb{T}^3}\theta_1(t,x)dx
=\frac{d}{dt}\int_{\mathbb{T}^3}\sigma_1(t,x)dx=0
\]
which assures the validity of (\ref{conservation}) for $m=1$. For
the LHS of the second equation in (\ref{second}), since
\[
\langle\frac{1}{\sqrt{\mu}}
E_1\cdot\nabla_v(\sqrt{\mu}f_1)-E_2\cdot v\sqrt{\mu},
\sqrt{\mu}\rangle=0\text{ and }\langle v\cdot\nabla_x
\mathbf{P_2}g_2, \sqrt{\mu} \rangle=0,
\]
first it is reduced to
\begin{align*}
&\langle \partial_t g_1+v\cdot\nabla_x g_2+\frac{1}{\sqrt{\mu}}
E_1\cdot\nabla_v(\sqrt{\mu}f_1)-E_2\cdot v\sqrt{\mu},
\sqrt{\mu}\rangle\\
&=\partial_t\sigma_1+ \langle v\cdot\nabla_x\{\mathbf{I-P_2}\}g_2,
\sqrt{\mu} \rangle.
\end{align*}
By (\ref{mi1}), the second term can be written as,
\begin{equation}
\begin{split}
\langle v\cdot\nabla_x\{\mathbf{I-P_2}\}g_2, \sqrt{\mu}
\rangle&=\langle
v\cdot\nabla_x\{\mathcal{L}^{-1}(v\sqrt{\mu})\cdot
(-\nabla_x\sigma_1+E_1)+\sigma_1
f_1\}, \sqrt{\mu} \rangle\\
\
&=-\alpha\Delta\sigma_1+\alpha\sigma_1+\nabla\cdot(\sigma_1u_1)\\
\ &=-\alpha\Delta\sigma_1+\alpha\sigma_1+u_1\cdot\nabla
\sigma_1,\label{vestimate}
\end{split}
\end{equation}
where
\[
\alpha=\int_{\mathbb{R}^3}\mathcal{L}^{-1}(v\sqrt{\mu})\cdot
v\sqrt{\mu}\;dv>0.
\]
Hence we obtain (\ref{vlasov}) and it completes the case $m=1$.

Next we move onto higher order systems; consider (\ref{key}) for
$m\geq 2$. We shall use an induction on $m$. First by collision
invariant property, we get the following:
\begin{equation}
\begin{split}
\langle\partial_t f_{m}+v\cdot\nabla_x
f_{m+1}+\frac{1}{\sqrt{\mu}}\sum_{\substack{i+j=m+1\\i,j\geq1}}(E_i+v\times
B_i)\cdot\nabla_v(\sqrt{\mu}g_j),\;[1,v,\frac{|v|^2}{2}]\sqrt{\mu}\rangle=0,&\\
\langle\partial_t g_{m}+v\cdot\nabla_x
g_{m+1}+\frac{1}{\sqrt{\mu}}\sum_{\substack{i+j=m+1\\i,j\geq1}}(E_i+v\times
B_i)\cdot\nabla_v(\sqrt{\mu}f_j)-E_{m+1}\cdot
v\sqrt{\mu},\sqrt{\mu} \rangle\;&\\
=0,& \label{mtheq}
\end{split}
\end{equation}
where $\langle\partial_t f_{m}+v\cdot\nabla_x
f_{m+1},[1,v,|v|^2/2]\sqrt{\mu}\rangle$ was computed for pure
Boltzmann case in \cite{Guo2} (p.21-25). Here we will compute
electromagnetic field related  terms and then combine them with
the result in \cite{Guo2}.

Notice that for each $i,j\geq 1$,
$$\langle\frac{1}{\sqrt{\mu}}(E_i+v\times
B_i)\cdot\nabla_v(\sqrt{\mu}g_j),\;\sqrt{\mu}\rangle =0,$$ since
each integrand is a perfect derivative in $v$. Hence the very
first equation of (\ref{mtheq}) yields the $(m+1)$-th order
incompressibility condition:
\begin{equation}
\partial_t\rho_m+\nabla\cdot u_{m+1}=0.\label{m+1incom}
\end{equation}
To derive the velocity equation, we look at the next
electromagnetic field term:
\begin{equation*}
\begin{split}
&\sum_{\substack{i+j=m+1\\i,j\geq1}}\langle\frac{1}{\sqrt{\mu}}(E_i+v\times
B_i)\cdot\nabla_v(\sqrt{\mu}g_j),\;v\sqrt{\mu}\rangle\\
&=-\sum_{\substack{i+j=m+1\\i,j\geq1}}\int(E_i+v\times
B_i)\;g_j\sqrt{\mu} dv\;\;(\text{by the integration by parts})\\
&=-\sum_{\substack{i+j=m+1\\i,j\geq1}}\{E_i\nabla\cdot
E_j-(\partial_t E_{j-1}-\nabla\times B_j)\times B_i\}\; (\text{by
(\ref{keyeb})})
\end{split}
\end{equation*}
So the second part of the first equation in (\ref{mtheq}) with
(4.8) in \cite{Guo2} is equivalent to
\begin{equation}
\begin{split}
& \partial_{t}u_{m}+\langle v\cdot\nabla_{x}L^{-1}\{-\partial_{t}%
f_{m-1}-v\cdot\nabla_{x}f_{m}-\frac
{1}{\sqrt{\mu}}\sum_{\substack{i+j=m\\i,j\geq1}}(E_i+v\times B_i)
\cdot\nabla_v(\sqrt{\mu}g_j)\\
&\;\;+\sum_{i+j=m+1}\Gamma(f_{i},f_{j})\},\;v\sqrt{\mu
}\rangle-\sum_{\substack{i+j=m+1\\i,j\geq1}}\{E_i\nabla\cdot
E_j-(\partial_t E_{j-1}-\nabla\times B_j)\times
B_i\}\\
&=-\nabla_{x}(\rho_{m+1}+ \theta_{m+1}),\label{velocity}
\end{split}
\end{equation}
where we have applied (\ref{hmicro1}) to solve for the
microscopic part $(\mathbf{I-P_1})f_{m+1}$. Any term of which
index is lower than $m$ in the above makes a contribution to the
remainder $R_m^u$. In particular, all other electromagnetic terms
are to be included in the remainder except for $E_m\nabla\cdot
E_1$ and $E_1\nabla\cdot E_m$. For the estimate of the rest, we
refer \cite{Guo2} (p.21-25). As splitting
$u_{m}=P_{0}u_{m}+(I-P_{0})u_{m}$, one can readily keep track of
each term in (\ref{hnavier-stokes}), (\ref{ru}) and
(\ref{press}). In addition, to obtain the $m$-th order Boussinesq
relation (\ref{hbous}), we use the pressure $p_{m-1}$ in
(\ref{press}) to solve for $\rho_{m}+\theta_{m}$. Notice that by
taking divergence of the $(m-1)$-th order Navier-Stokes equation
(\ref{hnavier-stokes}), we have another expression for $p_{m-1}$
as
\[
p_{m-1}=\Delta^{-1}\nabla\cdot
\{-u_{1}\cdot\nabla(P_{0}u_{m-1})-P_{0}u_{m-1}\cdot\nabla
u_{1}+E_1\sigma_m+E_m\sigma_1+R_{m-1}^{u}\}.
\]
For the temperature equation, we need to compute the following
electromagnetic term:
\begin{equation*}
\begin{split}
&\sum_{\substack{i+j=m+1\\i,j\geq1}}\langle\frac{1}{\sqrt{\mu}}(E_i+v\times
B_i)\cdot\nabla_v(\sqrt{\mu}g_j),\;\frac{|v|^2}{2}\sqrt{\mu}
\rangle\\
&=-\sum_{\substack{i+j=m+1\\i,j\geq1}}\int(E_i+v\times B_i)\cdot
v g_j\sqrt{\mu}dv\;\; (\text{by the integration by parts})\\
&=-\sum_{\substack{i+j=m+1\\i,j\geq1}}E_i\cdot\int v
g_j\sqrt{\mu}dv\;\; (\text{since $V\times W\cdot V=0$ for $V, W$
vectors in $\bf{R}^3$})\\
&=\sum_{\substack{i+j=m+1\\i,j\geq1}}E_i\cdot(\partial_t
E_{j-1}-\nabla\times B_j)\;\;(\text{by (\ref{keyeb})})
\end{split}
\end{equation*}
Then the last part of the first equation in (\ref{mtheq}) by using
(4.8) in \cite{Guo2} becomes
\begin{equation}
\begin{split}
& \frac{5}{2}\partial_{t}\theta_{m}+\langle v\cdot\nabla_{x}L^{-1}%
\{-\partial_{t}f_{m-1}-v\cdot\nabla_{x}f_{m}-\frac
{1}{\sqrt{\mu}}\sum_{\substack{i+j=m\\i,j\geq1}}(E_i+v\times B_i)
\cdot\nabla_v(\sqrt{\mu}g_j)\label{temperature}\\
&\;\;+\sum_{i+j=m+1}\Gamma(f_{i}%
,f_{j})\},\frac{|v|^{2}}{2}\sqrt{\mu}\rangle+
\sum_{\substack{i+j=m+1\\i,j\geq1}}E_i\cdot(\partial_t
E_{j-1}-\nabla\times B_j)\\
&=\partial_{t}\{\rho_{m}+\theta_{m}\},
\end{split}
\end{equation}
where (\ref{m+1incom}) has been used. One can easily deduce
(\ref{hfourier1}) as well as (\ref{rs}).

Now we come to the second equation in (\ref{mtheq}). It is
equivalent to
\begin{equation}
\partial_t\sigma_m+\langle v\cdot\nabla_x(\mathbf{I-P_2})g_{m+1},
\sqrt{\mu} \rangle=0,\label{varmth}
\end{equation}
because $\langle
\frac{1}{\sqrt{\mu}}\sum_{\substack{i+j=m+1\\i,j\geq1}}(E_i+v\times
B_i)\cdot\nabla_v(\sqrt{\mu}f_j)-E_{m+1}\cdot
v\sqrt{\mu},\sqrt{\mu} \rangle=0\;$ and %
$\langle v \cdot\nabla_x \mathbf{P}_2 g_{m+1},
\sqrt{\mu}\rangle=0$. Solving for the microscopic part
$(\mathbf{I-P_2})g_{m+1}$ by (\ref{hmicro2}) and plugging it into
(\ref{varmth}), we deduce that
\begin{equation}
\begin{split}
 \partial_{t}\sigma_{m}+\langle
v\cdot\nabla_{x}\mathcal{L}^{-1}\{-(\mathbf{I-P_2})&(v\cdot\nabla_{x}\mathbf{P}_2
g_{m})+E_m\cdot v\sqrt{\mu}+\Gamma(g_1,\mathbf{P}_1
f_m)\\
+\Gamma(\mathbf{P}_2 g_m,& f_1)\},\sqrt{\mu}\rangle \\
\ =\langle
v\cdot\nabla_{x}\mathcal{L}^{-1}\{\partial_t(\mathbf{I-P_2})
g_{m-1}&+(\mathbf{I-P_2})(v\cdot\nabla_x(\mathbf{I-P_2})g_m)-\Gamma
(g_1,(\mathbf{I-P_1})f_m)\\
-\Gamma((\mathbf{I-P}_2) g_m,& f_1) +\frac
{1}{\sqrt{\mu}}\sum_{\substack{i+j=m\\i,j\geq1}}(E_i+v\times B_i)
\cdot\nabla_v(\sqrt{\mu}f_j)\\
+\sum_{\substack{i+j=m+1\\i,j\geq
2}}\Gamma(g_{i},&f_{j})\},\;v\sqrt{\mu }\rangle.\label{hsigma}
\end{split}
\end{equation}
Note that each term in the RHS of (\ref{hsigma}) essentially has
an index of lower than $m$ and hence is part of (\ref{rv}). Terms
in the LHS can be computed in the same way as in
(\ref{vestimate}). Recalling $\mathbf{P}_2
g_m=\sigma_m\sqrt{\mu}$, we get
\begin{align*}
&\langle
v\cdot\nabla_{x}\mathcal{L}^{-1}\{(\mathbf{I-P_2})(v\cdot\nabla_{x}\mathbf{P}_2
g_{m})\},\sqrt{\mu}\rangle =\alpha\Delta\sigma_m,\\
\ &\langle v\cdot\nabla_{x}\mathcal{L}^{-1}\{E_m\cdot
v\sqrt{\mu}\},\sqrt{\mu}\rangle =\alpha\sigma_m,\\
\ &\langle
v\cdot\nabla_{x}\mathcal{L}^{-1}\{\Gamma(g_1,\mathbf{P}_1 f_m)
+\Gamma(\mathbf{P}_2 g_m,f_1)\},\sqrt{\mu}\rangle
=\nabla\cdot(\sigma_1 u_m)+\nabla\cdot(\sigma_m u_1),
\end{align*}
where we have used the fact $\Gamma(\sqrt{\mu},\mathbf{P}_1
f_m)=\mathcal{L}\mathbf{P}_1 f_m$ and
$\Gamma(\sqrt{\mu},f_1)=\mathcal{L}f_1$ which follow from the
following observation: for any $h$ with $Lh=0$,
\[
\mathcal{L}h=-\Gamma(h,\sqrt{\mu})=\Gamma(\sqrt{\mu},h),\; \text{
since }-Lh=\Gamma(h,\sqrt{\mu})+\Gamma(\sqrt{\mu},h).
\]
Therefore, (\ref{hvlasov}) can be derived with the remainder
$R_m^{\sigma}$ (\ref{rv}).

As for field equations,  (\ref{keyeb}) with the definition of
$\mathbf{P_2}g_m$ immediately leads to (\ref{helec}) and
(\ref{hmag}). To get the compatibility conditions
(\ref{compaEB}), take the integral in $x$ of field equations in
(\ref{keyeb}) and then the periodic boundary conditions give rise
to
\begin{equation*}
\begin{split}
\frac{d}{dt}\int_{\mathbb{T}^3}B_m dx=0,\;\;\;\;
\frac{d}{dt}\int_{\mathbb{T}^3}E_m dx=-\int_{\mathbb{T}^3}
\int_{\mathbb{R}^3} v g_{m+1}\sqrt{\mu}\;dv\; dx.
\end{split}
\end{equation*}
Thus $\int_{\mathbb{T}^3}B_m dx=0$ may be assumed for all time.
The dynamics of $\int_{\mathbb{T}^3}E_m dx$ is more complicated.
First note that $\int\int\mathbf{P_2}g_m v\sqrt{\mu}dvdx=0.$ By
using (\ref{hmicro2}), we have
\begin{equation*}
\begin{split}
\frac{d}{dt}\int_{\mathbb{T}^3}E_m dx
=-\int_{\mathbb{T}^3} \int_{\mathbb{R}^3}(\mathbf{I-P_2})g&_{m+1}\;v\sqrt{\mu}\;dv dx\\
=-\int_{\mathbb{T}^3} \int_{\mathbb{R}^3}
v\sqrt{\mu}\;\mathcal{L}^{-1}&\{-\partial_{t}g_{m-1}
-v\cdot\nabla_{x}
g_{m}+\sum_{\substack{i+j=m+1\\i,j\geq1}}\Gamma(g_{i},f_{j})\\
-\frac
{1}{\sqrt{\mu}}\sum_{\substack{i+j=m\\i,j\geq1}}&(E_i+v\times B_i)
\cdot\nabla_v(\sqrt{\mu}f_j)+E_{m}\cdot v\sqrt{\mu}\}dv dx\\
=\alpha\int_{\mathbb{T}^3} E_m\;dx+\ell_{m-1}&,
\end{split}
\end{equation*}
where
\begin{equation}
\begin{split}
\ell_{m-1}\equiv &\int_{\mathbb{T}^3}
\int_{\mathbb{R}^3}\partial_t
g_{m-1}\mathcal{L}^{-1}(v\sqrt{\mu})dvdx\\
&+\sum_{\substack{i+j=m\\i,j\geq1}}\int_{\mathbb{T}^3}
\int_{\mathbb{R}^3}\frac{1}{\sqrt{\mu}} (E_i+v\times
B_i)\cdot\nabla_v(\sqrt{\mu}f_j)
\mathcal{L}^{-1}(v\sqrt{\mu})dvdx.\label{ell}
\end{split}
\end{equation}
Note that $\ell_{m-1}$ only contains terms of index lower than
$m$ and $\int E_mdx$ at $t=0$ for $m\geq 2$ can be arbitrarily
given. Other terms have vanished after the integration because of
either a perfect derivative or collision invariant property.

It remains to verify conservation laws (\ref{conservation}) for
$m\geq 2$ to finish Lemma. Firstly, it is easy to see
\[
\frac{d}{dt}\int_{\mathbb{T}^3}\rho_m dx=0\;\text{ and }
\;\frac{d}{dt}\int_{\mathbb{T}^3}\sigma_mdx=0,
\]
by integrating (\ref{m+1incom}) and (\ref{varmth}) over
$\mathbb{T}^3$. We give the detailed computation for the average
of the velocity and the similar argument can be applied to the
temperature field. Integrate (\ref{velocity}) to get
\begin{equation}
\frac{d}{dt}\int_{\mathbb{T}^3}u_m dx=
\sum_{\substack{i+j=m+1\\i,j\geq1}}
\int_{\mathbb{T}^3}\{E_i\nabla\cdot E_j-(\partial_t
E_{j-1}-\nabla\times B_j)\times B_i\}\;dx.\label{avev}
\end{equation}
To simplify the RHS, we observe the following:
\begin{align*}
 \int_{\mathbb{T}^3}E_i(\nabla\cdot E_j) +E_j(\nabla\cdot E_i) dx
&=-\int_{\mathbb{T}^3}(\nabla\times E_i)\times E_j +(\nabla\times
E_j)\times E_i\; dx\\
&=\int_{\mathbb{T}^3}\partial_t B_{i-1}\times E_j +\partial_t
B_{j-1}
\times E_i\; dx,\\
 \int_{\mathbb{T}^3}(\nabla\times B_j)\times B_i +(\nabla\times
B_i)&\times B_j\; dx=0\;\;\;(\text{since }\nabla\cdot
B_i=\nabla\cdot
 B_j=0),
\end{align*}
where we have integrated by parts. Thus the RHS of (\ref{avev}) is
equivalent to
\[
\sum_{\substack{i+j=m+1\\i,j\geq1}} \int_{\mathbb{T}^3}\partial_t
B_{i-1}\times E_j-\partial_t E_{j-1}\times B_i \;dx
\]
and in turn by rearranging indices we deduce
\[
\frac{d}{dt}\int_{\mathbb{T}^3}u_m
dx=-\frac{d}{dt}\sum_{\substack{i+j=m\\i,j\geq1}}
\int_{\mathbb{T}^3}E_i\times B_{j}\;dx.
\]
Moreover, $\int_{\mathbb{T}^3}u_m dx=\int_{\mathbb{T}^3}P_0 u_m
dx$, since $\{I-P_0\}u_m=\nabla q_m$ for some scalar function
$q_m$ by the Helmholtz-Hodge decomposition. Therefore, the
desired result follows.
\end{proof}\

\section{The Diffusive Coefficients}\

We now show that we can determine the coefficients
$f_{1},f_{2},...f_{m};g_{1},g_{2},...g_{m}$ and
$E_1,...E_m;B_1,...B_m$ by giving the initial divergent free part
of velocity field $u_{i}^{0}(x)$, the temperature field
$\theta_{i}^{0}(x)$, the concentration difference field
$\sigma_i^0(x)$ and the average value $e_i$ of initial electric
field $E_i^0$ with $e_1=0$.

\begin{proof}
(\textit{of Theorem \ref{cecoe}:}) We use the induction over $r.$
First consider the case  $r=1.$ We need to solve the system
(\ref{navier-stokes})-(\ref{fourier}) with (\ref{incon}) and
$\int u_{1}=\int\theta_{1}=\int\sigma_1=0.$ Once  a priori
estimate is given, the existence and the uniqueness of a solution
$[u_{1}(t,x),\theta_{1}(t,x), \sigma_1(t,x),\phi_1(t,x)]$ follow
from the standard iteration argument via the fixed point theorem.
Thus it suffices to illustrate a priori energy estimate in the
subsequence. First we take $\partial_{\gamma}$ derivatives of
(\ref{navier-stokes})-(\ref{fourier}):
\begin{equation}
\begin{split}
&\partial_{t}\partial_{\gamma}u_{1}+u_{1}\nabla\partial_{\gamma}%
u_{1}+\nabla\partial_{\gamma}p_{1}-\eta\Delta\partial_{\gamma}u_{1}
=-\partial_{\gamma_{1}}u_{1}\nabla%
\partial_{\gamma_{2}}u_{1}+\partial_{\gamma}(\sigma_1\nabla\phi_1),\\
&\partial_{t}\partial_{\gamma}\sigma_{1}+u_{1}\nabla\partial_{\gamma}
\sigma_{1}-\alpha\Delta\partial_{\gamma}\sigma_{1}+\alpha
\partial_{\gamma}\sigma_1
=-\partial_{\gamma_{1}}u_{1}\nabla\partial_{\gamma_{2}}\sigma_{1},\\
&\partial_{t}\partial_{\gamma}\theta_{1}+u_{1}\nabla\partial_{\gamma}\theta
_{1}-\kappa\Delta\partial_{\gamma}\theta_{1}
=-\partial_{\gamma_{1}}u_{1}\nabla\partial_{\gamma_{2}}\theta_{1},
\label{gamma}\\
&\Delta\partial_{\gamma}\phi_1=\partial_{\gamma}\sigma_1,\;\;\nabla\cdot
u_1=0.
\end{split}
\end{equation}
The last three summations are over $\gamma_{1}+\gamma_{2}=\gamma$
with $|\gamma_{1}|\geq1$. Let $\gamma=0$. Multiply first three
equation in (\ref{gamma}) by $[u_1, \sigma_1,\theta_1]$,
integrate over $\mathbb{T}^3$ and by the incompressibility we get
\begin{equation*}
\begin{split}
\frac{1}{2}&\frac{d}{dt}[\int|u_1|^2+|\sigma_1|^2+|\theta_1|^2dx]+\eta\int|\nabla
u_1|^2dx+\alpha
\int|\nabla\sigma_1|^2dx+\kappa\int|\nabla\theta_1|^2dx\\
&+\alpha\int|\sigma_1|^2dx =\int\sigma_1\nabla\phi_1\cdot u_1 dx.
\end{split}
\end{equation*}
Note that the RHS can be absorbed into the LHS assuming $M$ is
sufficiently small at time $t$, since
\[
\int\sigma_1\nabla\phi_1\cdot u_1 dx \leq
(\sup_x|u_1|)||\sigma_1||\cdot||\nabla\phi_1||\leq
CM||\sigma_1||^2,
\]
where we have used the Sobolev embedding theorem and the $L^2$
estimate for the poisson equation $\Delta\phi_1=\sigma_1$.
Similarly, the energy estimate for $|\gamma|\leq 2$ shows
\begin{align*}
\frac{1}{2}\frac{d}{dt}[\int|\partial_{\gamma}
u_1|^2+|\partial_{\gamma}\sigma_1|^2+&|\partial_{\gamma}\theta_1|^2dx]
+\eta\int|\nabla
\partial_{\gamma} u_1|^2dx+\alpha
\int|\nabla\partial_{\gamma} \sigma_1|^2dx\\
+\kappa\int|\nabla\partial_{\gamma} \theta_1|^2dx
+&\alpha\int|\partial_{\gamma}\sigma_1|^2dx\\
\leq\sup_x(|u_1|+|\theta_1|+|\sigma_1|)(||\partial&_{\gamma} u_1
||^2+||\nabla\partial_{\gamma} u_1
||^2+||\partial_{\gamma}\sigma_1
||^2+||\nabla\partial_{\gamma}\sigma_1 ||^2+||\partial_{\gamma}
\theta_1 ||^2\\&+||\nabla\partial_{\gamma} \theta_1
||^2+||\nabla\partial_{\gamma}\phi_1 ||^2).
\end{align*}
Each term in the RHS can be absorbed into the LHS such that
\[
\frac{1}{2}\frac{d}{dt}(||u_1||_{H^2}^2+||\sigma_1||_{H^2}^2+
||\theta_1||_{H^2}^2)+\frac{\eta}{2}||\nabla u_1||_{H^2}^2+
\frac{\alpha}{2}|| \sigma_1||_{H^3}^2+\frac{\kappa}{2}||\nabla
\theta_1||_{H^2}^2\leq 0.
\]
Applying the Poincar$\acute{e}$ inequality, for
$\lambda=\frac{1}{4}\min\{\eta,\alpha,\kappa\}>0$ we have
\[
||u_{1}(t)||_{H^{2}}+||\sigma_{1}(t)||_{H^{2}}+
||\theta_{1}(t)||_{H^{2}}\leq Ce^{-\lambda t}%
\{||u_{1}^{0}||_{H^{2}}+||\sigma_{1}^0||_{H^{2}}+
||\theta_{1}^{0}||_{H^{2}}\},
\]
which immediately implies the existence, uniqueness and
exponential decay of
 $u_1(t,x),$ $\sigma_1(t,x)$ and $\theta_1(t,x)$.

Next we turn to high order derivative cases. We claim that for
$s\geq3,$ there exists a polynomial $U_{s}$ $\geq0$ with
$U_{s}(0)=0$ such that
\begin{equation}
||u_{1}(t)||_{H^{s}}+||\theta_{1}(t)||_{H^{s}}+||\sigma_{1}(t)||_{H^{s}}
\leq e^{-\lambda t}U_{s}%
(||u_{1}^{0}||_{H^{s}}+||\theta_{1}^{0}||_{H^{s}}+
||\sigma_{1}^{0}||_{H^{s}}). \label{claim}
\end{equation}
Separating the case of $|\gamma_{1}|=1$ or $|\gamma_{1}|=s,$ and
the case of $|\gamma_{1}|\leq s-1$, $|\gamma_{2}|\leq s-2$,  we
estimate the $L^{2}$ norm of the RHS' in (\ref{gamma}) by
\begin{align*}
\{C(||&u_{1}||_{H^{2}}+||\sigma_1||_{H^2})+\xi\}
\{\sum_{|\gamma|=s+1}(||\partial_{\gamma}%
u_{1}||+||\partial_{\gamma}\sigma_{1}||+
||\partial_{\gamma}\theta_{1}||)+\sum_{|\gamma|=s}
||\partial_{\gamma}\nabla\phi_1||\}\\
&+C_{\xi}\{||u_{1}||_{H^{s-1}%
}+||\sigma_{1}||_{H^{s-1}}+||\theta_{1}||_{H^{s-1}}+
||\nabla\phi_{1}||_{H^{s-1}}\}^{2},%
\end{align*}
for any small number $\xi>0,$ by an interpolation in the Sobolev
Space. We then use standard energy estimate for
$\partial_{\gamma}u_{1}$, $\partial_{\gamma}\sigma_{1}$  and
$\partial_{\gamma}\theta_{1}$ to get ($\nabla\cdot u_{1}=0$)%
\begin{align*}
 \frac{1}{2}\frac{d}{dt}\sum_{|\gamma|=s}\{||\partial_{\gamma}u_{1}%
||^{2}+||\partial&_{\gamma}\sigma_{1}||^{2}+||\partial_{\gamma}
\theta_{1}||^{2}\}\\+
 \sum_{|\gamma
|=s}&\{\eta||\nabla\partial_{\gamma}u_{1}||^{2}+\kappa||\nabla\partial_{\gamma
}\theta_{1}||^{2}+\alpha ||\partial_{\gamma}\sigma_1||^2
+\alpha||\nabla
\partial_{\gamma}\sigma_1||^2\}\leq\\
\{C(||u_{1}||_{H^{2}}+||\sigma_1||_{H^2})
+&\xi\}\{\sum_{|\gamma|=s+1}||\partial_{\gamma}u_{1}||+||\partial_{\gamma}\sigma_{1}||+
||\partial_{\gamma}\theta_{1}||+\sum_{|\gamma|=s}||\nabla
\partial_{\gamma}\phi_1||\}^2\\
&+C_{\xi}\{||u_{1}||_{H^{s-1}%
}+||\sigma_{1}||_{H^{s-1}}+||\theta_{1}||_{H^{s-1}}+
||\nabla\phi_{1}||_{H^{s-1}}\}^{4}.\\
\end{align*}
Note that $||\nabla\phi_{1}||_{H^{s-1}}\leq C
||\sigma_{1}||_{H^{s-1}}$. For
$C(||u_{1}||_{H^{2}}+||\sigma_{1}||_{H^{2}})$ and $\xi$
sufficiently small, we deduce from the Poincar$\acute{e}$
inequality and an induction for $s-1$ over the last lower order
term that
\begin{align*}
\sum_{|\gamma|\leq
s}||\partial_{\gamma}u_{1}(t)||+||\partial_{\gamma}
\sigma_{1}(t)||+ ||\partial_{\gamma} \theta_{1}(t)||\;\;&\\
 \leq e^{-\lambda t}\sum_{|\gamma|\leq s}||\partial
_{\gamma}u_{1}(0)||+||\partial_{\gamma}\sigma_{1}(0)||
+||&\partial_{\gamma}\theta_{1}(0)||\\
+\int_{0}^{t}e^{-\lambda\{t-\tau\}}e^{-2\lambda\tau}\{U_{s-1}
(||u_{1}^{0}||_{H^{s-1}}&+||\sigma_{1}^{0}||_{H^{s-1}}+
||\theta_{1}^{0}||_{H^{s-1}})\}^{2}d\tau\\
 \leq e^{-\lambda
t}U_{s}(||u_{1}(0)||_{H^{s}}+||\sigma_{1}(0)||_{H^{s}}
+||\theta&_{1}(0)||_{H^{s}}).
\end{align*}
We thus conclude our claim (\ref{claim}).

Now let us turn to general space-time derivatives
$\partial_{\tau}$. Notice that
\[
\Delta p_{1}=\nabla\cdot\{\sigma_1\nabla\phi_1 -u_{1}\cdot\nabla
u_{1}\}.
\]
We then use repeatedly the equations (\ref{navier-stokes}),
(\ref{vlasov}) and (\ref{fourier}) to solve for temporal
derivatives to get
\begin{align*}
\sum_{|\tau|\leq
s}\{||\partial_{\tau}u_{1}(t)||+||\partial_{\tau}\sigma
_{1}(t)||+||\partial_{\tau}\theta _{1}(t)||\}\leq e^{-\lambda
t}U_{2s}(||u_{1}(0)||_{H^{2s}}&+||\sigma
_{1}(0)||_{H^{2s}}\\
&+||\theta _{1}(0)||_{H^{2s}}).
\end{align*}
Notice that we need twice many spatial derivatives now for the
initial data. Finally, since
\[
f_{1}\equiv v\cdot
u_{1}\sqrt{\mu}+\{\frac{|v|^{2}}{2}-\frac{5}{2}\}\theta
_{1}\sqrt{\mu},\;\;\text{and }\;g_1\equiv\sigma_1\sqrt{\mu},
\]
it follows that for any $s\geq0$ and $\beta,$
\begin{align*}
\sum_{|\tau|\leq
s}||[\partial_{\tau}^{\beta}f_{1}(t),\partial_{\tau}^{\beta}g_1(t)]||_{\nu}
& \leq C\sum_{|\tau|\leq
s}\{||\partial_{\tau}u_{1}(t)||+||\partial_{\tau}\sigma
_{1}(t)||+||\partial_{\tau}\theta
_{1}(t)||\}\\
&  \leq e^{-\lambda t}U_{2s}(||u_{1}^{0}||_{H^{2s}}+
||\sigma_{1}^{0}||_{H^{2s}}+||\theta_{1}^{0}%
||_{H^{2s}}).
\end{align*}
Recall that $E_1=\nabla\phi_1$ with $\Delta\phi_1=\sigma_1$ and
$B_1=0$. Our theorem  thus is valid for $r=1.$

Assume that $f_{1},f_{2},...f_{r}$, $g_1,g_2,...g_r$,
$E_1,E_2,...E_r$, and $B_1(=0), B_2,...B_r$ have been constructed
to satisfy (\ref{hmicro1})-(\ref{compaEB}) for up to $r\geq1$. We
now construct $f_{r+1}$, $g_{r+1}$ in two steps.\\

\textbf{Step One} : \textit{Construct the microscopic part
$\{\mathbf{I-P_1}\}f_{r+1}$ and $\{\mathbf{I-P_2}\}g_{r+1}$ from
the microscopic equation (\ref{key}).}

In order to solve for $\{\mathbf{I-P_1}\}f_{r+1}$ and
$\{\mathbf{I-P_2}\}g_{r+1}$, we need (\ref{mtheq}) for $m=r-1$,
which are equivalent to (\ref{m+1incom}), (\ref{velocity}),
(\ref{temperature}) and (\ref{hsigma}) for $m=r-1,$ to be true.
By the induction hypothesis, the $r$-th incompressibility
condition (\ref{hincon}) and the ($r-1$)-th temperature equation
(\ref{hfourier1}) imply that (\ref{m+1incom}) and
(\ref{temperature}) hold for $m=r-1$. And ($r-1$)-th concentration
difference equation (\ref{hvlasov}) leads to (\ref{hsigma}) for
$m=r-1$. Finally, applying the $r$-th order Boussinesq relation
(\ref{hbous}) into the ($r-1)$-th order Vlasov-Navier-Stokes
system (\ref{hnavier-stokes}), we conclude that the pressure
$p_{r-1}$ is given by (\ref{press}) for $m=r-1,$ which also
implies that (\ref{velocity}) is valid for $r-1.$ Therefore, we
can solve $\{\mathbf{I-P_1}\}f_{r+1}$ and
$\{\mathbf{I-P_2}\}g_{r+1}$ from (\ref{key}) such that
\begin{equation*}
\begin{split}
L\{\mathbf{I-P_1}\}f_{r+1}=-\partial_{t}f_{r-1}-v\cdot\nabla_{x}f_{r}-
\frac{1}{\sqrt{\mu}}
\sum_{\substack{i+j=r\\i,j\geq1}}(E_i+v\times B_i)
\cdot\nabla_v(\sqrt{\mu}g_j)\\
 +\sum_{\substack{i+j=r+1\\i,j\geq1}}\Gamma(f_{i},f_{j});\\
\mathcal{L}\{\mathbf{I-P_2}\}g_{r+1}=-\partial_{t}g_{r-1}-
v\cdot\nabla_{x}g_{r}- \frac{1}{\sqrt{\mu}}
\sum_{\substack{i+j=r\\i,j\geq1}}(E_i+v\times B_i)
\cdot\nabla_v(\sqrt{\mu}f_j)\\
+E_r\cdot
v\sqrt{\mu}+\sum_{\substack{i+j=r+1\\i,j\geq1}}\Gamma(g_{i},f_{j}).
\end{split}
\end{equation*}
We next estimate such microscopic parts. Since two expressions in
the above have almost same structure, we compute only one term
$\{\mathbf{I-P_2}\}g_{r+1}$. The same argument can be applied to
$\{\mathbf{I-P_1}\}f_{r+1}$. Taking $\partial_{\tau}^{\beta}$
derivatives, we get
\begin{align*}
(\partial_{\tau}^{\beta}\mathcal{L}\{\mathbf{I-P_2}\}g_{r+1},
\partial_{\tau}^{\beta}\{\mathbf{I-P_2}\}g_{r+1})&\\
=(-\partial_{\tau}^{\beta}\partial_{t}g_{r-1}-
\partial_{\tau}^{\beta}\{v\cdot\nabla_{x}g_{r}\}-&\partial_{\tau}^{\beta}
\{\frac{1}{\sqrt{\mu}}
\sum_{\substack{i+j=r\\i,j\geq1}}(E_i+v\times B_i)
\cdot\nabla_v(\sqrt{\mu}f_j)\}\\
+ \partial_{\tau}^{\beta}\{E_r\cdot v\sqrt{\mu}\}
+&\partial_{\tau}^{\beta}\sum_{\substack{i+j=r+1\\i,j\geq1}}
\Gamma(g_{i},f_{j}),\;\partial_{\tau}^{\beta}\{\mathbf{I-P_2}\}g_{r+1}).
\end{align*}
Applying Lemma 3.3 in \cite{Guo2}, we first get
($C_{\beta}=0\text{ if }\beta=0$)
\begin{equation}
\begin{split}
\frac{1}{2}||\partial_{\tau}^{\beta}&
{\{\mathbf{I-P_2}\}g_{r+1}}||_{\nu}^{2}-C_{\beta} ||
{\{\mathbf{I-P_2}\}g_{r+1}}||_{\nu}^2\\
\leq&(-\partial_{t}\partial_{\tau}^{\beta}{g_{r-1}}-
\partial_{\tau}^{\beta}\{v\cdot\nabla_{x}{g_r}\}
+\partial_{\tau}^{\beta}\{E_r\cdot
v\sqrt{\mu}\},\;
\partial_{\tau}^{\beta
}{\{\mathbf{I-P_2}\}g_{r+1}})\\
&+(\partial_{\tau}^{\beta}\sum_{\substack{i+j=r+1\\i,j\geq1}}
{\Gamma(g_{i},f_{j})},\;
\partial_{\tau}^{\beta}{\{\mathbf{I-P_2}\}g_{r+1}})\\
&+(-\partial_{\tau}^{\beta} \{\frac{1}{\sqrt{\mu}}
\sum_{\substack{i+j=r\\i,j\geq1}}(E_i+v\times B_i) \cdot\nabla_v
{\sqrt{\mu}f_j}\},\;
\partial_{\tau}^{\beta
}{\{\mathbf{I-P_2}\}g_{r+1}}) \label{wr-1}\\
\equiv&(I)+(II)+(III).
\end{split}
\end{equation}
Since $|\partial^{\beta}\{v\sqrt{\mu}\}|$ is uniformly bounded,
by using (\ref{weight}), $(I)$ is bounded by
\[
\{||\partial_{t}\partial_{\tau}^{\beta}{g_{r-1}}||+||
\partial_{\tau}^{\beta}\nabla_{x}{g_r}||_{\nu}+
||\partial_{\tau}^{\beta
-\beta_{1}}\nabla_{x}{g_r}||+||\partial_{\tau}E_r||\}
\times||\partial_{\tau}^{\beta }\{\mathbf{I-P_2}\}g_{r+1}||_{\nu},
\]
where $|\beta_{1}|=1.$ Applying Lemma 3.3 in \cite{Guo2} again,
the nonlinear collision term $(II)$ is bounded by
\[
\sum_{\substack{i+j=r+1\\i,j\geq1}}\{||\partial_{\tau}^{\beta}
{g_i}||\cdot||\partial_{\tau}^{\beta}f_{j}||_{\nu}
+||\partial_{\tau
}^{\beta}{g_i}||_{\nu}\cdot||\partial_{\tau}^{\beta}f_{j}||\}
\times ||\partial_{\tau}^{\beta}
{\{\mathbf{I-P_2}\}g_{r+1}}||_{\nu}.
\]
To control the electromagnetic terms $(III)$, we introduce the
following inequality.\\

Claim. There exists $C>0$ such that for $|\tau|\leq s$ and
$\beta_1\leq \beta$,
\begin{equation}
\begin{split}
(\partial_{\tau}^{\beta}\{\frac{1}{\sqrt{\mu}}(E+v\times B
)&\cdot\nabla_v (\sqrt{\mu}f)\},\partial_{\tau}^{\beta}g)\\
\leq C
\sum_{|\widetilde{\tau}|\leq\max\{s,3\}}(&||\partial_{\widetilde{\tau}}E||\cdot
||\partial_{\widetilde{\tau}}^{\beta}\nabla_v
f||+||\partial_{\widetilde{\tau}}B||\cdot
||\partial_{\widetilde{\tau}}^{\beta_1}\nabla_v
f||_{\nu}\\
&+||\partial_{\widetilde{\tau}}E||\cdot
||\partial_{\widetilde{\tau}}^{\beta_1}f||)
\cdot||\partial_{\tau}^{\beta}g||_{\nu}.\label{eb1}\\
\\
\end{split}
\end{equation}
The direct result of the claim is that $(III)$ is bounded by
\[
\sum_{\substack{i+j=r\\i,j\geq 1}}\sum_{|\widetilde{\tau}|}\{||
[\partial_{\widetilde{\tau}}E_i,\partial_{\widetilde{\tau}}B_i]||\cdot||\partial_{\widetilde{\tau}}^{\beta}\nabla_v
{f_{j}}||_{\nu}+||\partial_{\widetilde{\tau}}E_i||\cdot||
\partial_{\widetilde{\tau}}^{\beta}{f_j}||\}
\times ||\partial_{\tau}^{\beta}
{\{\mathbf{I-P_2}\}g_{r+1}}||_{\nu},
\]
where $|\widetilde{\tau}| \leq\max\{|\tau|,3\}$. From the
induction hypothesis (\ref{rdecay}) for
$f_{1},...f_{r};g_1,...g_r$ and $E_1,...E_r;B_1,...B_r$,
consequently, we deduce that for $|\tau|\leq s,$ (there are $s+1$
derivatives for $f_{r},g_r$)
\begin{equation}
\begin{split}
||&[\partial_{\tau}^{\beta}\{\mathbf{I-P_1}\}f_{r+1},
\partial_{\tau}^{\beta}
\{\mathbf{I-P_2}\}g_{r+1}]||_{\nu}\\
\leq&\; C e^{-\lambda t}
U_{r+1}(\sum_{i=1}^{r}\{||u_{i}^{0}||_{H^{2s+2+4(r-i)}}
+||\theta_{i}^{0}||_{H^{2s+2+4(r-i)}}
+||\sigma_{i}^{0}||_{H^{2s+2+4(r-i)}}\}). \label{i-p}
\end{split}
\end{equation}
To complete (\ref{i-p}), it remains to prove
the above claim.\\

Proof of Claim: First we note that
\begin{equation*}
\frac{1}{\sqrt{\mu}}(E+v\times B )\cdot\nabla_v
(\sqrt{\mu}f)=(E+v\times B )\cdot\nabla_v
f-E\cdot\frac{v\sqrt{\mu}}{2}f.
\end{equation*}
We compute only one term and other terms can be treated
similarly.
\begin{align*}
(\partial_{\tau}^{\beta}\{E\cdot\nabla_v f\}
,\partial_{\tau}^{\beta}g)&=\int\int\partial_{\tau_1}E\cdot
(\partial_{\tau_2}^{\beta}\nabla_v f)
\partial_{\tau}^{\beta}g dvdx\;\;(\text{where }\tau_1+\tau_2=\tau)\\
&\leq(\int\int|\partial_{\tau_1}E\cdot
(\partial_{\tau_2}^{\beta}\nabla_v
f)|^2dvdx)^{\frac{1}{2}}||\partial_{\tau}^{\beta}g||\\
&\leq C
(\sum_{|\widetilde{\tau}|\leq\max\{s,3\}}||\partial_{\widetilde{\tau}}E||\cdot
||\partial_{\widetilde{\tau}}^{\beta}\nabla_v
f||)\cdot||\partial_{\tau}^{\beta}g||
\end{align*}
At the last step we have applied the Sobolev imbedding theorem.
Note that $\min\{|\tau_1|,|\tau_2|\}+2\leq\max\{|\tau|,3\}$.
Therefore, this finishes the first step.\\

\textbf{Step Two} : \textit{Construct the hydrodynamic field
$\mathbf{P_1}f_{r+1}$ and $\mathbf{P_2}g_{r+1}$ i.e.
$\rho_{r+1},\\u_{r+1},\theta_{r+1},\sigma_{r+1}$ and the
electromagnetic field $E_{r+1}$ and $B_{r+1}$.}

First of all, we recall that $\{I-P_{0} \mathbf{\}}u_{r+1}=\nabla
q_{r+1}$ for some scalar function $q_{r+1}$ and hence from
(\ref{hincon}) one obtains
\[
\Delta q_{r+1}=-\partial_{t}\rho_{r},
\]
where $\{I-P_{0}\}u_{r+1}$ has average zero. Hence the elliptic
$L^2$ estimates lead that $\sum_{|\tau|\leq
s}||\partial_{\tau}\{I-P_{0}\}u_{r+1}||$ is bounded by
$\sum_{|\tau|\leq s}||\partial_{\tau}\partial_{t}\rho_{r}||,$
which is bounded by the RHS of (\ref{i-p}) by the induction
hypothesis.

We determine $\rho_{r+1}+\theta_{r+1}$ from the Boussinesq
relation (\ref{hbous}):
\begin{align*}
\rho_{r+1}+\theta_{r+1}  &  =\Delta^{-1}\nabla\cdot
P_{0}\{-u_{1}\cdot
\nabla\{P_{0}u_{r}\}-u_{r}\cdot\nabla\{P_{0}u_{1}\}+R_{r}^{u}\}\\
&  +\langle\frac{|v|^{2}\sqrt{\mu}}{3},L^{-1}\{v\cdot\nabla_{x}f_{r-1}%
+\Gamma(f_{1},f_{r})\}\rangle-\frac{5}{2}\theta_{1}\theta_{r}-u_{r}\cdot
u_{1}.
\end{align*}
Notice that by the formula (\ref{ru}) and the induction
hypothesis, we easily conclude
 $||\Delta^{-1}\nabla\cdot
P_{0}\partial_{\tau}R_{r}^{u}|| $ is bounded by the RHS of
(\ref{i-p}). Hence $\sum_{|\tau|\leq s}||\partial_{\tau
}\{\rho_{r+1}+\theta_{r+1}\}||$ is bounded again by the RHS of
(\ref{i-p}).

Finally, to determine the remaining
$P_{0}u_{r+1},\sigma_{r+1},\theta_{r+1}$ and $E_{r+1}, B_{r+1}$
we solve the $(r+1)$-th order\ \textit{linear}
Vlasov-Navier-Stokes-Fourier system
(\ref{hnavier-stokes})-(\ref{compaEB}). By the standard energy
method ($\nabla\cdot P_{0}u_{r+1}=0$), there is the unique
solution
$[P_{0}u_{r+1},\\
\sigma_{r+1},\theta_{r+1},E_{r+1},B_{r+1}]$. We illustrate the
energy estimates in the subsequence. First, $B_{r+1}$ is
estimated and other terms will be done all together. From
(\ref{hmag}) and the Poincar$\acute{e}$ inequality, we get
\begin{equation*}
||\partial_{\tau}\nabla B_{r+1}||\leq
||\{\mathbf{I-P_2}\}\partial_{\tau}g_{r+1}||
+||\partial_{\tau}\partial_t E_{r}||,\;\; ||B_{r+1}||\leq C
||\nabla B_{r+1}||.
\end{equation*}
By (\ref{i-p}) and the induction hypothesis, $\sum_{|\tau|\leq s
}||\partial_{\tau}B_{r+1}||$ is bounded by the desired quantity.

Next we move onto other terms. Take $\partial_{\tau}$ derivatives
of (\ref{hnavier-stokes}) when $m=r+1$, multiply by
$\partial_{\tau}P_0u_{r+1}$ and integrate over $\mathbb{T}^3$ to
get
\begin{equation*}
\begin{split}
&\frac{1}{2}\frac{d}{dt}||\partial_{\tau}P_{0}u_{r+1}||^{2}+\eta
||\nabla\partial_{\tau}P_{0}u_{r+1}||^{2}\\
&=-(\partial_{\tau}\{u_1\cdot\nabla
P_0u_{r+1}+P_0u_{r+1}\cdot\nabla
u_1\},\partial_{\tau}P_0u_{r+1})\\
&\;\;\;+(\partial_{\tau}\{E_1\sigma_{r+1}+E_{r+1}\sigma_1\}
,\partial_{\tau}P_{0}u_{r+1})+(\partial_{\tau
}R_{r+1}^{u},\partial_{\tau}P_{0}u_{r+1}).
\end{split}
\end{equation*}
By the Cauchy-Schwartz inequality and Sobolev imbedding theorem,
first two terms in the RHS can be estimated as follows:
\[
\begin{split}
(&\partial_{\tau}\{u_1\cdot\nabla P_0u_{r+1}+P_0u_{r+1}\cdot\nabla
u_1\},\partial_{\tau}P_0u_{r+1})\\
&\leq \frac{\eta}{8}\sum_{|\tau|\leq s+1} ||\partial_{\tau}P_0
u_{r+1}||^2+C_{\eta}\sum_{|\tau|\leq s+2}
||\partial_{\tau}u_1||^4;\\
(&\partial_{\tau}\{E_1\sigma_{r+1}+E_{r+1}\sigma_1\},
\partial_{\tau}P_{0}u_{r+1})\\
&\leq \frac{\eta}{8}\sum_{|\tau|\leq s+1} ||\partial_{\tau}P_0
u_{r+1}||^2+\frac{\alpha}{4}\sum_{|\tau|\leq s}
||\partial_{\tau}\sigma_{r+1}||^2+\frac{\alpha}{8}\sum_{|\tau|\leq
s} ||\partial_{\tau}E_{r+1}||^2\\
&\;\;\;+C_{\eta,\alpha}\sum_{|\tau|\leq s+1}
||\partial_{\tau}\sigma_1||^4.
\end{split}
\]
In the same fashion, starting with (\ref{hvlasov}) and
(\ref{hfourier1}) respectively, one gets
\[
\begin{split}
&\frac{1}{2}\frac{d}{dt}||\partial_{\tau}\sigma_{r+1}||^{2}+\alpha
||\nabla\partial_{\tau}\sigma_{r+1}||^{2}
+\alpha||\partial_{\tau}\sigma_{r+1}||^2\\
&=-(\partial_{\tau}\{u_1\cdot\nabla
\sigma_{r+1}+P_0u_{r+1}\cdot\nabla\sigma_1\},
\partial_{\tau}\sigma_{r+1})+(\partial_{\tau
}R_{r+1}^{\sigma},\partial_{\tau}\sigma_{r+1}),\\
&\leq \frac{\alpha}{6}\sum_{|\tau|\leq s+1}
||\partial_{\tau}\sigma_{r+1}||^2+\frac{\eta}{8}\sum_{|\tau|\leq
s} ||\partial_{\tau}P_0 u_{r+1}||^2+C_{\alpha}\sum_{|\tau|\leq
s+2}||\partial_{\tau}\sigma_1||^4\\
&\;\;\;\;+(\partial_{\tau
}R_{r+1}^{\sigma},\partial_{\tau}\sigma_{r+1});\\
&\frac{1}{2}\frac{d}{dt}||\partial_{\tau}\theta_{r+1}||^{2}+\kappa
||\nabla\partial_{\tau}\theta_{r+1}||^{2}\\
&=-(\partial_{\tau}\{u_1\cdot\nabla
\theta_{r+1}+P_0u_{r+1}\cdot\nabla
\theta_1\},\partial_{\tau}\theta_{r+1})+(\partial_{\tau
}R_{r+1}^{\theta},\partial_{\tau}\theta_{r+1})\\
&\leq \frac{\kappa}{2}\sum_{|\tau|\leq s+1}
||\partial_{\tau}\theta_{r+1}||^2+\frac{\eta}{8}\sum_{|\tau|\leq
s} ||\partial_{\tau}P_0 u_{r+1}||^2+C_{\kappa}\sum_{|\tau|\leq
s+2}||\partial_{\tau}\theta_1||^4\\
&\;\;\;\;+(\partial_{\tau
}R_{r+1}^{\theta},\partial_{\tau}\theta_{r+1}).
\end{split}
\]
For the electric field, from (\ref{helec}), we obtain
\[
||\partial_{\tau}\nabla E_{r+1}||\leq||\partial_{\tau}\partial_t
B_{r}||+||\partial_{\tau}\sigma_{r+1}||.
\]
The conservation laws (\ref{conservation}) and (\ref{compaEB}) are
utilized to handle no derivative terms: by Poincar$\acute{e}$
inequality,
\[
\begin{split}
&||P_0u_{r+1}||\leq C ||\nabla P_0u_{r+1}||+C|\int P_0u_{r+1}
dx|\\
&\;\;\;\;\;\;\;\;\;\;\;\;\;\;\;\leq C ||\nabla P_0u_{r+1}||+
C\sum_{1\leq i\leq r}(||E_i||^2+||B_i||^2);\\
&||\sigma_{r+1}||\leq C ||\nabla \sigma_{r+1}||;\\
 &||\theta_{r+1}||\leq C ||\nabla \theta_{r+1}||+
C\sum_{1\leq i\leq r}(||E_i||^2+||B_i||^2);\\
&||E_{r+1}||\leq C||\nabla E_{r+1}||+C|\int E_{r+1}dx|\\
&\;\;\;\;\;\;\;\;\;\;\;\;\leq C||\nabla E_{r+1}||+C|e^{-\alpha
t}\int E_{r+1}^0dx+ \int_0^t e^{-\alpha (t-s)}\ell_r ds|.
\end{split}
\]
Based on the above estimates, following the argument in
\cite{Guo2} (p.30-31), by the induction hypothesis and the
Gronwall lemma, one can verify that (\ref{rdecay}) holds for
$m=r+1$ and thus it completes the proof of the theorem.
\end{proof}\

\section{Uniform Spatial Energy Estimate}\

In this section, we consider the following model problem:%
\begin{equation}
\begin{split}
\partial_{t}f^{\varepsilon}+\frac{v}{\varepsilon}\cdot\nabla_{x}%
f^{\varepsilon}+\frac{1}{\varepsilon^{2}}Lf^{\varepsilon}
=h_1^{\varepsilon},\\
\label{model}%
\partial_{t}g^{\varepsilon}+\frac{v}{\varepsilon}\cdot(\nabla_{x}%
g^{\varepsilon}-\sqrt{\mu}E^{\varepsilon})+\frac{1}{\varepsilon^{2}}
\mathcal{L}g^{\varepsilon}=h_2^{\varepsilon},
\end{split}
\end{equation}
coupled with the Maxwell equations
\begin{equation}
\begin{split}
\varepsilon\partial_t E^{\varepsilon}-\nabla\times B^{\varepsilon}
=-\int_{\mathbb{R}^3} g^{\varepsilon}
v\sqrt{\mu}dv+j_1^{\varepsilon}
,\;\;\nabla\cdot B^{\varepsilon}=0,\\
\varepsilon\partial_t B^{\varepsilon}+\nabla\times
E^{\varepsilon}=j_2^{\varepsilon},\;\;\nabla\cdot
E^{\varepsilon}=\int_{\mathbb{R}^3}
g^{\varepsilon}\sqrt{\mu}dv,\label{model2}
\end{split}
\end{equation}
where $h_1^{\varepsilon}$, $h_2^{\varepsilon}$,
$j_1^{\varepsilon}$ and $j_2^{\varepsilon}$ will be given. We
shall establish a uniform space-time energy estimate for
$f^{\varepsilon},g^{\varepsilon},E^{\varepsilon}\text{ and }
B^{\varepsilon}$.

We use a different representation for the hydrodynamic field
(fluid) parts $\mathbf{P_1}f^{\varepsilon}$,
$\mathbf{P_2}g^{\varepsilon}$ in a different way as:
\[
[{\mathbf{P_1}f^{\varepsilon}},\;{\mathbf{P_2}g^{\varepsilon}}]
=[{\{a^{\varepsilon}(t,x)+b^{\varepsilon}(t,x)\cdot
v+c^{\varepsilon}(t,x)|v|^{2}\}\sqrt{\mu}},\;{d^{\varepsilon}(t,x)
\sqrt{\mu}}].
\]
Our goal is to estimate
$a^{\varepsilon}(t,x),b^{\varepsilon}(t,x),c^{\varepsilon}(t,x),
d^{\varepsilon}(t,x)$ and
$E^{\varepsilon}(t,x),B^{\varepsilon}(t,x)$ in terms of
$\mathbf{(I-P_1)}f^{\varepsilon}$ and
$\mathbf{(I-P_2)}g^{\varepsilon}$.

We remark that our argument in this section is valid for all
$\varepsilon\leq 1$. In particular, the case of $\varepsilon=1$
yields another proof of the existence of the classical solution
for the Vlasov-Maxwell-Boltzmann system without using
the time derivative. \\

\begin{lem}
\label{abc}Assume
$f^{\varepsilon},g^{\varepsilon},E^{\varepsilon},B^{\varepsilon}$
are solutions to the system (\ref{model}), (\ref{model2}) such
that for all $t\geq 0$,
\begin{equation}
\begin{split}
&\int_{\mathbb{T}^3}
a^{\varepsilon}(t,x)dx=O(\varepsilon(||[E^{\varepsilon},
B^{\varepsilon}]||^2(t)))+O(\varepsilon \mathcal{A}),\\  \label{mme}%
&\int_{\mathbb{T}^3}
b^{\varepsilon}(t,x)dx=O(\varepsilon(||[E^{\varepsilon},
B^{\varepsilon}]||^2(t)))+O(\varepsilon \mathcal{A}),\\
&\int_{\mathbb{T}^3} c^{\varepsilon}(t,x)dx=
O(\varepsilon(||[E^{\varepsilon},
B^{\varepsilon}]||^2(t)))+O(\varepsilon \mathcal{A}),\\
&\int_{\mathbb{T}^3} d^{\varepsilon}(t,x)dx=0 \;\text{ and }
\int_{\mathbb{T}^3} B^{\varepsilon}(t,x)dx=0,\\
\end{split}
\end{equation}
where we have used the standard big $O$ notation in $\varepsilon$
and $\mathcal{A}(t)$ is a function of $t$. And suppose that
$||E^{\varepsilon}||^2+||B^{\varepsilon}||^2\leq
M$ for some constant $M$. Then there exists a constant $C_{1}>0$ such that%
\begin{equation}
\begin{split}
\sum_{|\gamma|\leq
N}&||[\partial_{\gamma}\mathbf{P_1}f^{\varepsilon},\partial_{\gamma}
{\mathbf{P_2}g^{\varepsilon}}]
||^2\leq\varepsilon\frac{dG(t)}{dt}+\frac{C_{1}}{\varepsilon^{2}}
\sum_{|\gamma|\leq
N}||[\partial_{\gamma}\mathbf{\{I-P_1\}}f^{\varepsilon},
\partial_{\gamma}
{\mathbf{\{I-P_2\}}g^{\varepsilon}}]||_{\nu}%
^{2}\\&+C_{1}\varepsilon^{2}\sum_{|\gamma|\leq
N-1}||\partial_{\gamma} h_{||}^{\varepsilon}||^{2}
+\varepsilon^2\mathcal{A}^2+C_1\varepsilon^2
(||j_1^{\varepsilon}||^2+||\nabla j_1^{\varepsilon}||^2+
||j_2^{\varepsilon}||^2+||\nabla j_2^{\varepsilon}||^2) \label{diff}%
\end{split}
\end{equation}
where $\varepsilon\leq 1$  and $G(t)$ is a function of $t$
satisfying
\begin{equation}
|G(t)|\leq \sum_{|\gamma|\leq
N}||[\partial_{\gamma}f^{\varepsilon},\partial_{\gamma}
{g^{\varepsilon}}]||^2(t)+ \sum_{|\gamma|\leq
N}||[\partial_{\gamma}{E^{\varepsilon}},\partial_{\gamma}
{B^{\varepsilon}}]||^2(t), \label{gbound}
\end{equation}
and $\partial_{\gamma}h_{||}^{\varepsilon}$ is the $L_{v}^{2}$
projection of
$[\partial_{\gamma}{h_1^{\varepsilon}},\partial_{\gamma}
{h_2^{\varepsilon}}](t,x,v)$ on the subspace generated by
\[
[\sqrt{\mu},v_{i}\sqrt{\mu},v_{i}v_{j}\sqrt{\mu},v_{i}|v|^{2}\sqrt{\mu
}].
\]
\end{lem}\

\begin{proof}
As well illustrated in \cite{Guo2}, there are two fundamental
ingredients in the proof. First of all, we use the \textit{Local
Conservation Laws:} By multiplying
$\sqrt{\mu},v\sqrt{\mu},|v|^{2}\sqrt{\mu}$ with the first
equation of (\ref{model}), $\sqrt{\mu}$ with the second equation
and then integrating in $v\in\mathbb{R}^{3}$, we obtain the
following:
\begin{equation}
\begin{split}
&  a_{t}^{\varepsilon}=\frac{1}{2\varepsilon}\langle v\cdot\nabla
_{x}\mathbf{(I-P_1)}f^{\varepsilon},|v|^{2}\sqrt{\mu}\rangle+\langle
h_1^{\varepsilon},\{\frac{5}{2}-\frac{|v|^{2}}{2}\}\sqrt{\mu}\rangle,\\
&  c_{t}^{\varepsilon}+\frac{1}{3\varepsilon}\nabla_{x}\cdot
b^{\varepsilon }=-\frac{1}{6\varepsilon}\langle
v\cdot\nabla_{x}\mathbf{(I-P_1)}f^{\varepsilon
},|v|^{2}\sqrt{\mu}\rangle+\langle h_1^{\varepsilon},\{\frac{|v|^{2}}{6}%
-\frac{1}{2}\}\sqrt{\mu}\rangle,\label{abct}\\
&
b_{t}^{\varepsilon}+\frac{1}{\varepsilon}\{\nabla_{x}a^{\varepsilon
}+5\nabla_{x}c^{\varepsilon}\}=-\frac{1}{\varepsilon}\langle
v\cdot\nabla
_{x}\mathbf{(I-P_1)}f^{\varepsilon},v\sqrt{\mu}\rangle+\langle
h_1^{\varepsilon },v\sqrt{\mu}\rangle,\\
& d_{t}^{\varepsilon}=-\frac{1}{\varepsilon}\langle v\cdot\nabla
_{x}\mathbf{(I-P_2)}g^{\varepsilon},\sqrt{\mu}\rangle+\langle
h_2^{\varepsilon},\sqrt{\mu} \rangle.\\
\end{split}
\end{equation}

The second ingredient is the study of the \textit{Macroscopic
Equations}: notice that by plugging $f^{\varepsilon}\equiv\mathbf{P_1}%
f^{\varepsilon}+\mathbf{(I-P_1)}f^{\varepsilon}$ and $g^{\varepsilon}\equiv\mathbf{P_2}%
g^{\varepsilon}+\mathbf{(I-P_2)}g^{\varepsilon}$ into
(\ref{model}),
\begin{align*}
&  \varepsilon\{a_{t}^{\varepsilon}+b_{t}^{\varepsilon}\cdot v+c_{t}%
^{\varepsilon}|v|^{2}\}\sqrt{\mu}+v\cdot\{\nabla_{x}a^{\varepsilon}+\nabla
_{x}b^{\varepsilon}\cdot v+\nabla_{x}c^{\varepsilon}|v|^{2}\}\sqrt{\mu}\\
&
=-\{\varepsilon\partial_{t}+v\cdot\nabla_{x}\}\mathbf{(I-P_1)}f^{\varepsilon
}-\frac{1}{\varepsilon}L\mathbf{(I-P_1)}f^{\varepsilon}+\varepsilon
h_1^{\varepsilon},\\
& \varepsilon\;
d_{t}^{\varepsilon}\sqrt{\mu}+v\cdot\{\nabla_{x}d^{\varepsilon}
-E^{\varepsilon}\}\sqrt{\mu}\\
&
=-\{\varepsilon\partial_{t}+v\cdot\nabla_{x}\}\mathbf{(I-P_2)}g^{\varepsilon
}-\frac{1}{\varepsilon}\mathcal{L}\mathbf{(I-P_2)}g^{\varepsilon}+\varepsilon
h_2^{\varepsilon}.
\end{align*}
Fixing $t$ and $x,$ and comparing the coefficients on both sides
in front of the
$[\sqrt{\mu},v\sqrt{\mu},\\v_{i}v_{j}\sqrt{\mu},v_{i}|v|^{2}\sqrt{\mu}],$
we obtain the \textit{macroscopic} equations as%
\begin{align}
\nabla_{x}c^{\varepsilon}  &  =l_{c}^{\varepsilon}+\varepsilon h_{c}%
^{\varepsilon},\label{c}\\
\varepsilon\partial_{t}c^{\varepsilon}+\partial_{i}b_{i}^{\varepsilon}
&
=l_{i}^{\varepsilon}+\varepsilon h_{i}^{\varepsilon},\label{bi}\\
\partial_{i}b_{j}^{\varepsilon}+\partial_{j}b_{i}^{\varepsilon}  &
=l_{ij}^{\varepsilon}+\varepsilon h_{ij}^{\varepsilon}\;\;\text{for }
i\neq j,\label{bij}\\
\varepsilon\partial_{t}b_{i}^{\varepsilon}+\partial_{i}a^{\varepsilon}
&
=l_{bi}^{\varepsilon}+\varepsilon h_{bi}^{\varepsilon},\label{ai}\\
\varepsilon\partial_{t}a^{\varepsilon}  &
=l_{a}^{\varepsilon}+\varepsilon h_{a}^{\varepsilon},
\label{adot}\\
 \varepsilon\partial_{t}d^{\varepsilon}  &
=l_{d}^{\varepsilon}+\varepsilon h_{d}^{\varepsilon},\label{d}\\
 \nabla_x d^{\varepsilon}-E^{\varepsilon}  &
=l_{e}^{\varepsilon}+\varepsilon h_{e}^{\varepsilon}.\label{e}
\end{align}
Here by elementary linear algebra, the linear parts $l_{c}^{\varepsilon}%
,l_{i}^{\varepsilon},l_{ij}^{\varepsilon},l_{bi}^{\varepsilon},
l_{a}^{\varepsilon},l_{d}^{\varepsilon}$ and $l_{e}^{\varepsilon}$
are all of the form
\begin{equation}
\begin{split}
\text{either }\;&\langle-\{\varepsilon\partial_{t}+v\cdot\nabla_{x}
\}\mathbf{(I-P_1)}%
f^{\varepsilon},\zeta\rangle-\frac{1}{\varepsilon}\langle L\mathbf{(I-P_1)}%
f^{\varepsilon},\zeta\rangle\\ \label{projection}%
\text{ or}\;\;\;\;&\langle-\{\varepsilon\partial_{t}+v\cdot\nabla_{x}\}
\mathbf{(I-P_2)}%
g^{\varepsilon},\zeta\rangle-\frac{1}{\varepsilon}\langle \mathcal{L}\mathbf{(I-P_2)}%
g^{\varepsilon},\zeta\rangle
\end{split}
\end{equation}
where $\zeta$ is a (different) linear combination of the basis
$[\sqrt{\mu
},v\sqrt{\mu},v_{i}v_{j}\sqrt{\mu},v_{i}|v|^{2}\sqrt{\mu}]$
accordingly$,$
while $h_{c}^{\varepsilon},h_{i}^{\varepsilon},h_{ij}^{\varepsilon}%
,h_{bi}^{\varepsilon},h_{a}^{\varepsilon},h_{d}^{\varepsilon}$
and $h_{e}^{\varepsilon}$ are defined as $\langle
h^{\varepsilon},\zeta\rangle$ with same choices of $\zeta.$ Notice that%
\[
||\partial_{\gamma}h_{c}^{\varepsilon}||+||\partial_{\gamma}h_{i}%
^{\varepsilon}||+||\partial_{\gamma}h_{ij}^{\varepsilon}||+||\partial_{\gamma
}h_{bi}^{\varepsilon}||+||\partial_{\gamma}h_{a}^{\varepsilon}||
+||\partial_{\gamma}h_{d}^{\varepsilon}||+||\partial_{\gamma}h_{e}^{\varepsilon}||\leq
C||\partial_{\gamma}h_{||}^{\varepsilon}||.
\]

The macroscopic equations (\ref{c})-(\ref{adot}) have the same
structure as the pure Boltzmann case. Following the proof of Lemma
6.1 in \cite{Guo2}, we can deduce the estimates on
$\nabla_{x}\partial_{\gamma_1}a^{\varepsilon},\;
\nabla_{x}\partial_{\gamma_1}b^{\varepsilon}$ and
$\nabla_{x}\partial_{\gamma_1}c^{\varepsilon}$: for
$|\gamma_1|\leq N-1$,
\begin{align*}
\frac{1}{2}||\nabla\partial_{\gamma_1}b^{\varepsilon}||^{2} &
\leq-\frac
{d}{dt}\int_{\mathbb{T}^{3}}\langle\varepsilon\mathbf{(I-P_1)}\partial_{\gamma_1
}f^{\varepsilon},\zeta_{ij}\rangle\cdot\partial_{j}\partial_{\gamma_1
}b^{\varepsilon}dx\\
&  +\frac{C}{\varepsilon^{2}}\{||\nabla_{x}\partial_{\gamma_1}\mathbf{(I-P_1)}%
f^{\varepsilon}||_{\nu}^{2}+||\partial_{\gamma_1}\mathbf{(I-P_1)}f^{\varepsilon
}||_{\nu}^{2}\}\\
&  +C\varepsilon^{2}\{||\nabla_{x}\partial_{\gamma_1}a^{\varepsilon}%
||^{2}+||\nabla_{x}\partial_{\gamma_1}c^{\varepsilon}||^{2}+||\partial_{\gamma_1
}h_{||}^{\varepsilon}||^{2}\},\\
\frac{1}{2}||\nabla\partial_{\gamma_1}c^{\varepsilon}||^{2}  &  \leq-\frac{d}{dt}%
\int_{\mathbb{T}^{3}}\langle\varepsilon\mathbf{(I-P_1)}\partial_{\gamma_1
}f^{\varepsilon},\zeta_{c}\rangle\cdot\nabla_{x}\partial_{\gamma_1
}c^{\varepsilon}dx\\
&  +\frac{C}{\varepsilon^{2}}\{||\nabla_{x}\partial_{\gamma_1}\mathbf{(I-P_1)}%
f^{\varepsilon}||_{\nu}^{2}+||\partial_{\gamma_1}\mathbf{(I-P_1)}f^{\varepsilon
}||_{\nu}^{2}\}\\
&  +C\varepsilon^{2}\{||\nabla_{x}\partial_{\gamma_1}a^{\varepsilon}%
||^{2}+||\nabla_{x}\partial_{\gamma_1}c||^{2}+||\partial_{\gamma_1}h_{||}^{\varepsilon
}||^{2}\},\\
\frac{1}{2}||\nabla\partial_{\gamma_1}a^{\varepsilon}||^{2}  &
\leq
-\frac{d}{dt}\{\int_{\mathbb{T}^{3}}\langle\varepsilon\mathbf{(I-P_1)}%
\partial_{\gamma_1}f^{\varepsilon},\zeta\rangle\cdot\nabla_{x}\partial_{\gamma_1
}a^{\varepsilon}dx+\int_{\mathbb{T}^{3}}\varepsilon\partial_{\gamma_1
}b^{\varepsilon}\cdot\nabla_{x}\partial_{\gamma_1}a^{\varepsilon}dx\}\\
&  +\frac{C}{\varepsilon^{2}}\{||\nabla_{x}\partial_{\gamma_1}\mathbf{(I-P_1)}%
f^{\varepsilon}||_{\nu}^{2}+||\partial_{\gamma_1}\mathbf{(I-P_1)}f^{\varepsilon
}||_{\nu}^{2}\}\\
&  +C\varepsilon^{2}\{||\partial_{\gamma_1}h_{||}^{\varepsilon}||^{2}%
+||\nabla\partial_{\gamma_1}b^{\varepsilon}||^{2}\}.
\end{align*}
We shall, however, estimate
$\nabla_{x}\partial_{\gamma_1}d^{\varepsilon}$ in the same spirit.
\begin{align*}
 \Delta\partial_{\gamma_1}d_{i}^{\varepsilon}&=\sum_{j}\partial_{jj}%
\partial_{\gamma_1}d_{i}^{\varepsilon}\\
&  =\sum_{j}\partial_{j}[\partial_{\gamma_1}E^{\varepsilon}_j
+\partial_{\gamma_1}l_{e}^{\varepsilon}+\varepsilon\partial_{\gamma_1}
h_e^{\varepsilon}]\;\;\text{by (\ref{e})}\\
&=\partial_{\gamma_1}d^{\varepsilon}
+\sum_{j}[\partial_{j}\partial_{\gamma_1}l_{e}^{\varepsilon}+
\varepsilon\partial_{j}\partial_{\gamma_1} h_e^{\varepsilon}]\;\;
(\text{since } \nabla\cdot E^{\varepsilon}= d^{\varepsilon})
\end{align*}
Multiply the above by $\partial_{\gamma_1}d^{\varepsilon}$ and
integrate to get
\begin{equation}
\int_{\mathbb{T}^{3}}|\nabla_x\partial_{\gamma_1}d^{\varepsilon}|^2
+|\partial_{\gamma_1}d^{\varepsilon}|^2dx=\sum_{j}\int_{\mathbb{T}^{3}}
(\partial_{\gamma_1}l_{e}^{\varepsilon}+
\varepsilon\partial_{\gamma_1}
h_e^{\varepsilon})\cdot\partial_{j}\partial_{\gamma_1}d^{\varepsilon}dx.
\label{ddd}
\end{equation}
Recall that
\begin{align*}
\partial_{\gamma_1}l_e^{\varepsilon}=&-\langle\varepsilon\partial_t
\mathbf{(I-P_2)}\partial_{\gamma_1}g^{\varepsilon},v\sqrt{\mu}\rangle\\
&-\langle v\cdot\nabla_x
\mathbf{(I-P_2)}\partial_{\gamma_1}g^{\varepsilon},v\sqrt{\mu}\rangle
-\frac{1}{\varepsilon}\langle \mathcal{L}
\mathbf{(I-P_2)}\partial_{\gamma_1}g^{\varepsilon},v\sqrt{\mu}\rangle.
\end{align*}
Last two terms are of the desired form. As for the first one, we
integrate it by parts in the $t$ variable:
\begin{align*}
\int&_{\mathbb{T}^{3}}\langle\varepsilon\partial_t
\mathbf{(I-P_2)}\partial_{\gamma_1}g^{\varepsilon},v\sqrt{\mu}\rangle
\cdot\partial_{j}\partial_{\gamma_1}d^{\varepsilon}dx\\
=&\frac{d}{dt}\int_{\mathbb{T}^{3}}\langle\varepsilon
\mathbf{(I-P_2)}\partial_{\gamma_1}g^{\varepsilon},v\sqrt{\mu}\rangle
\cdot\partial_{j}\partial_{\gamma_1}d^{\varepsilon}dx\\
&+ \int_{\mathbb{T}^{3}}\langle
\mathbf{(I-P_2)}\partial_{j}\partial_{\gamma_1}g^{\varepsilon},v\sqrt{\mu}\rangle
\cdot\varepsilon\partial_t\partial_{\gamma_1}d^{\varepsilon}dx.
\end{align*}
Replace $\varepsilon\partial_t\partial_{\gamma_1}d^{\varepsilon}$
with
$\partial_{\gamma_1}l_{d}^{\varepsilon}+\varepsilon\partial_{\gamma_1}
h_{d}^{\varepsilon}$ in the latter integral by the macroscopic
equation (\ref{d}). Thus (\ref{ddd}) leads to the following:
\begin{align}
\frac{1}{2}(||\nabla\partial_{\gamma_1}d^{\varepsilon}||^{2}+
||\partial_{\gamma_1}d^{\varepsilon}||^{2}) & \leq-\frac
{d}{dt}\int_{\mathbb{T}^{3}}\langle\varepsilon\mathbf{(I-P_2)}
\partial_{\gamma_1}g^{\varepsilon},v\sqrt{\mu}\rangle
\cdot\partial_{j}\partial_{\gamma_1}d^{\varepsilon}dx\label{dest1}\nonumber\\
&  +\frac{C}{\varepsilon^{2}}\{||\nabla_{x}\partial_{\gamma_1}\mathbf{(I-P_2)}%
g^{\varepsilon}||_{\nu}^{2}+||\partial_{\gamma_1}\mathbf{(I-P_2)}
g^{\varepsilon}||_{\nu}^{2}\}\\
&+C\varepsilon^{2}||\partial_{\gamma_1}h_{||}^{\varepsilon}||^{2}.
\nonumber
\end{align}

The field estimate, the new ingredient of the proof, comes into
play for the need of the estimates on no derivative terms
$a^{\varepsilon},$ $b^{\varepsilon}$ and $c^{\varepsilon}$ to
complete the proof of Lemma: from (\ref{mme}) we obtain by the
Poincar$\acute{e}$ inequality,
\begin{equation}
\begin{split}
||a^{\varepsilon}||^2+||b^{\varepsilon}||^2+||c^{\varepsilon}||^2
\leq ||\nabla a^{\varepsilon}||^2+||&\nabla
b^{\varepsilon}||^2+||\nabla
c^{\varepsilon}||^2\\&+\varepsilon^2(||E^{\varepsilon}||^2+
||B^{\varepsilon}||^2)^2+\varepsilon^2\mathcal{A}^2.\label{abceb}
\end{split}
\end{equation}
Also it will play an important role not only to close the energy
estimates in later sections but also to achieve decay rates
(\ref{rdecay}) and (\ref{decayrate}). The intriguing computation
is carried out in two steps: firstly, the electric field will be
estimated via the macroscopic equation (\ref{e}) and then the
Maxwell equations (\ref{model2}) will give rise to the estimate
on the magnetic field.

We start with (\ref{e}): $E^{\varepsilon}=\nabla_x d^{\varepsilon}
-l_{e}^{\varepsilon}-\varepsilon h_{e}^{\varepsilon}$.
\begin{equation}
\begin{split}
||\partial_{\gamma_1}E^{\varepsilon}||^2&=\int_{\mathbb{T}^3}(\nabla_x
\partial_{\gamma_1}d^{\varepsilon} -\partial_{\gamma_1}l_{e}^{\varepsilon}
-\varepsilon\partial_{\gamma_1}
h_{e}^{\varepsilon})\cdot \partial_{\gamma_1}E^{\varepsilon} dx\\
&\leq \int_{\mathbb{T}^3}(\nabla_x
\partial_{\gamma_1}d^{\varepsilon}-\varepsilon
\partial_{\gamma_1}h_{e}^{\varepsilon})\cdot \partial_{\gamma_1}E^{\varepsilon}
dx-\int_{\mathbb{T}^3}\partial_{\gamma_1}
l_{e}^{\varepsilon}\cdot \partial_{\gamma_1}E^{\varepsilon}dx\\
&\leq-||\partial_{\gamma_1}
d^{\varepsilon}||^2+\frac{||\partial_{\gamma_1}E^{\varepsilon}||^2}{4}+
\varepsilon^2||\partial_{\gamma_1}h_{e}^{\varepsilon}||^2-
\int_{\mathbb{T}^3} \partial_{\gamma_1}l_{e}^{\varepsilon}\cdot
\partial_{\gamma_1}E^{\varepsilon}dx\label{ee}
\end{split}
\end{equation}
The last term needs a special care.
\begin{equation*}
\begin{split}
-\int_{\mathbb{T}^3} \partial_{\gamma_1}l_{e}^{\varepsilon}\cdot
\partial_{\gamma_1}E^{\varepsilon}dx=&
\int_{\mathbb{T}^3}\langle\varepsilon\partial_t\partial_{\gamma_1}\mathbf{(I-P_2)}
g^{\varepsilon},v\sqrt{\mu}\rangle\cdot \partial_{\gamma_1}E^{\varepsilon}dx\\
+ \int_{\mathbb{T}^3}&\langle
v\cdot\nabla_x\partial_{\gamma_1}\mathbf{(I-P_2)}
g^{\varepsilon}+\frac{1}{\varepsilon}\mathcal{L}\mathbf{(I-P_2)}
\partial_{\gamma_1}g^{\varepsilon}
,v\sqrt{\mu}\rangle\cdot \partial_{\gamma_1}E^{\varepsilon}dx\\
\equiv&\;(I)\;+\;(II)
\end{split}
\end{equation*}
It is easy to see that the second term $(II)$ is bounded by
$$\frac{||\partial_{\gamma_1}E^{\varepsilon}||^2}{4}+\frac{C}{\varepsilon^2}
\{||\nabla_x\mathbf{(I-P_2)}\partial_{\gamma_1}
g^{\varepsilon}||_{\nu}^2 +||\mathbf{(I-P_2)}
\partial_{\gamma_1}g^{\varepsilon}||_{\nu}^2\}.$$
As for the first term $(I)$, we first integrate it by parts in
$t$:
\[
\frac{d}{dt}\int_{\mathbb{T}^3}\langle\varepsilon\mathbf{(I-P_2)}
\partial_{\gamma_1}g^{\varepsilon},v\sqrt{\mu}\rangle\cdot
\partial_{\gamma_1} E^{\varepsilon}dx-
\int_{\mathbb{T}^3}\langle\mathbf{(I-P_2)}
\partial_{\gamma_1}g^{\varepsilon},v\sqrt{\mu}\rangle\cdot\varepsilon
\partial_t\partial_{\gamma_1}
E^{\varepsilon}dx.\\
\]
By using the Maxwell equation (\ref{model2}) to eliminate
$\varepsilon\partial_t\partial_{\gamma_1} E^{\varepsilon}$ in the
second term, we get the following:
\begin{align*}
&- \int_{\mathbb{T}^3}\langle\mathbf{(I-P_2)}
\partial_{\gamma_1}g^{\varepsilon},v\sqrt{\mu}\rangle\cdot\varepsilon
\partial_t\partial_{\gamma_1}
E^{\varepsilon}dx\\
&=\int_{\mathbb{T}^3}|\langle\mathbf{(I-P_2)}
\partial_{\gamma_1}g^{\varepsilon},v\sqrt{\mu}\rangle|^2dx
+ \int_{\mathbb{T}^3}\langle\mathbf{(I-P_2)}
\partial_{\gamma_1}g^{\varepsilon},v\sqrt{\mu}\rangle\cdot
\partial_{\gamma_1}j_1^{\varepsilon}\;dx\\
&\;\;\;-\int_{\mathbb{T}^3}\langle\mathbf{(I-P_2)}
\partial_{\gamma_1}g^{\varepsilon},v\sqrt{\mu}\rangle\cdot\nabla\times
\partial_{\gamma_1}B^{\varepsilon}dx\\
&\leq\frac{3}{2}||\mathbf{(I-P_2)}
\partial_{\gamma_1}g^{\varepsilon}||_{\nu}^2+\frac{1}{2}
||\partial_{\gamma_1}j_1^{\varepsilon}||^2+C_{\xi}
||\nabla_x\mathbf{(I-P_2)}
\partial_{\gamma_1}g^{\varepsilon}||_{\nu}^2+
\xi||\partial_{\gamma_1}B^{\varepsilon}||^2,
\end{align*}
where $\xi$ is a fixed small number and we have integrated by
parts in $x$ to get the last inequality. Hence (\ref{ee}) leads to
\begin{equation}
\begin{split}
\frac{1}{2}||\partial_{\gamma_1}E^{\varepsilon}||^2+
||\partial_{\gamma_1}d^{\varepsilon}||^2\leq&\frac{d}{dt}\int_{\mathbb{T}^3}\langle\varepsilon\mathbf{(I-P_2)}
\partial_{\gamma_1}g^{\varepsilon},v\sqrt{\mu}\rangle\cdot
\partial_{\gamma_1}E^{\varepsilon}dx\label{e1}\\
+&C_{\xi} \{||\nabla_x\mathbf{(I-P_2)}
\partial_{\gamma_1}g^{\varepsilon}||_{\nu}^2+||\mathbf{(I-P_2)}
\partial_{\gamma_1}g^{\varepsilon}||_{\nu}^2\}\\+&
\varepsilon^2||\partial_{\gamma_1}h_{e}^{\varepsilon}||^2+\frac{1}{2}
||\partial_{\gamma_1}j_1^{\varepsilon}||^2+
\xi||\partial_{\gamma_1}B^{\varepsilon}||^2.
\end{split}
\end{equation}
As for the estimate on $\partial_{\gamma_1}B^{\varepsilon}$,
recall that $\nabla\times B^{\varepsilon}=\varepsilon\partial_t
E^{\varepsilon}+\langle\mathbf{(I-P_2)}g^{\varepsilon},
v\sqrt{\mu}\rangle-j_1^{\varepsilon}.$ Let $|\gamma_2|\leq N-2$.
\begin{equation*}
\begin{split}
||\nabla\times\partial_{\gamma_2}B^{\varepsilon}||^2=&\int_{\mathbb{T}^3}
(\varepsilon\partial_t\partial_{\gamma_2}
E^{\varepsilon}+\langle\mathbf{(I-P_2)}\partial_{\gamma_2}g^{\varepsilon},
v\sqrt{\mu}\rangle-\partial_{\gamma_2}j_1^{\varepsilon})\cdot\nabla\times\partial_{\gamma_2}B^{\varepsilon}dx\\
=&\frac{d}{dt}\int_{\mathbb{T}^3}\varepsilon\partial_{\gamma_2}
E^{\varepsilon}\cdot\nabla\times\partial_{\gamma_2}B^{\varepsilon}dx-
\int_{\mathbb{T}^3}\partial_{\gamma_2}
E^{\varepsilon}\cdot\nabla\times\varepsilon
\partial_t\partial_{\gamma_2}B^{\varepsilon}dx\\
&+\int_{\mathbb{T}^3}(\langle\mathbf{(I-P_2)}\partial_{\gamma_2}g^{\varepsilon},
v\sqrt{\mu}\rangle-\partial_{\gamma_2}j_1^{\varepsilon})\cdot\nabla\times\partial_{\gamma_2}B^{\varepsilon}dx
\end{split}
\end{equation*}
By (\ref{model2}), the second term containing $\varepsilon
\partial_t\partial_{\gamma_2}B^{\varepsilon}$ in the above can
be majorized as following:
\begin{align*}
& \int_{\mathbb{T}^3}\partial_{\gamma_2}
E^{\varepsilon}\cdot(\nabla\times\nabla\times\partial_{\gamma_2}
E^{\varepsilon}-\nabla\times\partial_{\gamma_2}j_2^{\varepsilon})dx\\
=&\int_{\mathbb{T}^3}\partial_{\gamma_2}
E^{\varepsilon}\cdot(\nabla(\nabla\cdot\partial_{\gamma_2}
E^{\varepsilon})-\Delta\partial_{\gamma_2}
E^{\varepsilon}-\nabla\times\partial_{\gamma_2}j_2^{\varepsilon})dx\\
=&\int_{\mathbb{T}^3}\partial_{\gamma_2}
E^{\varepsilon}\cdot\nabla\partial_{\gamma_2}d^{\varepsilon} dx+
\int_{\mathbb{T}^3}|\nabla\partial_{\gamma_2}E^{\varepsilon}|^2
dx-\int_{\mathbb{T}^3}\partial_{\gamma_2}
E^{\varepsilon}\cdot\nabla\times\partial_{\gamma_2}j_2^{\varepsilon}dx\\
\leq&\int_{\mathbb{T}^3}|\partial_{\gamma_2}E^{\varepsilon}|^2 +
|\nabla\partial_{\gamma_2}E^{\varepsilon}|^2
dx+\frac{1}{2}\int_{\mathbb{T}^3}|\nabla\partial_{\gamma_2}d^{\varepsilon}|^2
+|\nabla\times\partial_{\gamma_2} j_2^{\varepsilon}|^2 dx
\end{align*}
On the other hand, the last term is bounded by
\begin{align*}
C_{\xi}(||\mathbf{(I-P_2)}\partial_{\gamma_2}
g^{\varepsilon}||_{\nu}^2+||\partial_{\gamma_2}j_1^{\varepsilon}||^2)
+\xi ||\nabla\times\partial_{\gamma_2}B^{\varepsilon}||^2
\end{align*}
for any small number $\xi$. After absorbing $\xi
||\nabla\times\partial_{\gamma_2}B^{\varepsilon}||^2$ into the
LHS, we obtain for some $0<C_2=1-\xi<1$,
\begin{align}
C_2||\nabla\times\partial_{\gamma_2}B^{\varepsilon}||^2\leq&
\frac{d}{dt}\int_{\mathbb{T}^3}\varepsilon\partial_{\gamma_2}
E^{\varepsilon}\cdot\nabla\times\partial_{\gamma_2}B^{\varepsilon}dx+\int_{\mathbb{T}^3}|\partial_{\gamma_2}E^{\varepsilon}|^2
+ |\nabla\partial_{\gamma_2}E^{\varepsilon}|^2 dx
\nonumber\\
+\frac{1}{2}\int_{\mathbb{T}^3}|&\nabla\partial_{\gamma_2}d^{\varepsilon}|^2
+|\nabla\times\partial_{\gamma_2} j_2^{\varepsilon}|^2
dx+C_{\xi}(||\mathbf{(I-P_2)}\partial_{\gamma_2}
g^{\varepsilon}||_{\nu}^2+||\partial_{\gamma_2}j_1^{\varepsilon}||^2).\label{b1}
\end{align}
Letting $|\gamma_2|=|\gamma_1|-1\geq 0$, combine (\ref{e1}) with
(\ref{b1}) to get for some $0<C_3<1$,
\begin{align}
C_3||\partial_{\gamma_1}E^{\varepsilon}||^2
\leq&\frac{d}{dt}\int_{\mathbb{T}^3}\langle\varepsilon\mathbf{(I-P_2)}
\partial_{\gamma_1}g^{\varepsilon},v\sqrt{\mu}\rangle\cdot
\partial_{\gamma_1}E^{\varepsilon}+\xi C_2\varepsilon\partial_{\gamma_2}
E^{\varepsilon}\cdot\nabla\times\partial_{\gamma_2}B^{\varepsilon}dx
\nonumber\\
+&C_{\xi} \{||\nabla_x\mathbf{(I-P_2)}
\partial_{\gamma_1}g^{\varepsilon}||_{\nu}^2+||\mathbf{(I-P_2)}
\partial_{\gamma_1}g^{\varepsilon}||_{\nu}^2+||\mathbf{(I-P_2)}
\partial_{\gamma_2}g^{\varepsilon}||_{\nu}^2\}\label{e11}\\
+&\varepsilon^2||\partial_{\gamma_1}h_{e}^{\varepsilon}||^2
+||\partial_{\gamma_1}j_1^{\varepsilon}||^2+\xi||\nabla\times
\partial_{\gamma_2}j_2^{\varepsilon}||^2.\nonumber
\end{align}

Now let us go back to (\ref{abceb}). After applying (\ref{b1}) and
(\ref{e11}) for $|\gamma_1|=0,1$ and $|\gamma_2|=0$ with fixed
small $\xi$, we have
\begin{equation}
\begin{split}
C_4(&||E^{\varepsilon}||^2+||B^{\varepsilon}||^2)\leq
\frac{d}{dt}\int_{\mathbb{T}^3}\varepsilon
E^{\varepsilon}\cdot\nabla\times B^{\varepsilon}dx\\
&+\frac{d}{dt}\int_{\mathbb{T}^3}\langle\varepsilon\mathbf{(I-P_2)}
\nabla_x g^{\varepsilon},v\sqrt{\mu}\rangle\cdot \nabla
E^{\varepsilon}+\langle\varepsilon\mathbf{(I-P_2)}
g^{\varepsilon},v\sqrt{\mu}\rangle\cdot
E^{\varepsilon}dx\\
&+\frac{C}{\varepsilon^2} \{||\nabla_x\mathbf{(I-P_2)} \nabla_x
g^{\varepsilon}||_{\nu}^2+||\mathbf{(I-P_2)} \nabla_x
g^{\varepsilon}||_{\nu}^2+||\mathbf{(I-P_2)}
g^{\varepsilon}||_{\nu}^2\}\\
&+\varepsilon^2(||\nabla_x h_{e}^{\varepsilon}||^2
+||h_{e}^{\varepsilon}||^2)+||j_1^{\varepsilon}||^2+ ||\nabla
j_1^{\varepsilon}||^2+||j_2^{\varepsilon}||^2+||\nabla
j_2^{\varepsilon}||^2,
\end{split}
\end{equation}
where $0<C_4<1$ and $\varepsilon\leq 1$. This estimate with
(\ref{abceb}) finally finishes (\ref{diff}) and thus Lemma.
\end{proof}\

\begin{lem}
\label{stenergy}Assume
$f^{\varepsilon},g^{\varepsilon},E^{\varepsilon},B^{\varepsilon}$
are solutions to the system (\ref{model}), (\ref{model2}) and satisfy
(\ref{mme}). Then there exists constant $C_{1}%
\geq1$ such that the following energy estimate is valid:
\begin{equation}
\begin{split}
 \frac{d}{dt}&\{C_{1}\sum_{|\gamma|\leq N}\{||[\partial_{\gamma}
 {f^{\varepsilon}},\;\partial_{\gamma}{g^{\varepsilon}}]||^{2}+
||[\partial_{\gamma}E^{\varepsilon},\;
\partial_{\gamma}B^{\varepsilon}]||^2\}-\varepsilon
\delta
G(t)\}\\
&+\delta\sum_{|\gamma|\leq
N}\{\frac{1}{\varepsilon^{2}}||[\partial
_{\gamma}{\mathbf{\{I-P_1\}}f^{\varepsilon}},\partial_{\gamma}
{\mathbf{\{I-P_2\}}g^{\varepsilon}}]||_{\nu}^{2}+||[\partial_{\gamma
}{\mathbf{P_1}f^{\varepsilon}},\partial_{\gamma}
{\mathbf{P_2}g^{\varepsilon}}]||^{2}\}\label{ste}\\
\ \leq&\; 2C_{1}\sum_{|\gamma|\leq
N}\{\mathbf{(}\partial_{\gamma}h_1^{\varepsilon},\partial_{\gamma
}f^{\varepsilon})+\mathbf{(}\partial_{\gamma}h_2^{\varepsilon},\partial_{\gamma
}g^{\varepsilon})\}+\varepsilon^{2}\delta\sum_{|\gamma|\leq
N-1}||\partial_{\gamma}h_{||}^{\varepsilon }||^{2}\\
&+\frac{2C_1}{\varepsilon}\sum_{|\gamma|\leq
N}(||\partial_{\gamma}E^{\varepsilon}||\cdot
||\partial_{\gamma}j_1^{\varepsilon}||+||\partial_{\gamma}B^{\varepsilon}
||\cdot||\partial_{\gamma}j_2^{\varepsilon}||)\\
&+\varepsilon^{2}\delta\mathcal{A}^2+\varepsilon^2\delta(||j_1^{\varepsilon}||^2
+||\nabla j_1^{\varepsilon}||^2+||j_2^{\varepsilon}||^2+||\nabla
j_2^{\varepsilon}||^2).
\end{split}
\end{equation}
\end{lem}\

\begin{proof}
The standard $\partial_{\gamma}$ energy estimates with
(\ref{model}) and (\ref{model2}) give rise to
\begin{equation*}
\begin{split}
\frac{1}{2}\frac{d}{dt}\{||\partial_{\gamma}f^{\varepsilon}||^2
+||\partial_{\gamma}g^{\varepsilon}||^2\}+
 \frac{1}{\varepsilon^2}\{(L\partial_{\gamma}f^{\varepsilon},
 \partial_{\gamma}f^{\varepsilon})
 +(\mathcal{L}\partial_{\gamma}g^{\varepsilon},
 \partial_{\gamma}g^{\varepsilon})\}\\
-\frac{1}{\varepsilon}
 (v\sqrt{\mu}\cdot \partial_{\gamma}E^{\varepsilon},
 \partial_{\gamma}g^{\varepsilon}) =(\partial_{\gamma}
 h_1^{\varepsilon},\partial_{\gamma}f^{\varepsilon})
 +(\partial_{\gamma}h_2^{\varepsilon},\partial_{\gamma}g^{\varepsilon}).
\end{split}
\end{equation*}
Use (\ref{model2}) twice to deal with $-\frac{1}{\varepsilon}
 (v\sqrt{\mu}\cdot \partial_{\gamma}E^{\varepsilon},
 \partial_{\gamma}g^{\varepsilon})$:
\begin{equation*}
\begin{split}
-\frac{1}{\varepsilon} (v\sqrt{\mu}\cdot
\partial_{\gamma}E^{\varepsilon},\partial_{\gamma}g^{\varepsilon})=
\frac{1}{\varepsilon}\int_{\mathbb{T}^3}
\partial_{\gamma}E^{\varepsilon}\cdot\{\varepsilon \partial_t
\partial_{\gamma} E^{\varepsilon}-\nabla\times\partial_{\gamma}
B^{\varepsilon}-\partial_{\gamma}j_1^{\varepsilon}\}dx\\
=\frac{1}{2}\frac{d}{dt}\{||\partial_{\gamma}E^{\varepsilon}||^2
+||\partial_{\gamma}B^{\varepsilon}||^2\}-\frac{1}{\varepsilon}
\int_{\mathbb{T}^3}\{\partial_{\gamma}E^{\varepsilon}\cdot
\partial_{\gamma}j_1^{\varepsilon}
+\partial_{\gamma}B_R^{\varepsilon}\cdot\partial_{\gamma}
j_2^{\varepsilon}\}dx
\end{split}
\end{equation*}
But from Lemma \ref{abc}, there is a constant $C_{1}\geq1$ such
that
\[
\begin{split}
\frac{\delta}{2\varepsilon^{2}}\sum_{|\gamma|\leq
N}||[\partial_{\gamma
}{\mathbf{\{I-P_1\}}f^{\varepsilon}},&\partial_{\gamma
}{\mathbf{\{I-P_2\}}
g^{\varepsilon}}]||_{\nu}^{2}\\
\geq\frac{1}{2C_{1}}\sum_{|\gamma|\leq N}\{\delta
||[\partial_{\gamma
}{\mathbf{P_1}f^{\varepsilon}},&\partial_{\gamma}
{\mathbf{P_2}g^{\varepsilon}}]|| -\varepsilon\delta G(t)\}
-\varepsilon^{2}\delta\sum_{|\gamma|\leq
N-1}||\partial_{\gamma}h_{||}^{\varepsilon}||^{2}\\-\varepsilon^{2}
\delta\mathcal{A}^2-\varepsilon^2\delta(||j_1^{\varepsilon}||^2
+&||\nabla j_1^{\varepsilon}||^2+||j_2^{\varepsilon}||^2+||\nabla
j_2^{\varepsilon}||^2).
\end{split}
\]
By  (\ref{lower}), multiplying by $C_{1}$ and collecting terms, we
deduce our lemma.
\end{proof}\

\section{The First Order Remainder}\

In this section, we prove Theorem \ref{nslimit}. We study the
solution to the kinetic equation (\ref{kinetic}). We first
establish the spatial energy estimate for $\partial_{\gamma
}f^{\varepsilon}$ and $\partial_{\gamma }g^{\varepsilon}$.\\

\begin{lem}
\label{pure1}Let
$f^{\varepsilon},g^{\varepsilon},E^{\varepsilon},B^{\varepsilon}$
be solutions to the system (\ref{kinetic})-(\ref{kinetic1}) which
satisfy the following conservation laws:
\begin{equation}
\begin{split}
\int_{\mathbb{T}^3}
a^{\varepsilon}(t,x)dx&=\frac{\varepsilon}{2}\int_{\mathbb{T}^3}
|E^{\varepsilon}(t,x)|^2+|B^{\varepsilon}(t,x)|^2dx,\\  \label{mme1}%
\int_{\mathbb{T}^3}
b^{\varepsilon}(t,x)dx&=-{\varepsilon}\int_{\mathbb{T}^3}
E^{\varepsilon}(t,x)\times B^{\varepsilon}(t,x)dx,\\
\int_{\mathbb{T}^3}
c^{\varepsilon}(t,x)dx&=-\frac{\varepsilon}{6}\int_{\mathbb{T}^3}
|E^{\varepsilon}(t,x)|^2+|B^{\varepsilon}(t,x)|^2dx,\\
\int_{\mathbb{T}^3} d^{\varepsilon}(t,x)dx&=0\;\text{ and
}\;\int_{\mathbb{T}^3} B^{\varepsilon}(t,x)dx=0.
\end{split}
\end{equation}
Then for any instant energy functional
$\mathcal{E}_{N}(f^{\varepsilon},g^{\varepsilon},
E^{\varepsilon},B^{\varepsilon})$ where $N\geq 8$, there is a
constant $C>0$ such that
\begin{equation}
\begin{split}
 \frac{d}{dt}&\{C_{1}\{||[\partial_{\gamma}{f^{\varepsilon}},
 \partial_{\gamma}
{g^{\varepsilon}}]||^{2}+
||[\partial_{\gamma}E^{\varepsilon},\partial_{\gamma}B^{\varepsilon}]||^2\}-\varepsilon
\delta
G(t)\}\\
&+\delta\{\frac{1}{\varepsilon^{2}}||[\partial
_{\gamma}{\mathbf{\{I-P_1\}}f^{\varepsilon}},\partial_{\gamma}
{\mathbf{\{I-P_2\}}g^{\varepsilon}}]||_{\nu}^{2}+||[\partial_{\gamma
}{\mathbf{P_1}f^{\varepsilon}},\partial_{\gamma}{\mathbf{P_2}g^{\varepsilon}}]||^{2}\}\\
\leq&\; C\{\mathcal{E}_{N}(f^{\varepsilon},g^{\varepsilon},
E^{\varepsilon},B^{\varepsilon})^{\frac{1}{2}}+\mathcal{E}_{N}(f^{\varepsilon},g^{\varepsilon},
E^{\varepsilon},B^{\varepsilon})%
\}\mathcal{D}_{N}(f^{\varepsilon},g^{\varepsilon}).\label{purelen}
\end{split}
\end{equation}\
\end{lem}\

Notation: We use $\mathcal{E}_N\text{ and }\mathcal{D}_N$ instead
of $\mathcal{E}_{N}(f^{\varepsilon},g^{\varepsilon},
E^{\varepsilon},B^{\varepsilon})$ and
$\mathcal{D}_{N}(f^{\varepsilon},g^{\varepsilon})$ without
confusion.

\begin{proof}
Note that (\ref{mme1}) falls into the category of (\ref{mme})
with $\mathcal{A}=0$. We apply Lemma \ref{stenergy} with
\begin{align*}
&h_1^{\varepsilon}\equiv\frac{1}{\varepsilon}\Gamma(f^{\varepsilon}%
,f^{\varepsilon})+\frac{1}{2}E^{\varepsilon}\cdot v g^{\varepsilon}
-(E^{\varepsilon}+v\times B^{\varepsilon})\cdot\nabla_v g^{\varepsilon},\\
&h_2^{\varepsilon}\equiv\frac{1}{\varepsilon}\Gamma(g^{\varepsilon}%
,f^{\varepsilon})+\frac{1}{2}E^{\varepsilon}\cdot v
f^{\varepsilon} -(E^{\varepsilon}+v\times
B^{\varepsilon})\cdot\nabla_v f^{\varepsilon},\\
&j_1^{\varepsilon}=j_2^{\varepsilon}\equiv 0.
\end{align*}
It suffices to estimate the RHS of (\ref{ste}): we only need to
estimate
$\varepsilon^{2}||\partial_{\gamma}h_{||}^{\varepsilon}||^{2}$
and $\mathbf{(}\partial_{\gamma}h_1^{\varepsilon},\partial_{\gamma
}f^{\varepsilon})+\mathbf{(}\partial_{\gamma}h_2^{\varepsilon},\partial_{\gamma
}g^{\varepsilon})$. Terms related to collision kernel such as
$\frac{1} {\varepsilon}\Gamma(f^{\varepsilon},f^{\varepsilon})$,
$\frac{1} {\varepsilon}\Gamma(g^{\varepsilon},f^{\varepsilon})$
are bounded by $C\{(\mathcal{E}_{N})^{\frac{1}{2}}+\mathcal{E}_{N}
\}\mathcal{D}_{N}$ due to Lemma 7.1 in \cite{Guo2}. Here we show
that the ones related to the electromagnetic field also have the
same bound. First we look at the projection part
$||\partial_{\gamma}h_{||}^{\varepsilon}||^{2}$ including
$E^{\varepsilon}$ and $B^{\varepsilon}$:
\begin{equation}
\begin{split}
&\int_{\mathbb{T}^3}[\int_{\mathbb{R}^3}\partial_{\gamma}
\{\frac{1}{2}E^{\varepsilon}\cdot v g^{\varepsilon}
-(E^{\varepsilon}+v\times
B^{\varepsilon})\cdot\nabla_v g^{\varepsilon}\}\zeta dv]^2dx\\
&+\int_{\mathbb{T}^3}[\int_{\mathbb{R}^3}\partial_{\gamma}
\{\frac{1}{2}E^{\varepsilon}\cdot v f^{\varepsilon}
-(E^{\varepsilon}+v\times B^{\varepsilon})\cdot\nabla_v
f^{\varepsilon}\}\zeta dv]^2dx.\label{h1}
\end{split}
\end{equation}
We will only estimate the first term. Let $\gamma=\gamma_1
+\gamma_2$ where $|\gamma|\leq N-1$.
\begin{equation*}
(I)\equiv\int_{\mathbb{T}^3}[\int_{\mathbb{R}^3}
\{\frac{1}{2}\partial_{\gamma_1}E^{\varepsilon}\cdot v
\partial_{\gamma_2}g^{\varepsilon} -(\partial_{\gamma_1}E^{\varepsilon}
+v\times\partial_{\gamma_1}
B^{\varepsilon})\cdot\nabla_v\partial_{\gamma_2}
g^{\varepsilon}\}\zeta dv]^2dx
\end{equation*}
If $|\gamma_1|\leq|\gamma_2|$, we take the sup of
$\partial_{\gamma_1}E^{\varepsilon},\partial_{\gamma_1}B^{\varepsilon}$:
\begin{align*}
(I)\leq \int_{\mathbb{T}^3}|\partial_{\gamma_1}&E^{\varepsilon}|^2
|\int_{\mathbb{R}^3} v
\partial_{\gamma_2}g^{\varepsilon}\zeta dv|^2 dx+
\int_{\mathbb{T}^3}|\partial_{\gamma_1}E^{\varepsilon}|^2
|\int_{\mathbb{R}^3}
\partial_{\gamma_2}g^{\varepsilon}\nabla_v\zeta dv|^2 dx\\ +
\int_{\mathbb{T}^3}&|\partial_{\gamma_1}B^{\varepsilon}|^2
|\int_{\mathbb{R}^3}v\times
\partial_{\gamma_2}g^{\varepsilon}\nabla_v\zeta dv|^2 dx\\
\leq \mathcal{E}_{|\gamma_1|+2}\{&\int_{\mathbb{T}^3}
|\int_{\mathbb{R}^3} v
\partial_{\gamma_2}g^{\varepsilon}\zeta dv|^2 dx+\int_{\mathbb{T}^3}
|\int_{\mathbb{R}^3}
\partial_{\gamma_2}g^{\varepsilon}\nabla_v\zeta dv|^2 dx\\
&+\int_{\mathbb{T}^3}|\int_{\mathbb{R}^3}v\times
\partial_{\gamma_2}g^{\varepsilon}\nabla_v\zeta dv|^2 dx\}
\;\text{ (Sobolev imbedding theorem)}\\
\leq \mathcal{E}_{|\gamma_1|+2}(&\int_{\mathbb{R}^3}|v\zeta|^2+
|\nabla_v\zeta|^2+|v\nabla_v\zeta|^2 dv)\int_{\mathbb{R}^3}
\int_{\mathbb{R}^3}|\partial_{\gamma_2}g^{\varepsilon}|^2 dxdv\;
\text{ (H}\ddot{o}\text{lder ineq)} \\
\leq \mathcal{E}_{|\gamma_1|+2}\;&\mathcal{D}_{|\gamma_2|}\;\text{
(Note that }\zeta\text{ decays exponentially) }\\
\leq \mathcal{E}_{N}\;\mathcal{D}_{N}&\;\;\;(N\geq8);
\end{align*}
on the other hand, if $|\gamma_2|<|\gamma_1|$, we take the sup of
$\partial_{\gamma_2}g^{\varepsilon}$:
\begin{align*}
(I)\leq
\sup_{x,v}|\partial_{\gamma_2}g^{\varepsilon}|\;\{&\int_{\mathbb{T}^3}
|\partial_{\gamma_1}E^{\varepsilon}|^2 |\int_{\mathbb{R}^3} v
\zeta dv|^2 dx+
\int_{\mathbb{T}^3}|\partial_{\gamma_1}E^{\varepsilon}|^2
|\int_{\mathbb{R}^3} \nabla_v\zeta dv|^2 dx\\
&+ \int_{\mathbb{T}^3}|\partial_{\gamma_1}B^{\varepsilon}|^2
|\int_{\mathbb{R}^3}v\times \nabla_v\zeta dv|^2 dx\}\\
\leq \mathcal{D}_{|\gamma_2|+4}\;\mathcal{E}_{|\gamma_1|}&\\
\leq \mathcal{E}_{N}\;\mathcal{D}_{N}\;\;\;(N&\geq8).
\end{align*}
Thus we conclude that for $|\gamma|\leq N-1$,
\begin{equation}
\varepsilon^{2}||\partial_{\gamma}h_{||}^{\varepsilon}||^{2}\leq
C\mathcal{E}_{N}\;\mathcal{D}_{N}.\label{gammah}
\end{equation}
Now we turn to
$\mathbf{(}\partial_{\gamma}h_1^{\varepsilon},\partial_{\gamma
}f^{\varepsilon})+\mathbf{(}\partial_{\gamma}h_2^{\varepsilon},\partial_{\gamma
}g^{\varepsilon})$ for $|\gamma|\leq N$.
\begin{align}
&\int_{\mathbb{T}^3}\int_{\mathbb{R}^3}\partial_{\gamma}
\{\frac{1}{2}E^{\varepsilon}\cdot v g^{\varepsilon}
-(E^{\varepsilon}+v\times B^{\varepsilon})\cdot\nabla_v
g^{\varepsilon}\}\cdot
\partial_{\gamma}f^{\varepsilon} dvdx\nonumber\\
&+\int_{\mathbb{T}^3}\int_{\mathbb{R}^3}\partial_{\gamma}
\{\frac{1}{2}E^{\varepsilon}\cdot v f^{\varepsilon}
-(E^{\varepsilon}+v\times B^{\varepsilon})\cdot\nabla_v
f^{\varepsilon}\}\cdot\partial_{\gamma}g^{\varepsilon} dvdx\nonumber\\
=&\sum_{\gamma_1\neq0}\int_{\mathbb{T}^3}\int_{\mathbb{R}^3}
\{\frac{1}{2}\partial_{\gamma_1}E^{\varepsilon}\cdot v
\partial_{\gamma_2}g^{\varepsilon} -(\partial_{\gamma_1}
E^{\varepsilon}+v\times\partial_{\gamma_1}
B^{\varepsilon})\cdot\nabla_v
\partial_{\gamma_2}g^{\varepsilon}\}\cdot
\partial_{\gamma}f^{\varepsilon} dvdx\label{gammafg}\\
+&\sum_{\gamma_1\neq0}\int_{\mathbb{T}^3}\int_{\mathbb{R}^3}
\{\frac{1}{2}\partial_{\gamma_1}E^{\varepsilon}\cdot v
\partial_{\gamma_2}f^{\varepsilon} -(\partial_{\gamma_1}
E^{\varepsilon}+v\times\partial_{\gamma_1}
B^{\varepsilon})\cdot\nabla_v
\partial_{\gamma_2}f^{\varepsilon}\}\cdot
\partial_{\gamma}g^{\varepsilon} dvdx\nonumber\\
+&\int_{\mathbb{T}^3}\int_{\mathbb{R}^3} E^{\varepsilon}\cdot v
\partial_{\gamma}f^{\varepsilon}
\partial_{\gamma}g^{\varepsilon}dvdx -\int_{\mathbb{T}^3}\int_{\mathbb{R}^3}
(E^{\varepsilon}+v\times B^{\varepsilon})\cdot\nabla_v
(\partial_{\gamma}f^{\varepsilon}
\partial_{\gamma}g^{\varepsilon}) dvdx\nonumber
\end{align}
The worst last $(|\gamma|+1)$-th derivative term of  becomes zero
after the integration by parts. As we note that $\nu(v)\sim
(1+|v|)$ and $\int\int v f g
 dv dx \leq C ||f||_{\nu}||g||_{\nu}$, we apply the Sobolev
imbedding theorem to other terms like the previous argument and
then we can get the desired bound
$(\mathcal{E}_N)^{\frac{1}{2}}\mathcal{D}_N$.
\end{proof}\

In order to prove Theorem \ref{nslimit}, it remains to estimate
the velocity derivatives.
\begin{proof}\textit{(of Theorem \ref{nslimit}:) }We notice
that for the hydrodynamic part
$[{\mathbf{P_1}f^{\varepsilon}},{\mathbf{P_2}g^{\varepsilon}}],$
\begin{equation}
||[\partial_{\gamma}^{\beta}{\mathbf{P_1}f^{\varepsilon}},
\partial_{\gamma}^{\beta}{\mathbf{P_2}g^{\varepsilon}}]||\leq C||[\partial
_{\gamma}{\mathbf{P_1}f^{\varepsilon}},\partial_{\gamma}
{\mathbf{P_2}g^{\varepsilon}}]|| \label{novhydro}
\end{equation}
which has been estimated in Lemma \ref{pure1}. It suffices to
estimate the remaining microscopic part
\[
[\partial_{\gamma_1}^{\beta}{\mathbf{\{I-P_1\}}f^{\varepsilon}}
,\partial_{\gamma_1}^{\beta}{\mathbf{\{I-P_2\}}g^{\varepsilon}}]
\]
for $|\gamma_1|+|\beta|\leq N$, $\;|\gamma_1|\leq N-1$ and
$\beta>0.$ We take $\partial_{\gamma_1}^{\beta}$ of equations
(\ref{kinetic}) and sum over $|\gamma_1|+|\beta|\leq N$ to get
\begin{align}
&\partial_{t}\partial_{\gamma_1}^{\beta}\mathbf{(I-P_1)}f^{\varepsilon}+\frac
{1}{\varepsilon}v\cdot\nabla_{x}\partial_{\gamma_1}^{\beta}\mathbf{(I-P_1)}%
f^{\varepsilon}+\frac{1}{\varepsilon^{2}}\partial_{\gamma_1}^{\beta
}L\mathbf{(I-P_1)}f^{\varepsilon}\nonumber\\
&+\{\partial_{t}\partial_{\gamma_1}^{\beta}\mathbf{P_1}f^{\varepsilon}+\frac
{1}{\varepsilon}v\cdot\nabla_{x}\partial_{\gamma_1}^{\beta}\mathbf{P_1}%
f^{\varepsilon}+\frac{1}{\varepsilon}\binom{\beta_{1}}{\beta}\partial
^{\beta_{1}}v\cdot\nabla_{x}\partial_{\gamma_1}^{\beta-\beta_{1}}f^{\varepsilon
}\}\nonumber\\
& =\partial_{\gamma_1}^{\beta}\{\frac{1}{\varepsilon}\Gamma(f^{\varepsilon}%
,f^{\varepsilon})+\frac{1}{2}E^{\varepsilon}\cdot v
g^{\varepsilon}
-(E^{\varepsilon}+v\times B^{\varepsilon})\cdot\nabla_v g^{\varepsilon}\}
;\nonumber\\
&\partial_{t}\partial_{\gamma_1}^{\beta}\mathbf{(I-P_2)}g^{\varepsilon}
+\frac
{1}{\varepsilon}v\cdot[\nabla_{x}\partial_{\gamma_1}^{\beta}\mathbf{(I-P_2)}%
g^{\varepsilon}-\partial^{\beta}(\sqrt{\mu})\partial_{\gamma_1}E^
{\varepsilon}]+\frac{1}{\varepsilon^{2}}\partial_{\gamma_1}^{\beta
}\mathcal{L}\mathbf{(I-P_2)}g^{\varepsilon}\label{vderivative}\\
&+\{\partial_{t}\partial_{\gamma_1}^{\beta}\mathbf{P_2}g^{\varepsilon}
+\frac
{1}{\varepsilon}v\cdot\nabla_{x}\partial_{\gamma_1}^{\beta}\mathbf{P_2}%
g^{\varepsilon}+\frac{1}{\varepsilon}\binom{\beta_{1}}{\beta}\partial
^{\beta_{1}}v\cdot[\nabla_{x}\partial_{\gamma_1}^{\beta-\beta_{1}}g^{\varepsilon
}-\partial^{\beta-\beta_1}(\sqrt{\mu})\partial_{\gamma_1}E^{\varepsilon}]\}
\nonumber\\
& =\partial_{\gamma_1}^{\beta}\{\frac{1}{\varepsilon}\Gamma(g^{\varepsilon}%
,f^{\varepsilon})+\frac{1}{2}E^{\varepsilon}\cdot v
f^{\varepsilon} -(E^{\varepsilon}+v\times
B^{\varepsilon})\cdot\nabla_v f^{\varepsilon}\},\nonumber
\end{align}
where $|\beta_{1}|=1.$ We illustrate the estimate on only
$\mathbf{\{I-P_2\}}g^{\varepsilon}$ and the other one can be done
in the same way. Also we will use many results from \cite{Guo2}.
Taking the inner product with $\partial_{\gamma_1 }^{\beta}
{\mathbf{\{I-P_2\}}g^{\varepsilon}}$, we get
\begin{equation}
\begin{split}
 \frac{1}{2}&\frac{d}{dt}\{||\partial_{\gamma_1}^{\beta}
 {\mathbf{(I-P_2)}g^{\varepsilon}}||^{2}\}
+\frac{1}{\varepsilon^{2}} \mathbf{{\huge (}}
\partial_{\gamma_1}^{\beta}
{\mathcal{L}\mathbf{(I-P_2)}g^{\varepsilon}},\partial
_{\gamma_1}^{\beta}{\mathbf{(I-P_2)}g^{\varepsilon}})\\
-&\frac
{1}{\varepsilon}(v\cdot\partial^{\beta}(\sqrt{\mu})\partial_{\gamma_1}E^
{\varepsilon},\partial_{\gamma_1}^{\beta}
\mathbf{(I-P_2)}g^{\varepsilon})+(\partial_{t}
\partial_{\gamma_1}^{\beta}
{\mathbf{P_2}g^{\varepsilon}}+\frac
{1}{\varepsilon}v\cdot\nabla_{x}\partial_{\gamma_1}^{\beta}
{\mathbf{P_2}g^{\varepsilon}}\\
+&\frac{1}{\varepsilon}\binom{\beta_{1}}{\beta}\partial^{\beta_{1}}v\cdot\partial_{\gamma_1}^{\beta-\beta_{1}}\nabla_{x}%
\partial_{\gamma_1}{g^{\varepsilon}},\partial_{\gamma_1}^{\beta}%
{\mathbf{(I-P_2)}g^{\varepsilon}})\\
-&(\frac{1}{\varepsilon}\binom{\beta_{1}}{\beta}\partial
^{\beta_{1}}v\cdot\partial^{\beta-\beta_1}(\sqrt{\mu})\partial_{\gamma_1}
E^{\varepsilon},\partial_{\gamma_1}^{\beta}
\mathbf{(I-P_2)}g^{\varepsilon})\\
 =\;&\frac{1}{\varepsilon}\mathbf{{\huge
(}}\partial_{\gamma_1}^{\beta}
{\Gamma(g^{\varepsilon},f^{\varepsilon})},\partial_{\gamma_1}^{\beta
}{\mathbf{(I-P_2)}g^{\varepsilon}})\\
+&\mathbf{{\huge (}}\partial_{\gamma_1}^{\beta
}\{\frac{1}{2}E^{\varepsilon}\cdot v f^{\varepsilon}
-(E^{\varepsilon}+v\times B^{\varepsilon})\cdot\nabla_v
f^{\varepsilon}\},\partial_{\gamma_1}^{\beta }
{\mathbf{(I-P_2)}g^{\varepsilon}}). \label{hardv}\\
\end{split}
\end{equation}
Here we only estimate the terms including $E^{\varepsilon}$ or
$B^{\varepsilon}$. We utilize the field estimate from the
previous section. Other terms can be estimated in the same way as
presented in \cite{Guo2}. See the proof of Theorem 2.2 there
 (p.38-40). The third term on the LHS of (\ref{hardv}) can be
 taken care of as following:
\begin{align*}
\frac
{1}{\varepsilon}&(v\cdot\partial^{\beta}(\sqrt{\mu})\partial_{\gamma_1}E^
{\varepsilon},\partial_{\gamma_1}^{\beta}
\mathbf{(I-P_2)}g^{\varepsilon})\\
\leq&\; C_{\xi}||\partial_{\gamma_1}E^{\varepsilon}||^2+\frac{\xi}
{\varepsilon^2}||\partial_{\gamma_1}^{\beta}
\mathbf{(I-P_2)}g^{\varepsilon}||_{\nu}^2\;\;(\xi \text{ is a
small
positive fixed number})\\
\leq&\;\widetilde{C}_{\xi}\frac{d}{dt}\int_{\mathbb{T}^3}\langle\varepsilon\mathbf{(I-P_2)}
\partial_{\gamma_1}g^{\varepsilon},v\sqrt{\mu}\rangle\cdot
\partial_{\gamma_1}E^{\varepsilon}+C_2\varepsilon\partial_{\gamma_2}
E^{\varepsilon}\cdot\nabla\times\partial_{\gamma_2}B^{\varepsilon}dx\\
+&\widetilde{C}_{\xi}\{||\nabla_x\mathbf{(I-P_2)}
\partial_{\gamma_1}g^{\varepsilon}||_{\nu}^2+||\mathbf{(I-P_2)}
\partial_{\gamma_1}g^{\varepsilon}||_{\nu}^2+||\mathbf{(I-P_2)}
\partial_{\gamma_2}g^{\varepsilon}||_{\nu}^2\}+\mathcal{E}_N\mathcal{D}_N\\
+&\frac{\xi} {\varepsilon^2}||\partial_{\gamma_1}^{\beta}
\mathbf{(I-P_2)}g^{\varepsilon}||_{\nu}^2\;\;\;(\text{by }
(\ref{e11}) \text{ and } (\ref{gammah}))
\end{align*}
where $|\gamma_2|\leq N-2$. One gets the same estimate for
$(\frac{1}{\varepsilon}\binom{\beta_{1}}{\beta}\partial
^{\beta_{1}}v\cdot\partial^{\beta-\beta_1}(\sqrt{\mu})\partial_{\gamma_1}
E^{\varepsilon},\\ \partial_{\gamma_1}^{\beta}
\mathbf{(I-P_2)}g^{\varepsilon})$, since it is bounded by
$C_{\xi}||\partial_{\gamma_1}E^{\varepsilon}||^2+\frac{\xi}
{\varepsilon^2}||\partial_{\gamma_1}^{\beta}
\mathbf{(I-P_2)}g^{\varepsilon}||_{\nu}^2$. Define
$\widetilde{G}(t)$ by
\[
\widetilde{G}(t)\equiv
\widetilde{C}_{\xi}\int_{\mathbb{T}^3}\langle\mathbf{(I-P_2)}
\partial_{\gamma_1}g^{\varepsilon},v\sqrt{\mu}\rangle\cdot
\partial_{\gamma_1}E^{\varepsilon}+C_2\partial_{\gamma_2}
E^{\varepsilon}\cdot\nabla\times\partial_{\gamma_2}B^{\varepsilon}dx.
\]

Next, the last term in (\ref{hardv}) will be treated. We decompose
$f^{\varepsilon}\text{ and }g^{\varepsilon}$ into
$\mathbf{P_1}f^{\varepsilon}+\mathbf{(I-P_1)}f^{\varepsilon}$ and
$\mathbf{P_2}g^{\varepsilon}+\mathbf{(I-P_2)}g^{\varepsilon}$
respectively:
\begin{equation}
\begin{split}
&\mathbf{{\huge (}}\partial_{\gamma_1}^{\beta }
\{\frac{1}{2}E^{\varepsilon}\cdot v f^{\varepsilon}
-(E^{\varepsilon}+v\times B^{\varepsilon})\cdot\nabla_v
f^{\varepsilon}\},\partial_{\gamma_1}^{\beta }
{\mathbf{(I-P_2)}g^{\varepsilon}})\\
&=\mathbf{{\huge (}}\partial_{\gamma_1}^{\beta }
\{\frac{1}{2}E^{\varepsilon}\cdot v \mathbf{P_1}f^{\varepsilon}
-(E^{\varepsilon}+v\times B^{\varepsilon})\cdot\nabla_v
\mathbf{P_1}f^{\varepsilon}\},\partial_{\gamma_1}^{\beta }
{\mathbf{(I-P_2)}g^{\varepsilon}})\\
&\;\;+\mathbf{{\huge (}}\partial_{\gamma_1}^{\beta }
\{\frac{1}{2}E^{\varepsilon}\cdot v
\mathbf{(I-P_1)}f^{\varepsilon} -(E^{\varepsilon}+v\times
B^{\varepsilon})\cdot\nabla_v
\mathbf{(I-P_1)}f^{\varepsilon}\},\partial_{\gamma_1}^{\beta }
{\mathbf{(I-P_2)}g^{\varepsilon}})\label{ebsplit}\\
&\equiv(I)+(II)
\end{split}
\end{equation}
The both terms can be computed in the same spirit as in
(\ref{gammafg}). The main concern is that the number of
derivatives could be $N+1$ for the worst case due to the Vlasov
term $[\nabla_v{g^{\varepsilon}},\nabla_v{f^{\varepsilon}}]$. As
for $(I)$, we recall that the hydrodynamic parts are not affected
by the velocity derivative as noted in (\ref{novhydro}) and so we
can deduce that it is bounded by
$(\mathcal{E}_N)^{\frac{1}{2}}\mathcal{D}_N$. For the second
term  $(II)$ , the $(N+1)$-th order derivative term together with
the similar term stemming from the first equation in
(\ref{vderivative}) is gone by the integration by parts:
\[
\begin{split}
&((E^{\varepsilon}+v\times
B^{\varepsilon})\cdot\nabla_v\partial_{\gamma_1}^{\beta }
\mathbf{(I-P_1)}g^{\varepsilon},\partial_{\gamma_1}^{\beta }
{\mathbf{(I-P_2)}f^{\varepsilon}})\\
&+ ((E^{\varepsilon}+v\times
B^{\varepsilon})\cdot\nabla_v\partial_{\gamma_1}^{\beta }
\mathbf{(I-P_1)}f^{\varepsilon},\partial_{\gamma_1}^{\beta }
{\mathbf{(I-P_2)}g^{\varepsilon}})=0.
\end{split}
\]
and hence it is also bounded by
$(\mathcal{E}_N)^{\frac{1}{2}}\mathcal{D}_N$.

Combining the above estimates with Theorem 2.2 in \cite{Guo2}, we
obtain for $|\gamma_1|\leq N-1$ and $|\beta|+|\gamma_1|=N$,
\begin{align*}
  \frac{d}{dt}&\{\frac{1}{2}||[\partial_{\gamma_1}^{\beta}
{\mathbf{\{I-P_1\}}f^{\varepsilon}},\partial_{\gamma_1}^{\beta}
{\mathbf{\{I-P_2\}}g^{\varepsilon}}]||^{2}\}\\
&\;\;\;\;+\frac{\delta}{4\varepsilon^{2}}||[\partial_{\gamma_1}^{\beta}
{\mathbf{\{I-P_1\}}f^{\varepsilon}},\partial_{\gamma_1}^{\beta}
{\mathbf{\{I-P_2\}}g^{\varepsilon}}]||_{\nu}^{2}\\
 \leq\;&
C\{||[\partial_{\gamma_1}{f^{\varepsilon}},\partial_{\gamma_1}{g^{\varepsilon}}]||_{\nu}^{2}
+\frac{1}{\varepsilon^2}(||\nabla_x\mathbf{(I-P_2)}
\partial_{\gamma_1}g^{\varepsilon}||_{\nu}^2+||\mathbf{(I-P_2)}
\partial_{\gamma_1}g^{\varepsilon}||_{\nu}^2)\}\\
&\;\;\;\;+\varepsilon\frac{d}{dt}\widetilde{G} +C\{(\mathcal{E}
_{N})^{\frac{1}{2}}+\mathcal{E}
_{N}+\varepsilon^{2}\}\mathcal{D}_{N}.
\end{align*}
Multiplying the above by factor $4$ and adding a large multiple
$K$ of (\ref{purelen}), we get the following for $|\gamma|\leq N$,
 $|\beta|+|\gamma_1|\leq N$ and $|\gamma_1|\leq N-1$,
\begin{equation}
\begin{split}
  \frac{d}{dt}&[
K\{C_{1}\{||[\partial_{\gamma}{f^{\varepsilon}},\partial_{\gamma}
{g^{\varepsilon}}]||^{2}+ (||\partial_{\gamma}E^{\varepsilon}||^2+
||\partial_{\gamma}B^{\varepsilon}||^2)\}-\varepsilon \delta
G(t)\}\\&+2||[\partial_{\gamma_1}^{\beta
}{\mathbf{\{I-P_1\}}f^{\varepsilon}},\partial_{\gamma_1}^{\beta
}{\mathbf{\{I-P_2\}}g^{\varepsilon}}
]||^{2}-4\varepsilon\widetilde{G}(t)]
+\delta\mathcal{D}_{N}\\
 \leq\;&C_{K}\{(\mathcal{E}_{N})^{\frac{1}{2}}+\mathcal{E}
_{N}+\varepsilon^{2}\}\mathcal{D}_{N}.
\label{hardenergy}%
\end{split}
\end{equation}
Notice that
\begin{align}
G(t)\leq C\sum_{|\gamma|\leq
N}\{||[\partial_{\gamma}{\mathbf{P_1}f^{\varepsilon
}},\partial_{\gamma}{\mathbf{P_2}g^{\varepsilon
}}]||+||\partial_{\gamma}E^{\varepsilon}||
+||\partial_{\gamma}B^{\varepsilon}||\}\cdot
\{||[\partial_{\gamma}{\mathbf{P_1}f^{\varepsilon
}},\partial_{\gamma}{\mathbf{P_2}g^{\varepsilon
}}]||&\nonumber\\
+||[\partial_{\gamma}{\mathbf{\{I-P_1\}}f^{\varepsilon
}},\partial_{\gamma}{\mathbf{\{I-P_2\}}g^{\varepsilon
}}]||+||\partial_{\gamma}E^{\varepsilon}||
+||\partial_{\gamma}B^{\varepsilon}||\}.& \label{gproduct}%
\end{align}
And $\widetilde{G}(t)$ has the same bound. We thus can redefine
an instant energy by
\begin{align}
\mathcal{E}_{N}(f^{\varepsilon},g^{\varepsilon},E^{\varepsilon},B^{\varepsilon})
\equiv
K\{C_{1}\{||[\partial_{\gamma}{f^{\varepsilon}},\partial_{\gamma}
{g^{\varepsilon}}]||^{2}+ (||\partial_{\gamma}E^{\varepsilon}||^2+
||\partial_{\gamma}B^{\varepsilon}||^2)\}-\varepsilon \delta
G(t)\}&\nonumber\\+2||[\partial_{\gamma_1}^{\beta
}{\mathbf{\{I-P_1\}}f^{\varepsilon}},\partial_{\gamma_1}^{\beta
}{\mathbf{\{I-P_2\}}g^{\varepsilon}}
]||^{2}-4\varepsilon\widetilde{G}(t),\;\;\;\;\;\;\;\;\;\;&\label{hardin}%
\end{align}
for $\varepsilon$ sufficiently small. With such a small
$\varepsilon>0,$ adjusting the constants in (\ref{hardenergy})
and applying a standard continuity argument, we thus deduce
(\ref{en}) by letting $\mathcal{E}_{N}$ sufficiently small
initially. Indeed, our argument is still valid for
$\varepsilon=1$, since we can choose $K, \;C_1$ large enough. And
that guarantees the global existence of the
Vlasov-Maxwell-Boltzmann system.

Since no estimates for the highest $N$-th derivatives of
$E^{\varepsilon}$, $B^{\varepsilon}$ are not available as we can
see that full (up to $N$-th) derivatives of $g^{\varepsilon}$
 are needed in the RHS of (\ref{e11}) to control
 the $(N-1)$-th derivative of $E^{\varepsilon}$,  and hence the dissipation
rate is weaker than the instant energy, we do not expect
exponential decay unlike the pure Boltzmann case. To obtain a
decay rate (\ref{decayrate}), we use an interpolation argument
presented in \cite{St}. We point out that our method is slightly
easier because we do not have to deal with the time derivative.
The key idea is to get a bound for the $N$-th derivatives
$\partial_{\gamma}E^{\varepsilon}$ and
$\partial_{\gamma}B^{\varepsilon}$ with the bound of higher
energy. Let $k$ be a positive integer and $|\gamma|\leq N$. By an
interpolation between Sobolev spaces $H^{N-1}$ and $H^{N+k}$, we
get
\begin{equation*}
||\partial_{\gamma}E^{\varepsilon}||^2+||\partial_{\gamma}B^{\varepsilon}||^2
\leq C
\{||E^{\varepsilon}||_{H^{N-1}}^{\frac{2k}{k+1}}+||B^{\varepsilon}||
_{H^{N-1}}^{\frac{2k}{k+1}}\}\{||E^{\varepsilon}||_{H^{N+k}}
^{\frac{2}{k+1}}+||B^{\varepsilon}|| _{H^{N+k}}^{\frac{2}{k+1}}\}.
\end{equation*}
Due to (\ref{enbound}), (\ref{b1}) and (\ref{e11}) we have
\begin{equation*}
\begin{split}
||E^{\varepsilon}||_{H^{N}}^{\frac{2k+2}{k}}+||B^{\varepsilon}||
_{H^{N}}^{\frac{2k+2}{k}}&\leq
C_{N+k}\mathcal{E}_{N+k}(0)^{\frac{1}{k}}\{||E^{\varepsilon}||
_{H^{N-1}}^{2}+||B^{\varepsilon}|| _{H^{N-1}}^{2}\}\\
&\leq
C_{N+k}\mathcal{E}_{N+k}(0)^{\frac{1}{k}}\{\varepsilon\frac{d}{dt}
\widetilde{\widetilde{G}}+\mathcal{D}_N\},
\end{split}
\end{equation*}
where $\varepsilon$ sufficiently small and
\begin{equation}
\widetilde{\widetilde{G}}(t)\equiv \int_{\mathbb{T}^3}
\langle\mathbf{(I-P_2)}
\partial_{\gamma_1}g^{\varepsilon},v\sqrt{\mu}\rangle\cdot
\partial_{\gamma_1}E^{\varepsilon}+\partial_{\gamma_2}
E^{\varepsilon}\cdot\nabla\times\partial_{\gamma_2}B^{\varepsilon}dx,
\label{ggbound}
\end{equation}
up to constant multiple. In particular, we have
$|\widetilde{\widetilde{G}}(t)|\leq \mathcal{E}_N$. Noting that
the other part of $\mathcal{E}_N$, i.e. the nonelectromagnetic
part is bounded by $\mathcal{D}_N$, a lower bound for
$\mathcal{D}_N$ can be given as following:
\[
C_{N,k}\mathcal{E}_{N+k}(0)^{-\frac{1}{k}}\mathcal{E}_N^{\frac{k+1}{k}}-
\varepsilon\frac{d}{dt} \widetilde{\widetilde{G}}\leq
\mathcal{D}_N.
\]
It follows that
\begin{equation}
\frac{d}{dt}\{\mathcal{E}_N-
\varepsilon\widetilde{\widetilde{G}}\}+C_{N,k}\mathcal{E}_{N+k}(0)^{-\frac{1}{k}}
\mathcal{E}_N^{\frac{k+1}{k}}\leq 0.\label{egd}
\end{equation}
Letting $\widetilde{\mathcal{E}}_N=\mathcal{E}_N-
\varepsilon\widetilde{\widetilde{G}},$ we get
$\frac{1}{C}\mathcal{E}_N\leq \widetilde{\mathcal{E}}_N\leq C
\mathcal{E}_N$ for some $C>1$ and thus (\ref{egd}) becomes
\[
\frac{d}{dt}\widetilde{\mathcal{E}}_N+
C_{N,k}\mathcal{E}_{N+k}(0)^{-\frac{1}{k}}
\widetilde{\mathcal{E}}_N^{\frac{k+1}{k}}\leq 0.
\]
After dividing this by $\widetilde{\mathcal{E}}_N^{\frac{k+1}{k}}$
and integrating it over $[0,t]$, we get (\ref{decayrate}).
\end{proof}\

\section{High Order Remainder}\

In this section, we establish Theorem \ref{high} for $n\geq2.$ As
noted in \cite{Guo2}, the remainder equations (\ref{remainderf})
and (\ref{remaindereb})  contain singular terms such as zeroth
order terms or first order terms in $\varepsilon$ which make it
hard to apply the energy estimate directly. The difficulty is
resolved by introducing new unknowns $f_{R}^{\varepsilon},$
$g_{R}^{\varepsilon},$ $E_{R}^{\varepsilon}$ and
$B_{R}^{\varepsilon},$ which can be obtained from further
construction of ($n+1$)-th coefficients $f_{n+1}, g_{n+1},
E_{n+1},B_{n+1}$ and ($n+2$)-th coefficients $f_{n+2}, g_{n+2},
E_{n+2},B_{n+2}$ in the diffusive expansion (\ref{1/e}). The
compatibility conditions (\ref{key}) and (\ref{keyeb}) will then
eliminate such severe singularity in $\varepsilon.$

We reformulate the problem. Given $f_{1},...,f_{n}$;
$g_{1},...,g_{n}$; $E_{1},...,E_{n}$; $B_{1},...,B_{n}$ determined
by the initial data, we further construct artificial $f_{n+1},
g_{n+1}, E_{n+1},B_{n+1}$ and $f_{n+2}, g_{n+2}, E_{n+2},B_{n+2}$
by letting $m=n+1$ and $m=n+2$ in Theorem \ref{cecoe} with
constant initial data satisfying
\[
u_{m}^{0}\equiv
-\frac{1}{|\mathbb{T}^3|}\sum_{\substack{i+j=m\\i,j\geq1}}\int_{\mathbb{T}^{3}}
E_i^0(x)\times B_j^0(x)\;dx;\;\;\;\;\;\;\;\sigma_{m}^0\equiv 0;
\]
\[
\theta_{m}^0\equiv-\frac{1}{3|\mathbb{T}^3|}
\sum_{\substack{i+j=m\\i,j\geq1}}\int_{\mathbb{T}^{3}}
E_i^0(x)\cdot E_j^0(x)+B_i^0(x)\cdot B_j^0(x)\;dx.\label{ic}
\]
Let us  introduce  new unknowns $f_{R}^{\varepsilon},$
$g_{R}^{\varepsilon},$ $E_{R}^{\varepsilon}$ and
$B_{R}^{\varepsilon}$ such that
\begin{equation}
\begin{split}
f_R^{\varepsilon}&\equiv f_{n}^{\varepsilon}-f_{n}-\varepsilon f_{n+1}%
-\varepsilon^{2}f_{n+2},\;\;\label{ge}
E_R^{\varepsilon}\equiv E_{n}^{\varepsilon}-E_{n}-\varepsilon E_{n+1}%
-\varepsilon^{2}E_{n+2},\\
g_R^{\varepsilon}&\equiv g_{n}^{\varepsilon}-g_{n}-\varepsilon g_{n+1}%
-\varepsilon^{2}g_{n+2},\;\;
B_R^{\varepsilon}\equiv B_{n}^{\varepsilon}-B_{n}-\varepsilon B_{n+1}%
-\varepsilon^{2}B_{n+2}.
\end{split}
\end{equation}
Plugging (\ref{ge}) into the remainder equations
(\ref{remainderf}), (\ref{remaindereb}) and applying (\ref{key})
up to $m=n$ as well as (\ref{helec}), (\ref{hmag}) for $m=n+2$ we
get for $n\geq2$,\
\begin{equation}
\begin{split}
&\partial_{t}f_R^{\varepsilon}+\frac{v}{\varepsilon}\cdot\nabla_{x}%
f_R^{\varepsilon}+\frac{1}{\varepsilon^{2}}Lf_R^{\varepsilon}\equiv
h_1^{\varepsilon},\\
&\partial_{t}g_R^{\varepsilon}+\frac{v}{\varepsilon}\cdot\nabla_{x}%
(g_R^{\varepsilon}-\sqrt{\mu}E_R^{\varepsilon})+
\frac{1}{\varepsilon^{2}}\mathcal{L}g_R^{\varepsilon}\equiv
h_2^{\varepsilon},\label{highfg}
\end{split}
\end{equation}
\begin{equation}
\begin{split}
&\;\;\;\;\;\varepsilon \partial_t E_R^{\varepsilon}-\nabla\times
B_R^{\varepsilon}=-\int_{\mathbb{R}^3} v
g_R^{\varepsilon}\sqrt{\mu}dv-
\varepsilon^3\partial_t E_{n+2},\;\;\nabla\cdot B_R^{\varepsilon}=0,\\
&\;\;\;\;\;\varepsilon \partial_t B_R^{\varepsilon}+\nabla\times
E_R^{\varepsilon}=-\varepsilon^3\partial_t B_{n+2},\;\;\nabla\cdot
E_R^{\varepsilon}=\int_{\mathbb{R}^3}
g_R^{\varepsilon}\sqrt{\mu}dv,\label{higheb}
\end{split}
\end{equation}
where
\[
\begin{split}
{h_1^{\varepsilon}}\;\equiv{\;h_1(f_R^{\varepsilon})+h_1(f)+h_1(E_R^{\varepsilon},B_R^{\varepsilon})
+h_1(E,B)},\\
{h_2^{\varepsilon}}\;\equiv{\;h_2(g_R^{\varepsilon})+h_2(g)+h_2(E_R^{\varepsilon},B_R^{\varepsilon})
+h_2(E,B)},
\end{split}
\]
and each $h_i$ is given as follows:
\begin{equation}
\begin{split}
h_1(f_R^{\varepsilon})&\equiv\varepsilon^{n-2}
\Gamma(f_R^{\varepsilon},f_R^{\varepsilon})+\sum_{i=1}%
^{n+2}\varepsilon^{i-2}\{\Gamma(f_R^{\varepsilon},f_{i})+\Gamma(f_{i}%
,f_R^{\varepsilon})\},\\
h_1(f)&\equiv \sum_{i+j\geq
n+3}\varepsilon^{i+j-n-2}\Gamma(f_{i},f_{j})
-\varepsilon\{\partial_{t}f_{n+1}+v\cdot\nabla_{x}f_{n+2}\}-\varepsilon
^{2}\partial_{t}f_{n+2},\\
h_1(E_R^{\varepsilon},B_R^{\varepsilon})&\equiv-\frac{\varepsilon^{n-1}}
{\sqrt{\mu}}(E_R^{\varepsilon}+v\times
B_R^{\varepsilon})\cdot\nabla_v
(g_R^{\varepsilon}\sqrt{\mu}) \\
-&\frac{1}{\sqrt{\mu}}\sum_{i=1}^{n+2}\varepsilon^{i-1}\{
(E_i+v\times B_i)\cdot\nabla_v
(g_R^{\varepsilon}\sqrt{\mu})+(E_R^{\varepsilon}+v\times
B_R^{\varepsilon})\cdot\nabla_v
(g_i\sqrt{\mu})\},\\  \label{highh}%
 h_1(E,B)&\equiv-\frac{1} {\sqrt{\mu}}\sum_{i+j\geq
n+2}\varepsilon^{i+j-n-1}(E_i+v\times B_i)\cdot\nabla_v
(g_j\sqrt{\mu}),\\
 h_2(g_R^{\varepsilon})&\equiv\varepsilon^{n-2}
\Gamma(g_R^{\varepsilon},f_R^{\varepsilon})+\sum_{i=1}%
^{n+2}\varepsilon^{i-2}\{\Gamma(g_R^{\varepsilon},f_{i})+\Gamma(g_{i}%
,f_R^{\varepsilon})\},\\
 h_2(g)&\equiv \sum_{i+j\geq
n+3}\varepsilon^{i+j-n-2}\Gamma(g_{i},f_{j})
-\varepsilon\{\partial_{t}g_{n+1}+v\cdot\nabla_{x}g_{n+2}\}-\varepsilon
^{2}\partial_{t}g_{n+2},\\
h_2(E_R^{\varepsilon},B_R^{\varepsilon})&\equiv-\frac{\varepsilon^{n-1}}
{\sqrt{\mu}}(E_R^{\varepsilon}+v\times
B_R^{\varepsilon})\cdot\nabla_v
(f_R^{\varepsilon}\sqrt{\mu}) \\
-&\frac{1}{\sqrt{\mu}}\sum_{i=1}^{n+2}\varepsilon^{i-1}\{
(E_i+v\times B_i)\cdot\nabla_v
(f_R^{\varepsilon}\sqrt{\mu})+(E_R^{\varepsilon}+v\times
B_R^{\varepsilon})\cdot\nabla_v (f_i\sqrt{\mu})\},
\end{split}
\end{equation}
\begin{equation*}
\begin{split}
 h_2(E,B)&\equiv-\frac{1} {\sqrt{\mu}}\sum_{i+j\geq
n+2}\varepsilon^{i+j-n-1}(E_i+v\times B_i)\cdot\nabla_v
(f_j\sqrt{\mu})+\varepsilon E_{n+2}\cdot v\sqrt{\mu}.\\
\end{split}
\end{equation*}\

Our goal is to study $f_{R}^{\varepsilon},$
$g_{R}^{\varepsilon},$ $E_{R}^{\varepsilon}$ and
$B_{R}^{\varepsilon}$ which are equivalent to
$f_{n}^{\varepsilon},$ $g_{n}^{\varepsilon},$
$E_{n}^{\varepsilon}$ and $B_{n}^{\varepsilon}$ via (\ref{ge}).
The procedure closely follows the proofs of the last section. We
first establish pure spatial energy estimate in Lemma
\ref{pure2}. We write the hydrodynamic field parts
$\mathbf{P_1}f_R^{\varepsilon}$ and
$\mathbf{P_2}g_R^{\varepsilon}$ as:
\[
[{\mathbf{P_1}f_R^{\varepsilon}},\;{\mathbf{P_2}g_R^{\varepsilon}}]
=[{\{a_R^{\varepsilon}(t,x)+b_R^{\varepsilon}(t,x)\cdot
v+c_R^{\varepsilon}(t,x)|v|^{2}\}\sqrt{\mu}},\;{d_R^{\varepsilon}(t,x)
\sqrt{\mu}}].
\]\

\begin{lem}
\label{pure2}Assume that (\ref{exp}) are the solution to the
kinetic equation (\ref{rvmb}) such that the following
conservation laws hold:
\begin{equation}
\begin{split}
\int_{\mathbb{T}^3}
a_R^{\varepsilon}dx=&\frac{\varepsilon^n}{2}\int_{\mathbb{T}^3}
|E_R^{\varepsilon}|^2+|B_R^{\varepsilon}|^2dx+
\sum_{i=1}^{n+2}\varepsilon^i\int_{\mathbb{T}^3}
\{E_i\cdot E_R^{\varepsilon}+B_i\cdot B_R^{\varepsilon}\} dx\\
&+\frac{1}{2}\sum_{\substack{i+j\geq n+3
\\1\leq i,j\leq n+2}} \varepsilon^{i+j-n}\int_{\mathbb{T}^3}
E_i\cdot E_j+B_i\cdot B_jdx,\\ \label{hmme}%
\int_{\mathbb{T}^3}
b_R^{\varepsilon}dx=&-\varepsilon^n\int_{\mathbb{T}^3}
E_R^{\varepsilon}\times
B_R^{\varepsilon}dx-\sum_{i=1}^{n+2}\varepsilon^i
\int_{\mathbb{T}^3}
\{E_i\times B_R^{\varepsilon}+E_R^{\varepsilon}\times B_i\} dx\\
&-\sum_{\substack{i+j\geq n+3
\\1\leq i,j\leq n+2}} \varepsilon^{i+j-n}\int_{\mathbb{T}^3}
E_i\times B_jdx,\\
\int_{\mathbb{T}^3}
c_R^{\varepsilon}dx=&-\frac{\varepsilon^n}{6}\int_{\mathbb{T}^3}
|E_R^{\varepsilon}|^2+|B_R^{\varepsilon}|^2dx-\frac{1}{3}
\sum_{i=1}^{n+2}\varepsilon^i\int_{\mathbb{T}^3}
\{E_i\cdot E_R^{\varepsilon}+B_i\cdot B_R^{\varepsilon}\} dx\\
&-\frac{1}{6}\sum_{\substack{i+j\geq n+3
\\1\leq i,j\leq n+2}} \varepsilon^{i+j-n}\int_{\mathbb{T}^3}
E_i\cdot E_j+B_i\cdot B_jdx,\\
\int_{\mathbb{T}^3} d_R^{\varepsilon}dx=&0\;\;\;\text{ and }\;\;\;
\int_{\mathbb{T}^3} B_R^{\varepsilon}dx=0.\\
\end{split}
\end{equation}
Then there is $\lambda>0$ so that for any $\xi>0,$ there exists
some polynomial $U_{\xi }(0)=0$ such that
\begin{equation}
\begin{split}
 \frac{d}{dt}&\{C_{1}\sum_{|\gamma|\leq N}\{||[\partial_{\gamma}
 {f_R^{\varepsilon}},\partial_{\gamma}
{g_R^{\varepsilon}}]||^{2}+
||[\partial_{\gamma}E_R^{\varepsilon},
\partial_{\gamma}B_R^{\varepsilon}]||^2\}-\varepsilon
\delta
G(t)\}\\
+&\delta\sum_{|\gamma|\leq
N}\{\frac{1}{\varepsilon^{2}}||[\partial
_{\gamma}{\mathbf{\{I-P_1\}}f_R^{\varepsilon}},\partial_{\gamma}
{\mathbf{\{I-P_2\}}g_R^{\varepsilon}}]||_{\nu}^{2}+||[\partial_{\gamma
}{\mathbf{P_1}f_R^{\varepsilon}},\partial_{\gamma}
{\mathbf{P_2}g_R^{\varepsilon}}]||^{2}\}\\
\leq&\; e^{-\lambda
t}U_{\xi}(|||\{\mathbf{u}_{n},\mathbf{\theta}_{n},\mathbf{\sigma}_{n}
\}|||_{N}^{2})\{\varepsilon^2+\mathcal{E}_N\}\\
&+C\{(\mathcal{E}_{N})^{\frac{1}{2}}+ \mathcal{E}_{N}+\varepsilon
U_{\xi}(|||\{\mathbf{u}_{n},\mathbf{\theta} _{n},\mathbf{\sigma}
_{n}\}|||_{N})+\xi\} \mathcal{D}_{N},\label{pure2en}
\end{split}
\end{equation}
for $\varepsilon$ sufficiently small,  $\mathcal{E}_N\equiv
\mathcal{E}_N(f_R^{\varepsilon},g_R^{\varepsilon},E_R^{\varepsilon},
B_R^{\varepsilon})(t)$ and $\mathcal{D}_N\equiv
\mathcal{D}_N(f_R^{\varepsilon},g_R^{\varepsilon})(t)$.
\end{lem}\

Notation: $\mathcal{E}_N[i]\equiv
\mathcal{E}_N(f_i,g_i,E_i,B_i)(t)$.

\begin{proof} Note that (\ref{hmme}) falls into the category of
 (\ref{mme}) with $$\mathcal{A}\equiv\sum_{i=1}^{n+2}(||E_i||^2+||B_i||^2).$$
 We apply Lemma \ref{abc} and Lemma \ref{stenergy}
with $j_1^{\varepsilon}=- \varepsilon^3\partial_t E_{n+2},\;
j_2^{\varepsilon}=- \varepsilon^3\partial_t B_{n+2}$ to obtain the
following, an analogue of (\ref{ste}):
\begin{equation}
\begin{split}
 \frac{d}{dt}&\{C_{1}\sum_{|\gamma|\leq N}\{||[\partial_{\gamma}
 {f_R^{\varepsilon}},\partial_{\gamma}
{g_R^{\varepsilon}}]||^{2}+ ||[\partial_{\gamma}E_R^{\varepsilon},
\partial_{\gamma}B_R^{\varepsilon}]||^2\}-\varepsilon
\delta
G(t)\}\\
&+\delta\sum_{|\gamma|\leq
N}\{\frac{1}{\varepsilon^{2}}||[\partial
_{\gamma}{\mathbf{\{I-P_1\}}f_R^{\varepsilon}},\partial_{\gamma}
{\mathbf{\{I-P_2\}}g_R^{\varepsilon}}]||_{\nu}^{2}+||[\partial_{\gamma
}{\mathbf{P_1}f_R^{\varepsilon}},\partial_{\gamma}
{\mathbf{P_2}g_R^{\varepsilon}}]||^{2}\}
\label{ste1}\\
\leq&\; 2C_{1}\sum_{|\gamma|\leq
N}\{\mathbf{(}\partial_{\gamma}h_1^{\varepsilon},\partial_{\gamma
}f_R^{\varepsilon})+\mathbf{(}\partial_{\gamma}h_2^{\varepsilon},\partial_{\gamma
}g_R^{\varepsilon})\}+\varepsilon^{2}\delta\sum_{|\gamma|\leq
N-1}||\partial_{\gamma}h_{||}^{\varepsilon
}||^{2}\\
&+2C_1{\varepsilon^2} (||\partial_{\gamma}E_R^{\varepsilon}||\cdot
||\partial_t\partial_{\gamma}
E_{n+2}||+||\partial_{\gamma}B_R^{\varepsilon}||\cdot
||\partial_t\partial_{\gamma} B_{n+2}||)\\
&+\varepsilon^2\delta\sum_{i=1}^{n+2}(||[E_i,B_i]||^2)^2
+\varepsilon^8\delta (||[\partial_t{
E_{n+2}},\partial_t{B_{n+2}}]||^2+||[\partial_t\nabla
{E_{n+2}},\partial_t\nabla{B_{n+2}}]||^2)\\
\equiv&(I)+(II)+(III).\\
\end{split}
\end{equation}
Now it suffices to estimate $(I),\;(II)\text{ and }(III)$ in the
above to finish Lemma. Observe that by Theorem \ref{cecoe}
\begin{equation}
(III)\leq C\varepsilon^2 e^{-2\lambda
t}U(|||\{\mathbf{u}_n,\theta_n,\sigma_n\}|||_N^2).\label{III}
\end{equation}
As for $(II)$, we use the Cauchy-Schwartz inequality:
\begin{equation}
\begin{split}
(II)\leq &C\varepsilon^2(\mathcal{E}_N)^{\frac{1}{2}} e^{-\lambda
t}U(|||\{\mathbf{u}_n,\theta_n,\sigma_n\}|||_N)\\
\leq&Ce^{-\lambda t}U(|||\{\mathbf{u}_n,\theta_n,\sigma_n\}|||_N)
\{\varepsilon^4+\mathcal{E}_N\}.\label{II}
\end{split}
\end{equation}
The estimate for $(I)$ can be done similarly as in the pure
Boltzmann case. Electromagnetic related terms do not cause any
technical difficulty. Following the proof of Lemma 8.1 in
\cite{Guo2}, we can obtain
\begin{align}
\;\;\;\;\;\varepsilon^{2}\sum_{|\gamma|\leq
N-1}||\partial_{\gamma}h_{||}^{\varepsilon }||^{2}\leq\;
e^{-2\lambda
t}U(|||\{\mathbf{u}_{n},\mathbf{\theta}_{n},\mathbf{\sigma}_{n}\}|||_N^{2})
\{\varepsilon^2+\mathcal{E}_N\}\nonumber\\
\;\;\;\;\;+C\{\varepsilon^2\mathcal{E}_N
+\varepsilon^2e^{-2\lambda
t}U(|||\{\mathbf{u}_{n},\mathbf{\theta}_{n},\mathbf{\sigma}_{n}\}|||_N^{2})\}
\mathcal{D}_N;\label{hproj}\\
\;\;\;\;\;\sum_{|\gamma|\leq
N}\{(\partial_{\gamma}h_1^{\varepsilon},\partial_{\gamma
}f_R^{\varepsilon})\}+(\partial_{\gamma}h_2^{\varepsilon},\partial_{\gamma
}g_R^{\varepsilon}) \leq\; e^{-\lambda
t}U_{\xi}(|||\{\mathbf{u}_{n},\mathbf{\theta}_{n},\mathbf{\sigma}_{n}
\}|||_{N}^{2})\{\varepsilon^2+\mathcal{E}_N\}\nonumber\\
\;\;\;\;\;+C\{\varepsilon(\mathcal{E}_{N})^{\frac{1}{2}}+
\varepsilon U_{\xi}(|||\{\mathbf{u}_{n},\mathbf{\theta}
_{n},\mathbf{\sigma} _{n}\}|||_{N})+\xi\}
\mathcal{D}_{N}.\nonumber
\end{align}
Hence, as adjusting constants and collecting terms, we deduce our
lemma.
\end{proof}\

In order to prove Theorem \ref{high}, it remains to estimate the
velocity derivatives.

\begin{proof}
(\textit{of Theorem \ref{high}:}) It suffices to estimate just the
remaining part
\[[\partial_{\gamma}^{\beta}
{\mathbf{\{I-P_1\}}f_R^{\varepsilon}},
\partial_{\gamma}^{\beta}{\mathbf{\{I-P_2\}}g_R^{\varepsilon}}].\]
We follow exactly the same argument in the proof of Theorem
\ref{nslimit}. Comparing $\partial_{\gamma_1}^{\beta}$ derivative
of (\ref{highfg}) with (\ref{vderivative}), we notice the only
difference is that now
$\partial_{\gamma_1}^{\beta}h_1^{\varepsilon}$ and
$\partial_{\gamma_1}^{\beta}h_2^{\varepsilon}$ are more
complicated. Thus it is enough to estimate only
\[
(\partial_{\gamma_1}^{\beta}h_1^{\varepsilon},
\partial_{\gamma_1}^{\beta}\mathbf{\{I-P_1\}}f_R^{\varepsilon})
+(\partial_{\gamma_1}^{\beta}h_2^{\varepsilon},
\partial_{\gamma_1}^{\beta}\mathbf{\{I-P_2\}}
g_R^{\varepsilon}).
\]
Here we illustrate how it works for electromagnetic related
terms. The terms containing
$h_1(f),h_1(f_R^{\varepsilon}),h_2(g)\text{ or }
h_2(g_R^{\varepsilon})$ can be done in the same way as in
\cite{Guo2}. See the proof of Theorem 2.3 there (p.50-51).
Recalling (\ref{rdecay}), we first get
\begin{align*}
(\partial_{\gamma_1}^{\beta}&h_1(E,B),
\partial_{\gamma_1}^{\beta}\mathbf{\{I-P_1\}}f_R^{\varepsilon})
+(\partial_{\gamma_1}^{\beta}h_2(E,B),
\partial_{\gamma_1}^{\beta}\mathbf{\{I-P_2\}}
g_R^{\varepsilon})\\
&  \leq\varepsilon e^{-\lambda
t}U(|||\{\mathbf{u}_{n},\mathbf{\theta }_{n},\mathbf{\sigma
}_{n}\}|||_{N})||[\partial_{\gamma_1}^{\beta}
{\mathbf{\{I-P_1\}}f_R^{\varepsilon}},\partial_{\gamma_1}^{\beta}
{\mathbf{\{I-P_2\}}g_R^{\varepsilon}}]||_{\nu}\\
&  \leq \varepsilon^2 e^{-\lambda
t}U(|||\{\mathbf{u}_{n},\mathbf{\theta }_{n},\mathbf{\sigma
}_{n}\}|||_{N})(\mathcal{D}_{N})^{\frac{1}{2}}\\
&\leq\varepsilon^2 e^{-2\lambda
t}U_{\xi}(|||\{\mathbf{u}_{n},\mathbf{\theta }_{n},\mathbf{\sigma
}_{n}\}|||_{N}^2)+\frac{\xi}{2}\mathcal{D}_{N}.
\end{align*}
Note that we have also taken  into account $\varepsilon
\partial_{\gamma_1}^{\beta}(E_{n+2}\cdot v\sqrt{\mu})$ in $h_2(E,B)$. On the
other hand, as for $h_1(E_R^{\varepsilon},B_R^{\varepsilon})\text{
and } h_2(E_R^{\varepsilon},B_R^{\varepsilon})$, by splitting
$f_R^{\varepsilon},g_R^{\varepsilon}$ into
$\mathbf{P_1}f_R^{\varepsilon}+(\mathbf{I-P_1})f_R^{\varepsilon}
,\;\mathbf{P_2}g_R^{\varepsilon}+(\mathbf{I-P_2})g_R^{\varepsilon}$
 as in (\ref{ebsplit}), one
can get
\begin{align*}
(\partial&_{\gamma_1}^{\beta}h_1(E_R^{\varepsilon},B_R^{\varepsilon}),
\partial_{\gamma_1}^{\beta}\mathbf{\{I-P_1\}}f_R^{\varepsilon})
+(\partial_{\gamma_1}^{\beta}h_2(E_R^{\varepsilon},B_R^{\varepsilon}),
\partial_{\gamma_1}^{\beta}\mathbf{\{I-P_2\}}
g_R^{\varepsilon})\\
&\leq\varepsilon(\mathcal{E}_{N})^{\frac{1}{2}}||[\partial_{\gamma_1}
^{\beta}{g_R^{\varepsilon}},\partial_{\gamma_1}
^{\beta}{f_R^{\varepsilon}}]||_{\nu}\cdot||[\partial_{\gamma_1}^{\beta}
{\mathbf{\{I-P_1\}}f_R^{\varepsilon}},\partial_{\gamma_1}^{\beta}
{\mathbf{\{I-P_2\}}g_R^{\varepsilon}}]||_{\nu}\\
&\;\;\;+\varepsilon
\sum_{i=1}^{n+2}\varepsilon^{i-1}\{(\mathcal{E}_{N}[i])
^{\frac{1}{2}} ||[\partial_{\gamma_1}
^{\beta}{g_R^{\varepsilon}},\partial_{\gamma_1}
^{\beta}{f_R^{\varepsilon}}]||_{\nu}+(\mathcal{E}_{N})^{\frac{1}{2}}
||[\partial_{\gamma_1}^{\beta}\nabla_v{g_{i}},
\partial_{\gamma_1}^{\beta}\nabla_v{f_i}]||_{\nu}\}\\
&\;\;\;\;\;\;\;\;\;\;\;\;\;\;\;\;\;\;\;\;\;\cdot
\frac{1}{\varepsilon}||[\partial_{\gamma_1}^{\beta}
{\mathbf{\{I-P_1\}}f_R^{\varepsilon}},\partial_{\gamma_1}^{\beta}
{\mathbf{\{I-P_2\}}g_R^{\varepsilon}}]||_{\nu}\\
& \leq\varepsilon(\mathcal{E}_{N})^{\frac{1}{2}}\mathcal{D}
_{N}+\varepsilon e^{-\lambda
t}U(|||\{\mathbf{u}_{n},\mathbf{\theta} _{n},\mathbf{\sigma}
_{n}\}|||_{N})\{\mathcal{D}_{N}+(\mathcal{E}_{N})^{\frac{1}{2}}
(\mathcal{D}_{N})^{\frac{1}{2}}\}\\
&\leq\varepsilon e^{-\lambda
t}U(|||\{\mathbf{u}_{n},\mathbf{\theta} _{n},\mathbf{\sigma}
_{n}\}|||_{N})\mathcal{E}_{N}+\{\varepsilon(\mathcal{E}_{N})
^{\frac{1}{2}}+e^{-\lambda
t}U(|||\{\mathbf{u}_{n},\mathbf{\theta} _{n},\mathbf{\sigma}
_{n}\}|||_{N})\}\mathcal{D}_{N}.
\end{align*}\

Therefore, as in the proof of Theorem \ref{nslimit}$,$ with all
the above estimates as well as the ones in \cite{Guo2} (p.50-51),
we deduce that for any $\xi>0,$
\begin{equation}
\begin{split}
  \frac{d}{dt}&[
K\{C_{1}\{||[\partial_{\gamma}{f_R^{\varepsilon}},\partial_{\gamma}
{g_R^{\varepsilon}}]||^{2}+ ||[\partial_{\gamma}E_R^{\varepsilon},
\partial_{\gamma}B_R^{\varepsilon}]||^2\}-\varepsilon
\delta G(t)\}\\&+2\delta||[\partial_{\gamma}^{\beta}
{\mathbf{\{I-P_1\}}f_R^{\varepsilon}},\partial_{\gamma}^{\beta}
{\mathbf{\{I-P_2\}}g_R^{\varepsilon}}
||^{2}-4\varepsilon\delta\widetilde{G}(t)]
+\delta\mathcal{D}_{N}\\
 \leq&\; e^{-\lambda
t}U_{\xi}(|||\{\mathbf{u}_{n},\mathbf{\theta}_{n},\mathbf{\sigma}_{n}
\}|||_{N}^{2})\{\varepsilon^2+\mathcal{E}_N\}\\
&+C\{(\mathcal{E}_{N})^{\frac{1}{2}}+ \mathcal{E}_{N}+\varepsilon
U_{\xi}(|||\{\mathbf{u}_{n},\mathbf{\theta} _{n},\mathbf{\sigma}
_{n}\}|||_{N})+\xi\} \mathcal{D}_{N}.\label{longenergy}
\end{split}
\end{equation}
We redefine an equivalent instant energy functional $\mathcal{E}
_{N}$ as in (\ref{hardin}) to get
\begin{align*}
 \frac{d}{dt}\mathcal{E}_{N}+\mathcal{D}_{N}
\leq&\; e^{-\lambda
t}U_{\xi}(|||\{\mathbf{u}_{n},\mathbf{\theta}_{n},\mathbf{\sigma}_{n}
\}|||_{N}^{2})\{\varepsilon^2+\mathcal{E}_N\}\\
&+C\{(\mathcal{E}_{N})^{\frac{1}{2}}+ \mathcal{E}_{N}+\varepsilon
U_{\xi}(|||\{\mathbf{u}_{n},\mathbf{\theta} _{n},\mathbf{\sigma}
_{n}\}|||_{N})+\xi\} \mathcal{D}_{N}.
\end{align*}
Choose and then fix $\varepsilon$ and $\xi$ so that $\varepsilon
U_{\xi}(|||\{\mathbf{u}_{n},\mathbf{\theta} _{n},\mathbf{\sigma}
_{n}\}|||_{N})$ is sufficiently small. We also assume that
\[
\mathcal{E}_{N}\leq M
\]
sufficiently small such that the coefficient in front of
$\mathcal{D} _{N}$ satisfies
\[
C\{(\mathcal{E}_{N})^{\frac{1}{2}}+ \mathcal{E}_{N}+\varepsilon
U_{\xi}(|||\{\mathbf{u}_{n},\mathbf{\theta} _{n},\mathbf{\sigma}
_{n}\}|||_{N})+\xi\}<\frac{1}{2}.
\]
Thus we  obtain%
\begin{equation}
\frac{d}{dt}\mathcal{E}_{N}+\frac{1}{2}\mathcal{D}%
_{N}\leq Ce^{-\lambda
t}U(|||\{\mathbf{u}_{n},\mathbf{\theta}_{n},\mathbf{\sigma}_{n}
\}|||_{N}^{2})\{\varepsilon^{2}+\mathcal{E}_{N}%
\}. \label{highen}%
\end{equation}
In turn, we have
\begin{align*}
&\frac{d}{dt}\{e^{-CU(|||\{\mathbf{u}_{n},\mathbf{\theta}_{n},
\mathbf{\sigma}_{n}
\}|||_{N}^{2})\int_{0}^{t}e^{-\lambda s}ds}\mathcal{E}_{N}\}\\
&  \leq
Ce^{-CU(|||\{\mathbf{u}_{n},\mathbf{\theta}_{n},\mathbf{\sigma}_{n}
\}|||_{N}^{2})\int_{0}%
^{t}e^{-\lambda s}ds}e^{-\lambda
t}\varepsilon^{2}U(|||\{\mathbf{u}_{n},\mathbf{\theta}_{n},\mathbf{\sigma}_{n}
\}|||_{N}^{2}).
\end{align*}
Integrating over $t$, we deduce
\[
\sup_{0\leq t\leq\infty}\mathcal{E}_{N}(t)\leq
e^{CU(|||\{\mathbf{u}_{n},\mathbf{\theta}_{n},\mathbf{\sigma}_{n}
\}|||_{N}^{2})\int_{0}^{\infty }e^{-\lambda
s}ds}\{\mathcal{E}_{N}(0)+C\varepsilon
^{2}U(|||\{\mathbf{u}_{n},\mathbf{\theta}_{n},\mathbf{\sigma}_{n}
\}|||_{N}^{2})\}.
\]
We conclude for $\varepsilon$ sufficiently small and for some
other polynomial $U,$
\begin{equation}
\begin{split}
\sup_{0\leq t\leq\infty}\mathcal{E}_{N}(t)&\leq
e^{U(|||\{\mathbf{u}_{n},\mathbf{\theta}_{n},\mathbf{\sigma}_{n}
\}|||_{N}^{2})}\{\mathcal{E}%
_{N}(0)+\varepsilon^{2}
U(|||\{\mathbf{u}_{n},\mathbf{\theta}_{n},\mathbf{\sigma}_{n}
\}|||_{N})\}\\
&<M/2. \label{gebound}
\end{split}
\end{equation}
A standard continuity argument shows that the hypothesis $\mathcal{E}%
_{N}(t)\leq M$ is valid and (\ref{gebound}) is proven.

Recalling
\[
\begin{split}
f_R^{\varepsilon}=\{f_{n}^{\varepsilon}-f_{n}\}+\varepsilon
f_{n+1}+\varepsilon^{2}f_{n+2},\;\;E_R^{\varepsilon}=\{E_{n}^{\varepsilon}-E_{n}\}+\varepsilon
E_{n+1}+\varepsilon^{2}E_{n+2},\\
g_R^{\varepsilon}=\{g_{n}^{\varepsilon} -g_{n}\}+\varepsilon
g_{n+1}+\varepsilon^{2}g_{n+2},\;\;B_R^{\varepsilon}=\{B_{n}^{\varepsilon}-B_{n}\}+\varepsilon
B_{n+1}+\varepsilon^{2}B_{n+2},
\end{split}
\]
and by Theorem \ref{cecoe},
\begin{equation}
\begin{split}
\mathcal{E}_{N}&\mathcal{(}\varepsilon
f_{n+1}+\varepsilon^{2}f_{n+2},\varepsilon
g_{n+1}+\varepsilon^{2}g_{n+2},\varepsilon
E_{n+1}+\varepsilon^{2}E_{n+2},\varepsilon
B_{n+1}+\varepsilon^{2}B_{n+2})\\
&\leq C\varepsilon^{2}e^{-2\lambda
t}U(|||\{\mathbf{u}_{n},\mathbf{\theta}_{n},\mathbf{\sigma}_{n}\}
|||_{N}^{2}),\label{leftoverdecay}
\end{split}
\end{equation}
we thus deduce (\ref{hen}) for
$f_{n}^{\varepsilon}-f_{n},\;g_{n}^{\varepsilon}-g_{n},\;
E_{n}^{\varepsilon}-E_{n},\;B_{n}^{\varepsilon}-B_{n}$.

To get a decay rate (\ref{hdecayrate}), we use the interpolation
argument as done in the last part of the previous section. Let
$|\gamma|\leq N$ and $k\geq 1$. Recall that
\begin{equation*}
||\partial_{\gamma}E_R^{\varepsilon}||^2+||\partial_{\gamma}B_R^{\varepsilon}||^2
\leq C
\{||E_R^{\varepsilon}||_{H^{N-1}}^{\frac{2k}{k+1}}+||B_R^{\varepsilon}||
_{H^{N-1}}^{\frac{2k}{k+1}}\}\{||E_R^{\varepsilon}||_{H^{N+k}}
^{\frac{2}{k+1}}+||B_R^{\varepsilon}||
_{H^{N+k}}^{\frac{2}{k+1}}\}.
\end{equation*}
Denote a bound for $\mathcal{E}_N$ by $I_N$, i.e.
\[
I_N\equiv
e^{U(|||\{\mathbf{u}_{n},\mathbf{\theta}_{n},\mathbf{\sigma}_{n}
\}|||_{N}^{2})}\{\mathcal{E}%
_{N}(0)+\varepsilon^{2}
U(|||\{\mathbf{u}_{n},\mathbf{\theta}_{n},\mathbf{\sigma}_{n}
\}|||_{N})\}.
\]
Due to  (\ref{b1}), (\ref{e11}), (\ref{hproj}) and
(\ref{gebound}) we have
\begin{equation*}
\begin{split}
||E_R^{\varepsilon}||_{H^{N}}^{\frac{2k+2}{k}}+||B_R^{\varepsilon}||
_{H^{N}}^{\frac{2k+2}{k}}&\leq
C_{N+k}(I_{N+k})^{\frac{1}{k}}\{||E_R^{\varepsilon}||
_{H^{N-1}}^{2}+||B_R^{\varepsilon}|| _{H^{N-1}}^{2}\}\\
&\leq C_{N+k}(I_{N+k})^{\frac{1}{k}}\{\varepsilon\frac{d}{dt}
\widetilde{\widetilde{G}}+\mathcal{D}_N+e^{-2\lambda
t}U(|||\{\mathbf{u}_{n},\mathbf{\theta}_{n},\mathbf{\sigma}_{n}
\}|||_{N}^2)I_N\},
\end{split}
\end{equation*}
where $\varepsilon$ sufficiently small and
$\widetilde{\widetilde{G}}(t)$ is defined in (\ref{ggbound}) so
that $|\widetilde{\widetilde{G}}(t)|\leq \mathcal{E}_N$. Noting
that the other part of $\mathcal{E}_N$, i.e. the
nonelectromagnetic part is bounded by $\mathcal{D}_N$, a lower
bound for $\mathcal{D}_N$ can be given as following:
\[
C_{N,k}(I_{N+k})^{-\frac{1}{k}}\mathcal{E}_N^{\frac{k+1}{k}}-
\varepsilon\frac{d}{dt} \widetilde{\widetilde{G}}-e^{-2\lambda
t}U(|||\{\mathbf{u}_{n},\mathbf{\theta}_{n},\mathbf{\sigma}_{n}
\}|||_{N}^2)I_N\leq \mathcal{D}_N.
\]
It follows from (\ref{highen}) that
\begin{equation}
\frac{d}{dt}\{\mathcal{E}_N-
\varepsilon\widetilde{\widetilde{G}}\}+C_{N,k}(I_{N+k})^{-\frac{1}{k}}
\mathcal{E}_N^{\frac{k+1}{k}}\leq C e^{-\lambda t}
I_{N}\;U(|||\{\mathbf{u}_{n},\mathbf{\theta}_{n},\mathbf{\sigma}_{n}
\}|||_{N}^{2}).\label{hegd}
\end{equation}
Letting $\widetilde{\mathcal{E}}_N=\mathcal{E}_N-
\varepsilon\widetilde{\widetilde{G}},$ we get
$\frac{1}{C}\mathcal{E}_N\leq \widetilde{\mathcal{E}}_N\leq C
\mathcal{E}_N$ for some $C>1$ and thus (\ref{hegd}) becomes
\[
\frac{d}{dt}\widetilde{\mathcal{E}}_N+
C_{N,k}(I_{N+k})^{-\frac{1}{k}}
\widetilde{\mathcal{E}}_N^{\frac{k+1}{k}}\leq C e^{-\lambda t}
I_{N}\;U(|||\{\mathbf{u}_{n},\mathbf{\theta}_{n},\mathbf{\sigma}_{n}
\}|||_{N}^{2}).
\]
Multiply the above by $(1+\frac{t}{k})^k$, we get
\begin{equation}
\frac{d}{dt}\{(1+\frac{t}{k})^k\widetilde{\mathcal{E}}_N(t)\}\leq
Q(t)+C  e^{-\lambda t}(1+\frac{t}{k})^k
I_{N}\;U(|||\{\mathbf{u}_{n},\mathbf{\theta}_{n},\mathbf{\sigma}_{n}
\}|||_{N}^{2}),\label{poly}
\end{equation}
where
\[
\begin{split}
Q(t)&\equiv(1+\frac{t}{k})^{k-1}\widetilde{\mathcal{E}}_N(t)-C_{N,k}
(I_{N+k})^{-\frac{1}{k}}(1+\frac{t}{k})^k
\widetilde{\mathcal{E}}_N^{\frac{k+1}{k}}(t)\\
&=(1+\frac{t}{k})^{k-1}\widetilde{\mathcal{E}}_N(t)\{1-C_{N,k}
(I_{N+k})^{-\frac{1}{k}}\{(1+\frac{t}{k})^k
\widetilde{\mathcal{E}}_N(t)\}^{\frac{1}{k}}\}.
\end{split}
\]\
To conclude our theorem, it suffices to verify the following
statement:

Claim. There exists $\widetilde{C}_{N+k}>0$ such that
\[
\sup_{t}\{(1+\frac{t}{k})^k \widetilde{\mathcal{E}}_N(t)\}\leq
\widetilde{C}_{N+k}I_{N+k},
\]
since the same conclusion is valid for $\mathcal{E}_N$ and
recalling (\ref{ge}), combining with (\ref{leftoverdecay}), we can
deduce (\ref{hdecayrate}).

Proof of Claim: Let $S$ be the set of $t$ such that $Q(t)<0$. Note
that $Q(t)>0$ for sufficiently small $M$ and small $t$, which
implies that $S^c$ is nonempty. If $S$ is an empty set, namely
$Q(t)\geq 0$ for all $t$, we can set
$\widetilde{C}_{N+k}=(C_{N,k})^{-k}$. Let $t_1\in S$. We can find
$t_0\leq t_1$ so that $t_0\in S^c$ and $Q(t)\leq 0$ for $t_0\leq t
\leq t_1$. Integrate (\ref{poly}) from $t_0$ to $t_1$ to get
\[
\begin{split}
(1+\frac{t_1}{k})^k \widetilde{\mathcal{E}}_N(t_1)&\leq
(1+\frac{t_0}{k})^k \widetilde{\mathcal{E}}_N(t_0)+
CI_{N}U(|||\{\mathbf{u}_{n},\mathbf{\theta}_{n},\mathbf{\sigma}_{n}
\}|||_{N}^{2})\int_{t_0}^{t_1}  e^{-\lambda s}(1+\frac{s}{k})^k
ds\\
&\leq (C_{N,k})^{-k}I_{N+k}
+CI_{N}U(|||\{\mathbf{u}_{n},\mathbf{\theta}_{n},\mathbf{\sigma}_{n}
\}|||_{N}^{2})\int_{0}^{\infty}  e^{-\lambda s}(1+\frac{s}{k})^k
ds.
\end{split}
\]
Choosing
$\widetilde{C}_{N+k}=(C_{N,k})^{-k}+CU(|||\{\mathbf{u}_{n},
\mathbf{\theta}_{n},\mathbf{\sigma}_{n}
\}|||_{N}^{2})\int_{0}^{\infty}  e^{-\lambda s}(1+\frac{s}{k})^k
ds$, we conclude the proof.
\end{proof}\

\textbf{Acknowledgments:} The author would like to deeply thank
\textsc{Yan Guo}
for many stimulating discussions.\\

\end{document}